\pdfoutput=1
\RequirePackage{ifpdf}
\ifpdf 
\documentclass[pdftex]{sigma}
\else
\documentclass{sigma}
\fi

\usepackage{tabularx}
\usepackage{booktabs}
\usepackage{array}
\usepackage{multirow}
\newcolumntype{?}{!{\vrule width 1pt}}
\usepackage{makecell}

\usepackage{stmaryrd, amsbsy,
 mathrsfs, wasysym, cancel,textgreek}

 \usepackage{tikz, tikz-cd}
 \usetikzlibrary{arrows, matrix,angles,patterns}
 \usetikzlibrary{decorations.pathreplacing,decorations.markings,decorations.text}
\tikzstyle{n}=[circle, draw, fill, minimum size=6, inner sep=0]
\tikzstyle{ngr}=[circle, draw, fill, gray, minimum size=6, inner sep=0]
\tikzstyle{int}=[draw, circle, fill, minimum size=4, inner sep=1]
\tikzstyle{root}=[circle, draw, fill, minimum size=0, inner sep=0]
\tikzstyle{vertex}=[circle, draw, fill, minimum size=0, inner sep=0.7]
\tikzstyle{leaf}=[circle, draw, fill, minimum size=0, inner sep=0]
\tikzstyle{intgr}=[draw, circle, fill, gray, minimum size=3, inner sep=1]
\tikzstyle{ext}=[draw, circle, minimum size=1, inner sep=1]
\tikzstyle{exttiny}=[draw, circle, minimum size=1, inner sep=0.65]
\tikzset{lab/.style={draw, circle, minimum size=5, inner sep=1.1},
n/.style={draw, circle, fill, minimum size=5, inner sep=1}}
\tikzstyle{fleche}=[->,>= latex?,thick]
\tikzset{
 rightgreen/.style={
 decoration={markings,mark=at position .45 with {\arrow[scale=1.2,caribbeangreen]{latex}}},
 postaction={decorate},
 shorten >=0.4pt}}
\tikzset{
 leftgreen/.style={
 decoration={markings,mark=at position .2 with {\arrowreversed[scale=1.2,caribbeangreen]{latex}}},
 postaction={decorate},
 shorten >=0.4pt}}

 \tikzset{
 carrotshifted/.style={
 decoration={markings,mark=at position .7 with {\arrow[scale=1.2,carrotorange]{latex}}},
 postaction={decorate},
 shorten >=0.4pt}}

 \tikzset{
 rightredshifted/.style={
 decoration={markings,mark=at position .2 with {\arrow[scale=1.2,linkcolor]{latex}}},
 postaction={decorate},
 shorten >=0.4pt}}

 \tikzset{
 rightred/.style={
 decoration={markings,mark=at position .55 with {\arrow[scale=1.2,linkcolor]{latex}}},
 postaction={decorate},
 shorten >=0.4pt}}
\tikzset{
 leftred/.style={
 decoration={markings,mark=at position .25 with {\arrowreversed[scale=1.2,linkcolor]{latex}}},
 postaction={decorate},
 shorten >=0.4pt}}

\definecolor{Green}{RGB}{147,162,153}
\definecolor{Green2}{RGB}{26,148,49}
\definecolor{BrownL}{RGB}{173,143,103}
\definecolor{Red}{RGB}{210,83,60}
\definecolor{BrownD}{RGB}{114,96,86}
\definecolor{GreyD}{RGB}{76,90,106}
\definecolor{GreyB}{RGB}{128,141,160}
\definecolor{Maroon}{RGB}{121,70,61}
\definecolor{Blue}{RGB}{148,184,210}
\definecolor{Blue2}{RGB}{108,144,170}
\definecolor{Blue3}{RGB}{42, 107, 172}
\definecolor{linkcolor}{RGB}{150,50,60}
\definecolor{MedianBrown}{RGB}{119,95,85}
\definecolor{MedianLightBrown}{RGB}{235,221,195}
\definecolor{MedianLightBlue}{RGB}{148,184,210}
\definecolor{MedianOrange}{RGB}{230,140,71}
\definecolor{MedianLightOrange}{RGB}{216,178,92}
\definecolor{MedianUltraLightOrange}{RGB}{250,246,235}
\definecolor{MedianUltraLightBlue}{RGB}{230,230,235}
\definecolor{amber}{rgb}{1.0, 0.75, 0.0}
\definecolor{caribbeangreen}{rgb}{0.0, 0.8, 0.6}
\definecolor{chocolate}{rgb}{0.48, 0.25, 0.0}
\definecolor{carrotorange}{rgb}{0.99, 0.57, 0.13}
\definecolor{cerulean}{rgb}{0.0, 0.48, 0.65}
\definecolor{colororange}{HTML}{E65100} 
\definecolor{colordgray}{HTML}{795548} 
\definecolor{colorhgray}{HTML}{212121} 
\definecolor{colorgreen}{HTML}{009688} 
\definecolor{colorlgray}{HTML}{FAFAFA} 
\definecolor{colorblue}{HTML}{0277BB} 

\numberwithin{equation}{section}

\newtheorem{Theorem}{Theorem}[section]
\newtheorem*{Theorem*}{Theorem}
\newtheorem{Corollary}[Theorem]{Corollary}
\newtheorem{Lemma}[Theorem]{Lemma}
\newtheorem{Proposition}[Theorem]{Proposition}
\newtheorem{Conjecture}[Theorem]{Conjecture}
 { \theoremstyle{definition}

\newtheorem{Example}[Theorem]{Example}
\newtheorem{Remark}[Theorem]{Remark}
}

\newcommand{\grt}{\mathfrak{grt}}
\newcommand{\GRT}{\mathsf{GRT}}
\newcommand{\mLie}{\mathsf{Lie}}
\newcommand{\Lieinf}{$\mathsf{Lie}_\infty$}
\newcommand{\End}{\mathsf{End}}
\newcommand{\sgn}{{\text{\rm sgn}}}
\newcommand{\Id}{\mathrm{id}}
\newcommand{\Or}{{\mathsf{O}\vec{\mathsf r}}}
\newcommand{\Def}{\text{\rm Def}}
\newcommand{\CE}{{\rm\mathsf{CE}}}
\newcommand{\NR}{{\mathsf {NR}}}
\newcommand{\Rep}{{\rm\mathsf{Rep}}}
\newcommand{\mGer}{\mathsf{Ger}}

\newcommand{\ad}{\text{\rm ad}}
\newcommand\Hom{\text{\rm Hom}}
\newcommand\Ob{\text{\rm Ob}}

\newcommand\fonc[1]{\mathscr C^\infty\pl#1\pr}
\newcommand{\foncm}{{\mathscr C}^\infty(\M)}
\newcommand\form[1]{\Omega^1\pl #1\pr}

\newcommand{\vf}{\Gamma(T\M)}
\newcommand{\ff}{\form{\M}}
\newcommand{\alg}{\mathfrak g}
\newcommand{\w}{\wedge}
\newcommand{\foncg}{\mathscr C^\infty{((\alg\oplus\alg^*)[1])}}
\newcommand{\p}{\partial}
\newcommand{\half}{\frac{1}{2}}
\newcommand{\om}{\omega}
\newcommand{\Om}{\Omega}
\newcommand*{\longhookrightarrow}{\ensuremath{\lhook\joinrel\relbar\joinrel\rightarrow}}
\newcommand{\iso}{\overset{\sim}{\longrightarrow}}
\newcommand{\deltaS}{\delta_{\rm\mathsf{S}}}
\newcommand{\deltabr}{\delta}

\newcommand{\A}{\mathcal{A}}
\newcommand{\mE}{\mathcal{E}}
\newcommand{\mN}{\mathcal{N}}
\newcommand{\U}{\mathcal{U}}

\newcommand{\M}{\mathscr{M}}

\newcommand{\corps}{{\mathbb K}}
\newcommand{\mS}{{\mathbb S}}

\newcommand{\NPQ}{{\rm\mathsf{NPQ}}}
\newcommand{\NQ}{{\rm\mathsf{NQ}}}
\newcommand{\Q}{{\mathsf Q}}
\newcommand{\cH}{{\mathscr H}}
\newcommand{\Lag}{\mathcal{L}}
\newcommand{\symp}{\Omega}
\newcommand{\br}[2]{\crl #1,#2 \crr}
\newcommand{\pb}[2]{\big\{{#1},{#2}\big\}}
\newcommand{\brdot}{\br{\cdot}{\cdot}}
\newcommand{\pbdot}{\left\{{\cdot},{\cdot}\right\}}
\newcommand{\TE}{{\mathscr E}}
\newcommand{\fE}{\Gamma(E)}
\newcommand{\fmE}{\Gamma(\mE)}
\newcommand{\foncTE}{\mathscr C^\infty{(\TE)}}

\newcommand{\pl}{\left(}
\newcommand{\pr}{\right)}
\newcommand{\crl}{\left[}
\newcommand{\crr}{\right]}
\newcommand{\dl}{\left \langle}
\newcommand{\dr}{\right \rangle}
\newcommand{\dlr}[2]{\dl#1,#2\dr}
\newcommand{\dldot}{{\dl\cdot,\cdot\dr}}

\newcommand{\bi}{\begin{itemize}}
\newcommand{\ei}{\end{itemize}}
\newcommand{\be}{\begin{enumerate}}
\newcommand{\ee}{\end{enumerate}}
\newcommand\pset[1]{\left \lbrace #1\right \rbrace}

\newcommand\Prop[1]{Proposition \ref{#1}}

\newcommand{\cC}{\mathcal{C}}
\newcommand{\dD}{\mathcal{D}}

\newcommand{\Tpoly}{{\mathcal T_\text{\rm poly}}}
\newcommand{\Dpoly}{{\mathcal D_\text{\rm poly}}}
\newcommand{\eps}{{\vcenter{\hbox{{$[[\hbar]]$}}}}}

\newcommand{\dime}{{d}}
\newcommand{\graphGam}{{\gamma}}
\newcommand{\Gamgraph}{{\gamma}}
\newcommand{\fcGC}{\mathsf{fcGC}}
\newcommand{\odfGC}[3]{\mathsf{o}_{#3}\mathsf{d}_{#2}\mathsf{fGC}_{#1}}
\newcommand{\cGra}{\mathsf{cGra}}
\newcommand{\ugra}{\mathsf{gra}}
\newcommand{\dGraplain}{\mathsf{dGra}}
\newcommand{\Graplain}{\mathsf{Gra}}
\newcommand{\odgraplain}{\mathsf{ogra}}
\newcommand{\dgra}[2]{\mathsf{d}_{#2}\mathsf{gra}_{#1}}
\newcommand{\dcGra}[2]{\mathsf{d}_{#2}\mathsf{cGra}_{#1}}
\newcommand{\odGra}[3]{\mathsf{o}_{#3}\mathsf{d}_{#2}\mathsf{Gra}_{#1}}
\newcommand{\odcGra}[3]{\mathsf{o}_{#3}\mathsf{d}_{#2}\mathsf{cGra}_{#1}}
\newcommand{\odgra}[2]{\mathsf{o}_{#2}\mathsf{d}_{#1}\mathsf{gra}}
\newcommand{\sodGra}[3]{\mathsf{s}_{#3}\mathsf{d}_{#2}\mathsf{Gra}_{#1}}
\newcommand{\sodcGra}[3]{\mathsf{s}_{#3}\mathsf{d}_{#2}\mathsf{cGra}_{#1}}
\newcommand{\odfcGC}[3]{\mathsf{o}_{#3}\mathsf{d}_{#2}\mathsf{fcGC}_{#1}}
\newcommand{\sodfGC}[3]{\mathsf{s}_{#3}\mathsf{d}_{#2}\mathsf{fGC}_{#1}}
\newcommand{\sodfcGC}[3]{\mathsf{s}_{#3}\mathsf{d}_{#2}\mathsf{fcGC}_{#1}}
\newcommand{\sodGrablack}[3]{\mathsf{s}_{#3}\mathsf{d}_{#2}\mathsf{Gra}_{#1}^{\text{\rm black}}}
\newcommand{\sodGrared}[3]{\mathsf{s}_{#3}\mathsf{d}_{#2}\mathsf{Gra}_{#1}^{\text{\rm red}}}
\newcommand{\odGrablack}[3]{\mathsf{o}_{#3}\mathsf{d}_{#2}\mathsf{Gra}_{#1}^{\text{\rm black}}}
\newcommand{\odGrared}[3]{\mathsf{o}_{#3}\mathsf{d}_{#2}\mathsf{Gra}_{#1}^{\text{\rm red}}}
\newcommand{\odfGCblack}[3]{\mathsf{o}_{#3}\mathsf{d}_{#2}\mathsf{fGC}_{#1}^{\text{\rm black}}}
\newcommand{\odfGCred}[3]{\mathsf{o}_{#3}\mathsf{d}_{#2}\mathsf{fGC}_{#1}^{\text{\rm red}}}
\newcommand{\dfcGC}[2]{\mathsf{d}_{#2}\mathsf{fcGC}_{#1}}
\newcommand{\dfGC}[2]{\mathsf{d}_{#2}\mathsf{fGC}_{#1}}
\newcommand{\dGra}[2]{\mathsf{d}_{#2}\mathsf{Gra}_{#1}}
\newcommand{\gra}{\mathsf{dgra}}
\newcommand{\fGC}{\mathsf{fGC}}
\newcommand{\dfcGCplain}{\mathsf{dfcGC}}
\newcommand{\GC}{\mathsf{GC}}

\begin{document}
\allowdisplaybreaks

\newcommand{\arXivNumber}{2102.07593}

\renewcommand{\PaperNumber}{020}

\FirstPageHeading

\ShortArticleName{A Note on Multi-Oriented Graph Complexes}

\ArticleName{A Note on Multi-Oriented Graph Complexes \\ and Deformation Quantization of Lie Bialgebroids}

\Author{Kevin MORAND~$^{\rm ab}$}

\AuthorNameForHeading{K.~Morand}

\Address{$^{\rm a)}$~Department of Physics, Sogang University, Seoul 04107, South Korea}
\EmailD{\href{mailto:morand@sogang.ac.kr}{morand@sogang.ac.kr}}
\Address{$^{\rm b)}$~Center for Quantum Spacetime, Sogang University, Seoul 04107, South Korea}

\ArticleDates{Received July 07, 2021, in final form March 09, 2022; Published online March 20, 2022}

\Abstract{Universal solutions to deformation quantization problems can be conveniently classified by the cohomology of suitable graph complexes. In particular, the deformation quantizations of (finite-dimensional) Poisson manifolds and Lie bialgebras are characterised by an action of the Grothendieck--Teichm\"uller group via one-colored directed and oriented graphs, respectively. In this note, we study the action of multi-oriented graph complexes on Lie bialgebroids and their ``quasi'' generalisations. Using results due to T.~Willwacher and M.~\v{Z}ivkovi\'c on the cohomology of (multi)-oriented graphs, we show that the action of the Grothendieck--Teichm\"uller group on Lie bialgebras and quasi-Lie bialgebras can be generalised to quasi-Lie bialgebroids via graphs with two colors, one of them being oriented. However, this action generically fails to preserve the subspace of Lie bialgebroids. By resorting to graphs with two oriented colors, we instead show the existence of an obstruction to the quantization of a generic Lie bialgebroid in the guise of a new $\mathsf{Lie}_\infty$-algebra structure non-trivially deforming the ``big bracket'' for Lie bialgebroids. This exotic $\mathsf{Lie}_\infty$-structure can be interpreted as the equivalent in $d=3$ of the Kontsevich--Shoikhet obstruction to the quantization of infinite-dimensional Poisson manifolds (in $d=2$). We discuss the implications of these results with respect to a conjecture due to P.~Xu regarding the existence of a~quantization map for Lie bialgebroids.}


\Classification{53D55; 18G85; 17B62}

\section{Introduction}\label{Introduction}
Graph complexes play an essential r\^ole in the understanding of the deformation quantization of various algebraic and geometric structures, the paradigmatic example thereof being the Kontsevich graph complex and its relation to the deformation quantization problem for (finite-dimensional) Poisson manifolds \cite{Kontsevich1997}. In particular, the space of Kontsevich quantization maps is acted upon by the pro-unipotent group exponentiating the zeroth cohomology of the Kontsevich graph complex of directed graphs \cite{Dolgushev2011}. As shown by T.~Willwacher \cite{Willwacher2015}, the latter is isomorphic to the Grothendieck--Teichm\"uller group $\GRT_1$~-- introduced by V.~Drinfeld\footnote{Based on a suggestion due to A.~Grothendieck in his {\it Esquisse d'un Programme} \cite{Grothendieck} who proposed to study the combinatorial properties of the absolute Galois group $\text{Gal}\big(\bar{\mathbb Q}/\mathbb Q\big)$ via its natural action on the {tower} of Teichm\"uller groupoids.} \cite{Drinfeld1991} in the context of the absolute Galois group $\text{Gal}\big(\bar{\mathbb Q}/\mathbb Q\big)$ and the theory of quasi-Hopf algebras~-- so that the space of Kontsevich maps\footnote{More precisely, the space of homotopy classes \cite{Dolgushev2007b} of stable \cite{Dolgushev2011} formality morphisms is a $\GRT_1$-torsor.} is a~$\GRT_1$-torsor \cite{KontsevichLett.Math.Phys.48:35-721999}. Since its inception, the Grothendieck--Teichm\"uller group appeared in a variety of mathematical contexts such as the Kashiwara--Vergne conjecture, multiple zeta values, rational homotopy of the $\mathbf{E}_2$-operad, etc.

In the present context of graph complexes, another incarnation of the Grothendieck--Teich\-m\"ul\-ler group can be found in relation to the deformation quantization problem for Lie bialgebras via the action of the graph complex of {\it oriented} graphs \cite{Willwacher2015c} in dimension $d=3$.\footnote{\label{footdim}Recall that the parameter $d$ corresponds to the dimension of the source manifold of the relevant AKSZ $\sigma$-model \cite{AlexandrovKontsevichSchwarzEtAl1997}. The latter is related to the degree~$n$ of the corresponding target manifold via $d=n+1$ and is therefore independent of the dimension of the associated algebro-geometric structure (Poisson manifold, Lie bialgebra, etc.). Consistently, it relates to the dimension of the compactified configuration spaces of points of the associated de Rham field theories \cite{Merkulov2008,Merkulov2010}. Therefore, any graph complex related to Poisson manifolds has dimension $d=2$ while the ones related to Lie bialgebras and generalisations thereof have dimension $d=3$.} The latter action generalises to Lie-quasi bialgebras\footnote{As well as their dual, referred to as quasi-Lie bialgebras in the following, see footnote~\ref{footterminology} for terminology.} and furthermore provides a rationale for the classifying r\^ole played by the Grothendieck--Teichm\"uller group on the space of quantization maps for Lie and Lie-quasi bialgebras \`a la Etingof--Kazhdan \cite{Etingof1995,Sakalos2013}. The oriented graph complex also plays a crucial r\^ole regarding the obstruction theory to the existence of a universal quantization of infinite-dimensional Poisson manifolds \cite{Willwacher2015c}. The corresponding obstruction lives in the first order cohomology of the oriented graph complex in $d=2$ which is a one-dimensional space spanned by the so-called Kontsevich--Shoikhet cocycle. When represented on the space of (infinite-dimensional) polyvector fields, the latter yields an exotic \Lieinf-structure \cite{Shoikhet2008a} deforming non-trivially the Schouten bracket. Further, the zeroth order cohomology of the oriented graph complex (in $d=2$) vanishes thus preventing the Grothendieck--Teichm\"uller group to play a classifying r\^ole for quantizations of infinite-dimensional Poisson manifolds. The deformation quantization problem for infinite-dimensional Poisson manifolds thus differs essentially from the finite-dimensional case and this discrepancy can be traced back to the fact that their respective deformation theory is acted upon by a different graph complex (directed vs.\ oriented). Deformation quantization problems can then be partitioned into different classes according to the cohomology of the graph complex acting on them. We distinguish between three main classes:
\begin{itemize}\itemsep=0pt
\item A \textit{no-go} class comprising the deformation quantization problems for infinite-dimensional Poisson manifolds:
\begin{enumerate}\itemsep=0pt
\item The Grothendieck--Teichm\"uller group plays no classifying r\^ole regarding the universal deformations
(and hence quantizations).
\item There exists a potential obstruction to the existence of universal quantizations.
\end{enumerate}
\item A \textit{yes-go} class comprising the deformation quantization problems for (finite-dimensional) Poisson manifolds and Lie-(quasi) bialgebras for which:
\begin{enumerate}\itemsep=0pt
\item The Grothendieck--Teich\-m\"ul\-ler group plays a classifying r\^ole.
\item There is (conjecturally) no generic obstruction to the existence of universal quantizations.
\end{enumerate}
\item A \textit{middle way} class comprising the deformation quantization problems for Courant algebroids:
\begin{enumerate}\itemsep=0pt
\item The Grothendieck--Teichm\"uller group plays no classifying r\^ole; rather deformations are generated by the triangle cocycle as well as by conformal rescalings associated with trivalent graphs.
\item There is no generic obstruction to the existence of universal quantizations.
\end{enumerate}
\end{itemize}
In the present note, we add two threads to this on-going story by introducing some new universal models for the deformation theory of Lie bialgebroids and their ``quasi'' versions. Lie bialgebroids have been introduced by Mackenzie--Xu \cite{Mackenzie1994} as linearisations of Poisson groupoids and constitute a common generalisation of the notions of Poisson manifolds and Lie bialgebras. The corresponding quantization problem unifies the quantization problems for (finite-dimensional) Poisson manifolds and Lie bialgebras. It was spelled out by P.~Xu \cite{Xu1998a,Xu1999} who then conjectured that any Lie bialgebroid is quantizable. The main result of this note consists in providing some arguments for the non-existence of \textit{universal} quantizations of Lie bialgebroids. This is done by exhibiting a potential obstruction to the existence of a universal quantization map for Lie bialgebroids in the guise of an exotic \Lieinf-structure on the deformation complex of Lie bialgebroids. The latter is a non-trivial deformation of the so-called ``big bracket'' for Lie bialgebroids and can be considered as an avatar in $d=3$ of the Kontsevich--Shoikhet obstruction to the quantization of infinite-dimensional Poisson manifolds. Our main result is stated as follows:
\begin{Theorem}[no-go]\label{TheoremMain}The deformation complex of Lie bialgebroids is endowed with an exotic \Lieinf-structure deforming non-trivially the big bracket of Lie bialgebroids.
\end{Theorem}
It follows from the above considerations that the deformation quantization problem for Lie bialgebroids differs essentially from its Lie bialgebra counterpart and is in fact more akin to the one for infinite-dimensional Poisson manifolds, i.e., it belongs to the {no-go} class.
The origin of this obstruction can be traced back to an action of the graph complex of {\it bi-oriented} graphs (i.e., graphs with two oriented colors) on the deformation theory of Lie bialgebroids. Relaxing the orientation of one of the colors yields an action on the deformation theory of Lie-quasi bialgebroids (and their dual). As a corollary, we find an action of the Grothendieck--Teichm\"uller group on Lie-quasi bialgebroids generalising the one on Lie-(quasi) bialgebras.
\begin{Theorem}[yes-go]\label{thmGro}The Grothendieck--Teichm\"uller group acts via \Lieinf-automorphisms on the deformation complex of both Lie-quasi bialgebroids and quasi-Lie bialgebroids.
\end{Theorem}
Hence, the deformation quantization problem for Lie-quasi bialgebroids differs from its Lie bialgebroid counterpart and resembles more closely the one for Lie bialgebras, i.e., it belongs to the Yes-go class. We conjecture on the basis of these results the existence~-- given a Drinfeld associator~-- of a universal quantization for Lie-quasi bialgebroids (and their dual).

The graph complex approach to deformation quantization\footnote{Note that the case $d=1$ is somehow special among Table~\ref{figgraphsymsum} as there is no associated quantization problem. The relevant cohomological class $H^1(\fcGC_1)\simeq\corps$ is therefore not viewed as an obstruction to quantization but rather as a non-trivial deformation of the symplectic Poisson bracket as a Lie algebra, yielding the Moyal bracket on symplectic manifolds (cf.\ Section~\ref{section:Cohomology} for additional details).} is summed up in Table~\ref{figgraphsymsum}, where the original contribution of the present paper lies at $c=2$.\footnote{The graph complex $\odfGC{d}{j}{i}$ featured in Table~\ref{figgraphsymsum} lives in dimension $d$ (cf.\ footnote~\ref{footdim}) and involves graphs with $c=i+j$ directed colors, $i$ of them are oriented, cf.\ Section~\ref{section:Multi)-oriented graph complexes} for details.}
\begin{table}[h]\centering\small
 \caption{Classification of deformation quantization problems via graph cohomology.} \label{figgraphsymsum}\vspace{1mm}

\begin{tabular}{?c?c?c?c?c?}
\Xhline{2\arrayrulewidth}
\multicolumn{2}{?c?}{}&\multicolumn{3}{c?}{\rule[-0.3cm]{0cm}{0.8cm}\textit{Graph complexes} $\odfGC{d}{j}{i}$}\\
\Xhline{2\arrayrulewidth}
 \multicolumn{2}{?c?}{\rule[-0.3cm]{0cm}{0.8cm}\textit{Class name}}&{no-go}&{yes-go}&{middle way}\\
\Xhline{2\arrayrulewidth}
 \multicolumn{2}{?c?}{\rule[-0.2cm]{0cm}{0.6cm} \multirow{3}{*}{\textit{Model cohomology}}}&\rule[-0.3cm]{0cm}{0.8cm}$H^0(\fcGC_1)\simeq \mathbf{0}$&$H^0(\fcGC_2)\simeq\grt_1$&$H^0(\fcGC_3)\simeq\corps$\\
 \cline{3-5}
\multicolumn{2}{?c?}{} &$H^1(\fcGC_1)\simeq\corps$&\rule[-0.3cm]{0cm}{0.8cm}\hspace{-4mm}$H^1(\fcGC_2)\overset{?}{\simeq}\mathbf{0}$&\hspace{-1mm}$H^1(\fcGC_3)\simeq\mathbf{0}$\\
 \cline{1-5}
 \multicolumn{2}{?c?}{\rule[-0.3cm]{0cm}{0.8cm}\textit{Oriented directions $i$}}&$d-1$&$d-2$&$d-3$\\
\Xhline{2\arrayrulewidth}
$d$&$c=i+j$&\multicolumn{3}{c?}{\rule[-0.3cm]{0cm}{0.8cm}\textit{Actions}}\\
\Xhline{2\arrayrulewidth}
 \rule[-0.3cm]{0cm}{0.8cm}$d=1$&\multirow{7}{*}{$c=1$}&symplectic manifolds&&\\
 \cline{1-1}\cline{3-5}
 \rule[-0.3cm]{0cm}{0.8cm}$d=2$&&{Poisson} ($\dim=\infty$)&{Poisson} ($\dim<\infty$)&\\
 \cline{1-1}\cline{3-5}
 \multirow{6}{*}{$d=3$}&&&\rule[-0.2cm]{0cm}{0.6cm}{Lie bialgebras}&\multirow{2}{*}{proto-Lie bialgebras}\\
 \rule[-0.2cm]{0cm}{0.6cm} &&&\rule[-0.2cm]{0cm}{0.6cm}{Lie-quasi bialgebras}&\\
 \rule[-0.2cm]{0cm}{0.6cm} &&&\rule[-0.2cm]{0cm}{0.6cm}{quasi-Lie bialgebras}&{Courant algebroids}\\
 \cline{2-5}
 \rule[-0.2cm]{0cm}{0.6cm} & \multirow{2}{*}{$c=2$}&\multirow{2}{*}{Lie bialgebroids}&\rule[-0.2cm]{0cm}{0.6cm}Lie-quasi bialgebroids&\multirow{2}{*}{proto-Lie bialgebroids}\\
 \rule[-0.2cm]{0cm}{0.6cm} &&&\rule[-0.2cm]{0cm}{0.6cm}{quasi-Lie bialgebroids}&\\
 \Xhline{2\arrayrulewidth}
 \end{tabular}
\end{table}

{\bf Organisation of this paper.} The original universal model introduced by M.~Kontsevich \cite{Kontsevich1997} takes advantage of the graded geometric interpretation of Poisson manifolds as\footnote{In the remaining of the text, we use the prefix dg to refer to differential graded objects. } dg symplectic manifolds of degree~1. Correspondingly, Section \ref{section:Graded geometry} reviews the formulation of Lie bialgebras and Lie bialgebroids (as well as their generalisations Lie-quasi, quasi-Lie and proto-Lie) as particular dg symplectic manifolds of degree 2. This graded geometric description of the deformation theory of Lie bialgebroids and generalisations thereof will be instrumental in formulating associated universal models in Section \ref{section:Actions}.

The class of universal models introduced in this note involves {\it multi-oriented} graphs, as introduced in~\cite{Zivkovic2017a} and studied in \cite{Merkulov2017} in the context of multi-oriented props and their representations on homotopy algebras with branes. The main definitions and results regarding the cohomology of multi-oriented graph complexes are reviewed in Section~\ref{section:Multi)-oriented graph complexes}.

Following these two review sections, we introduce our main results in Section \ref{section:Actions}. We start by reviewing the known action of the (one-colored) oriented graph complex on Lie-(quasi) bialgebras in Section~\ref{section:Lie bialgebras2} and then move on to the Lie bialgebroid case in Section~\ref{section:Lie bialgebroids2}. Using the cohomological results reviewed in Section~\ref{section:Multi)-oriented graph complexes}, we prove our main results regarding the existence of an exotic \Lieinf-structure for Lie bialgebroids (Theorem~\ref{TheoremMain}) and the action of the Grothendieck--Teichm\"uller group on Lie-quasi bialgebroids (Theorem~\ref{thmGro}).

In view of the results of Section~\ref{section:Lie bialgebroids2}, we formulate two conjectures in Section \ref{section:Application to quantization}: a no-go (Conjecture~\ref{congno}) regarding the existence of (universal) quantizations for Lie bialgebroids and a~yes-go (Conjecture~\ref{congyes}) regarding the one of Lie-quasi bialgebroids (and their dual).

Two appendices conclude the present note. Appendix \ref{section:Appendix} reviews the (ungraded) geometric formulation of Lie bialgebroids and related notions, to be compared with the graded geometric interpretation of Section \ref{section:Graded geometry}. Appendix~\ref{section:Incarnation} contains explicit formulae and additional results regarding the graph cocycle generating the exotic \Lieinf-structure for Lie bialgebroids.

{\bf Conventions.} Throughout the text, we work over a ground field $\corps$ of characteristic zero. The operads introduced in the text live in the category of (graded) vector spaces over $\corps$. Given a~graded $\corps$-vector space $\alg:=\bigoplus_{k\in\mathbb Z}\alg^k$, the $n$-suspended graded vector space $\alg[n]$ is defined via its homogeneous components $\alg[n]^k:=\alg^{k+n}$. We will denote $s\colon \alg[n]\to \alg$ the corresponding suspension map of intrinsic degree $n$.

\section{Graded geometry}\label{section:Graded geometry}
It is a well-known result that a Lie algebra structure on a vector space $\alg$ yields a differential structure on the exterior algebra $\w^\bullet\alg^*$ in the guise of the Chevalley--Eilenberg differential. The exterior algebra can equivalently be recast as the algebra of functions on the shifted vector space~$\alg[1]$ seen as a graded manifold of degree~1 on which the Chevalley--Eilenberg differential defines a homological vector field. Such a supergeometric formulation of Lie algebras was generalised by A.Yu.~Vaintrob \cite{Vaintrob1997} who showed a bijective correspondence between dg manifolds (or $\NQ$-manifolds) of degree~1 and Lie algebroids. On the other hand, it was shown by D.~Roytenberg \cite{Roytenberg2007,Roytenberged.Contemp.Math.Vol.315Amer.Math.Soc.ProvidenceRI2002} that dg symplectic manifolds (or $\NPQ$-manifolds) of degree~1 (resp.\ of degree~2) are in bijective correspondence with Poisson manifolds (resp. Courant algebroids). Our aim in this section is to review how Lie bialgebra and Lie bialgebroid structures (and generalisations\footnote{As reviewed in Appendix~\ref{section:Appendix}.}) can be naturally recast as Hamiltonian functions for a specific graded Poisson algebra of functions on a graded manifold (as pioneered in \cite{Lecomte1990}, cf.\ \cite{Kosmann-Schwarzbach, Roytenberg2007} for details and related constructions). We start by reviewing this graded geometric construction for Lie bialgebras (including proto-Lie, Lie-quasi and quasi-Lie bialgebras) in Section~\ref{section:Lie bialgebras} and then move on to the Lie bialgebroid counterparts of these notions in Section~\ref{section:Lie bialgebroids}.

\subsection{Lie bialgebras}\label{section:Lie bialgebras}

Lie bialgebra structures (and generalisations thereof) on a vector space $\alg$ can be conveniently encoded into particular Hamiltonian functions on the graded manifold $T^*(\alg[1])\simeq(\alg\oplus\alg^*)[1]$ with homogeneous coordinates\footnote{The subscript denotes the corresponding degree.} $\big\{\underset{1}{\xi^a},\underset{1}{\zeta_a}\big\}$, with $a\in\pset{1,\dots,\dim\alg}$. The latter is a graded symplectic manifold with symplectic 2-form $\symp={\rm d}\xi^{a}\w {\rm d}\zeta_{a}$ of degree $2$. The associated Poisson bracket of degree $-2$ acts on homogeneous functions in $\foncg$ as
\begin{gather}
\pb{f}{g}^\alg_\symp=(-1)^{|f|}\left(\frac{\p f}{\p \xi^a} \frac{\p g}{\p\zeta_a}+\frac{\p f}{\p \zeta_a} \frac{\p g}{\p\xi^a}\right).\label{Poissonbi}
\end{gather}
The graded Poisson bracket \eqref{Poissonbi} can be seen as the graded geometric formulation of the ``big bracket'' (introduced by Y.~Kosmann--Schwarzbach \cite{Kosmann-Schwarzbach1992}) acting on $\w^\bullet(\alg\oplus\alg^*)\simeq\foncg$.

Upgrading the graded symplectic manifold $(\alg\oplus\alg^*)[1]$ to a dg symplectic manifold (or $\NPQ$-manifold)\footnote{Or equivalently, endowing the graded Poisson algebra of functions $\big(\foncg,\cdot,\pbdot^\alg_\symp\big)$ with a compatible differential.} allows to define various algebraic structures. A differential structure on a graded symplectic manifold is given by a vector field $\Q$ of degree 1 being homological (i.e., $\br{\Q}{\Q}_\mathsf{Lie}=0$) with respect to the graded Lie bracket of vector fields and compatible with the symplectic 2-form (i.e., $\Lag_\Q\symp=0$). This last compatibility relation ensures\footnote{Via Cartan's homotopy formula, cf.\ \cite[Lemma~2.2]{Roytenberged.Contemp.Math.Vol.315Amer.Math.Soc.ProvidenceRI2002}.} that $\Q$ is necessarily a Hamiltonian vector field, i.e., there exists a function of degree~3 called the Hamiltonian satisfying the structure equation $\pb{\cH}{\cH}^\alg_\symp=0$ and such that $\Q:=\pb{\cH}{\cdot}^\alg_\symp$. The most general function of degree 3 on $(\alg\oplus\alg^*)[1]$ reads explicitly as\footnote{The signs and coefficients are chosen for later convenience.}
\begin{gather}
\cH=-\tfrac{1}{2} f_{a b}{}^c \xi^a\xi^b\zeta_c-\tfrac{1}{2} C_c{}^{a b} \zeta_a\zeta_b\xi^c+\tfrac16\varphi^{abc}\zeta_a\zeta_b\zeta_c+\tfrac16\psi_{abc} \xi^a\xi^b\xi^c,\label{Hamiltonianbialgebra}
\end{gather}
where\footnote{Here and in the following, round (resp.\ square) brackets of indices will denote (skew)symmetrisation.} $f_{a b}{}^c=f_{[a b]}{}^c$, $C_c{}^{a b}=C_c{}^{[a b]}$, $\varphi^{abc}=\varphi^{[abc]}$ and $\psi_{abc}=\psi_{[abc]}$.

The Hamiltonian condition $\pb{\cH}{\cH}^\alg_\symp=0$ translates as a set of 5 constraints on the defining maps $\pset{f,C,\varphi,\psi}$:
\begin{gather}
\bullet \ \mathcal D_1{}_{a b c}{}^d:=-f_{e[a}{}^df_{bc]}{}^e-\psi_{e[ab}C_{c]}{}^{ed}=0,\label{consHam1bialgebra}\\
\bullet \ \mathcal D_2{}_d{}^{a b c}:=-C_d{}^{e[a}C_e{}^{bc]}-\varphi^{e[ab}f_{ed}{}^{c]}{}=0, \label{consHam2bialgebra}\\
\bullet \ \mathcal D_3{}_{ab}{}^{c d}:=2 f_{e[a}{}^{[c}C_{b]}{}^{d]e}-\tfrac{1}{2} f_{ab}{}^eC_{e}{}^{cd}-\tfrac{1}{2} \psi_{eab}\varphi^{ecd}=0, \label{consHam3bialgebra}\\
\bullet \ \mathcal D_4{}^{abcd}:=\tfrac{1}{2} \varphi^{e[ab}C_e{}^{cd]}=0, \label{consHam4bialgebra}\\
\bullet \ \mathcal D_5{}_{abcd}:=\tfrac{1}{2} \psi_{e[ab}f_{cd]}{}^e=0. \label{consHam5bialgebra}
\end{gather}
A set of maps $\pset{f,C,\varphi,\psi}$ satisfying the constraints \eqref{consHam1bialgebra}--\eqref{consHam5bialgebra} form the components of a \textit{proto-Lie bialgebra} on $(\alg,\alg^*)$ (cf.\ Appendix~\ref{section:Appendix} for a definition) whose deformation theory is therefore controlled by the dg Lie algebra\footnote{\label{footterminology}Note that the graded Poisson bracket has intrinsic degree $-2$ on $\foncg$. To recover the usual grading, one needs to consider the 2-suspension $\foncg[2]$.} $\big(\foncg,\Q,\pbdot^\alg_\symp\big)$. Proto-Lie bialgebras thus constitute the most general notion in the bialgebra realm and other structures (Lie-quasi, quasi-Lie and Lie bialgebras) will be defined as particular cases thereof.

The remainder of this section will therefore introduce several particular graded Poisson subalgebras of $\foncg$ whose Hamiltonian functions will encode various sub-classes of proto-Lie bialgebras. Let us start by defining the subspace $\A^\alg_\text{\rm Lie-quasi}\subset\foncg$ as $\displaystyle\A^\alg_\text{\rm Lie-quasi}:=\pset{f\in\foncg\, \big|\, f|_{\zeta=0}=0}$. In plain words, the subspace $\A^\alg_\text{\rm Lie-quasi}$ is obtained by discarding all functions of the form $\psi_{a_1\cdots a_m} \xi^{a_1}\cdots\xi^{a_m}$, for arbitrary values of $m\geq0$. It can be easily checked that $\A^\alg_\text{\rm Lie-quasi}$ is preserved by both the pointwise product of functions and the graded Poisson bracket~\eqref{Poissonbi} and thus defines a graded Poisson subalgebra of $\foncg$. The most general Hamiltonian function of $\A^\alg_\text{\rm Lie-quasi}$ reads as~\eqref{Hamiltonianbialgebra} with $\psi\equiv0$, where the maps $\pset{f,C,\varphi}$ satisfy \eqref{consHam1bialgebra}--\eqref{consHam4bialgebra} with $\psi\equiv0$. In particular, imposing $\psi\equiv0$ in equation~\eqref{consHam1bialgebra} ensures that the map $f$ defines a genuine Lie algebra structure on $\alg$ (while the structure defined on $\alg^*$ is still ``quasi'' due to the presence of $\varphi$). The resulting equations reproduce the defining conditions of a \textit{Lie-quasi bialgebra} on $(\alg,\alg^*)$ as introduced by Drinfeld in~\cite{Drinfeld1989} as semi-classicalisation of the notion of quasi-bialgebra.\footnote{Remark that Lie-quasi bialgebras were denoted ``quasi-Lie bialgebras'' in~\cite{Drinfeld1989}. We follow the terminology used in \cite{Kosmann-Schwarzbach} where the term Lie-quasi bialgebras is used for Lie algebras which fail to be Lie bialgebras (so that they are only ``quasi'' bialgebras) while the term quasi-Lie was reserved for the dual counterpart (not considered in \cite{Drinfeld1989}) where the Jacobi identity for $f$ is only ``quasi'' satisfied.}

Dually to the previous case, one defines the graded Poisson subalgebra $\A^\alg_\text{\rm quasi-Lie}$ as the subspace obtained by discarding all functions of the form $\varphi^{a_1\cdots a_n}\zeta_{a_1}\cdots\zeta_{a_n}$, for all $n\geq0$, i.e., $\A^\alg_\text{\rm quasi-Lie}:=\pset{f\in\foncg\, \big|\, f|_{\xi=0}=0}\subset\foncg$. The most general Hamiltonian function of $\A^\alg_{\rm quasi-Lie}$ reads as \eqref{Hamiltonianbialgebra} with $\varphi\equiv0$, where the functions $\pset{f,C,\psi}$ satisfy \eqref{consHam1bialgebra}--\eqref{consHam3bialgebra} and \eqref{consHam5bialgebra} with~$\varphi\equiv0$. Dually to the Lie-quasi case, setting $\varphi\equiv0$ in equation~\eqref{consHam2bialgebra} ensures that the map~$C$ defines a genuine Lie algebra structure on $\alg^*$ (while the structure defined on $\alg$ is only ``quasi'' Lie due to the presence of~$\psi$). The resulting equations reproduce the defining conditions of a \textit{quasi-Lie bialgebra} on $(\alg,\alg^*)$ as introduced and studied in~\cite{Bangoura1993, Kosmann-Schwarzbach1992}.

Finally, let us define the subspace $\A^\alg_{\rm Lie}:=\pset{f\in\foncg\, \big|\, f|_{\xi=0}=0\text{ and } f|_{\zeta=0}=0}$, i.e., $\A^\alg_{\rm Lie}$ is defined as the intersection $\A^\alg_{\rm Lie}:= \A^\alg_\text{\rm Lie-quasi}\cap \A^\alg_\text{\rm quasi-Lie}$ between the two previous Poisson subalgebras. The latter subspace can again be checked to be a graded Poisson subalgebra of $\foncg$ obtained by discarding all functions of the form $\psi_{a_1\cdots a_m} \xi^{a_1}\cdots\xi^{a_m}$ and $\varphi^{a_1\cdots a_n}\zeta_{a_1}\cdots\zeta_{a_n}$ for all $m,n\geq0$. In particular, the most general Hamiltonian function of~$\A^\alg_{\rm Lie}$ reads as~\eqref{Hamiltonianbialgebra} with $\varphi\equiv0$ and $\psi\equiv0$, where the functions $\pset{f,C}$ satisfy \eqref{consHam1bialgebra}--\eqref{consHam3bialgebra} with $\varphi\equiv0$ and $\psi\equiv0$. In particular, constraint~\eqref{consHam1bialgebra} \big(resp. \eqref{consHam2bialgebra}\big) ensures that the map~$f$ (resp.~$C$) defines a genuine Lie algebra structure on $\alg$ (resp.~$\alg^*$). These two Lie algebras are furthermore compatible with each other due to~\eqref{consHam3bialgebra} and hence define a \textit{Lie bialgebra} on $(\alg,\alg^*)$ (cf.~\eqref{combi}). We sum up the previous discussion in the following proposition:

\begin{Proposition}\label{propbialgHam}
Let $\alg$ be a vector space. The following correspondences hold:
\begin{itemize}\itemsep=0pt
\item Hamiltonians in $\foncg$ are in bijective correspondence with proto-Lie bialgebra structures on $(\alg,\alg^*)$.
\item \hspace{1.5cm}``\hspace{1.0cm}$\A^\alg_\text{\rm Lie-quasi}$\hspace{1.5cm}``\hspace{1.3cm} Lie-quasi bialgebra \hspace{0.45cm}``\hspace{0.45cm}.
\item \hspace{1.5cm}``\hspace{1.0cm}$\A^\alg_\text{\rm quasi-Lie}$ \hspace{1.4cm}``\hspace{1.4cm}quasi-Lie bialgebra \hspace{0.45cm}``\hspace{0.45cm}.
\item \hspace{1.5cm}``\hspace{1.0cm}$\A^\alg_\text{\rm Lie}$ \hspace{2.2cm}``\hspace{1.9cm}Lie bialgebra \hspace{0.91cm}``\hspace{0.45cm}.
\end{itemize}
\end{Proposition}

This interpretation of the deformation theory for Lie bialgebras and generalisations as graded Poisson algebras will be put to use in Section \ref{section:Actions} where will be discussed universal models thereof.
In the next section, we turn to the generalisation of this graded geometric interpretation to the larger class of Lie bialgebroids and variations thereof.

\subsection{Lie bialgebroids}
\label{section:Lie bialgebroids}
Letting $E\overset{\pi}{\to} \M$ be a vector bundle over the smooth (finite-dimensional) manifold\footnote{\label{footpoint}The restriction from the ``algebroid'' case to the ``algebra'' case can be done by assuming that $\M$ is the one-point manifold so that $E\simeq \alg$ becomes a $\corps$-vector space.}$\M$, the relevant graded Poisson algebra is the algebra of functions of the graded symplectic manifold $T^*[2]E[1]$, defined as the (2-shifted) cotangent bundle of the (1-shifted\footnote{Here and in the following, $A[n]$ will denote the vector bundle obtained by shifting the grading of the fiber of the vector bundle $A$ by $n$. }) vector bundle $E$, and denoted $\TE\equiv T^*[2]E[1]$ in the following.
\noindent The graded manifold $\TE$ is of degree 2 and can be locally\footnote{Correspondingly, all formulae appearing in this note will be local.} coordinatised by the set of homogeneous coordinates $\big\{\underset{0}{x^\mu},\underset{1}{\xi^a},\underset{1}{\zeta_a},\underset{2}{p_\mu}\big\}$ so that the symplectic 2-form of degree $2$ can be written as
\[
\symp={\rm d}x^\mu\w {\rm d}p_\mu+{\rm d}\xi^{a}\w {\rm d}\zeta_{a}.
\]
\noindent The associated Poisson bracket of degree $-2$ acts as follows on homogeneous functions $f,g\in\foncTE$:
\begin{gather}
\pb{f}{g}^E_\symp=\frac{\p f}{\p x^\mu} \frac{\p g}{\p p_\mu}-\frac{\p f}{\p p_\mu} \frac{\p g}{\p x^\mu}+(-1)^{|f|}\left(\frac{\p f}{\p \xi^a} \frac{\p g}{\p\zeta_a}+\frac{\p f}{\p \zeta_a} \frac{\p g}{\p\xi^a}\right).\label{eqPois}
\end{gather}
The latter is sometimes referred to as the ``big bracket'' for Lie bialgebroids. Upgrading the graded symplectic manifold~$\TE$ to a dg symplectic manifold will allow to define various geometric structures. Following the same path as in the bialgebra case, we introduce a compatible differential through a Hamiltonian function. The most general function of degree~3 on~$\TE$ reads
\begin{gather}
\cH=\rho_a{}^\mu(x) \xi^ap_\mu-\tfrac{1}{2} f_{a b}{}^c(x) \xi^a\xi^b\zeta_c+R^{a|\mu}(x) \zeta_ap_\mu-\tfrac{1}{2} C_c{}^{a b}(x) \zeta_a\zeta_b\xi^c\nonumber\\
\hphantom{\cH=}{}
+\tfrac16\varphi^{abc}(x) \zeta_a\zeta_b\zeta_c+\tfrac16\psi_{abc}(x) \xi^a\xi^b\xi^c, \label{Hamiltonian}
\end{gather}
where $\pset{\rho,f,R,C,\varphi,\psi}$ are functions on the base space $\M$, with symmetries $f_{a b}{}^c=f_{[a b]}{}^c$, $C_c{}^{a b}=C_c{}^{[a b]}$, $\varphi^{abc}=\varphi^{[abc]}$ and $\psi_{abc}=\psi_{[abc]}$. Imposing the Hamiltonian constraint $\pb{\cH}{\cH}^E_\symp=0$ yields a set of 9 conditions on the defining functions $\pset{\rho,f,R,C,\varphi,\psi}$ that we denote as follows
\begin{gather}
\bullet \ \mathcal C_1{}_{a b}{}^\mu:=2 \rho_{[a}{}^\lambda\p_\lambda\rho_{b]}{}^\mu-\rho_c{}^\mu f_{ab}{}^c+R^{c|\mu}\psi_{cab}=0,\label{consHam1}\\
\bullet \ \mathcal C_2{}_{a b c}{}^d:=\rho_{[a}{}^\lambda\p_\lambda f_{bc]}{}^d-f_{e[a}{}^df_{bc]}{}^e+\tfrac13 R^{d|\lambda}\p_\lambda\psi_{abc}-\psi_{e[ab}C_{c]}{}^{ed}=0,\label{consHam2}\\
\bullet \ \mathcal C_3{}^{a b|\mu}:=2 R^{[a|\lambda}\p_\lambda R^{b]|\mu}-R^{c|\mu}C_{c}{}^{ab}+\rho_c{}^\mu\varphi^{cab}=0,\label{consHam3}\\
\bullet \ \mathcal C_4{}_d{}^{a b c}:=R^{[a|\lambda}\p_\lambda C_d{}^{bc]}-C_d{}^{e[a}C_e{}^{bc]}+\tfrac13 \rho_{d}{}^{\lambda}\p_\lambda\varphi^{abc}-\varphi^{e[ab}f_{ed}{}^{c]}{}=0,\label{consHam4}\\
\bullet \ \mathcal C_5{}^{\mu \nu}:=R^{a(\mu}\rho_a{}^{\nu)}=0,\label{consHam5}\\
\bullet \ \mathcal C_6{}_{a}{}^{b|\mu}:=\rho_a{}^\lambda\p_\lambda R^{b|\mu}-R^{b|\lambda}\p_\lambda\rho_a{}^\mu-\rho_c{}^\mu C_a{}^{bc}-R^{c|\mu}f_{ca}{}^b=0,\label{consHam6}\\
\bullet \ \mathcal C_7{}_{ab}{}^{c d}:=\rho_{[a}{}^\lambda\p_\lambda C_{b]}{}^{c d}+R^{[c|\lambda}\p_\lambda f_{ab}{}^{d]}+2 f_{e[a}{}^{[c}C_{b]}{}^{d]e}-\tfrac{1}{2} f_{ab}{}^eC_{e}{}^{cd}-\tfrac{1}{2} \psi_{eab}\varphi^{ecd}=0,\label{consHam7}\\
\bullet \ \mathcal C_8{}^{abcd}:=\tfrac13R^{[d|\lambda}\p_\lambda\varphi^{abc]}+\tfrac{1}{2} \varphi^{e[ab}C_e{}^{cd]}=0,\label{consHam8}\\
\bullet \ \mathcal C_9{}_{abcd}:=\tfrac13\rho_{[d}{}^\lambda\p_\lambda\psi_{abc]}+\tfrac{1}{2} \psi_{e[ab}f_{cd]}{}^e=0. \label{consHam9}
\end{gather}
The latter constraints identify with the component expressions of the defining conditions of a~\textit{proto-Lie bialgebroid} on $(E,E^*)$ (compare with \eqref{curvproto1}--\eqref{curvproto5}). The graded Poisson algebra~$\foncTE$ admits several subalgebras defining in turn various sub-classes of proto-Lie bialgebroids. A~convenient way to characterise these Poisson subalgebras is as vanishing ideals of particular Lagrangian submanifolds.\footnote{Letting ${\mathcal N}$ be a (graded) manifold, the \textit{vanishing ideal} of a (graded) submanifold $\cC\subset{\mathcal N}$ is defined as the subalgebra of functions $I_\cC:=\pset{f\in\fonc{{\mathcal N}}\, \big |\, f|_\cC=0}$.
Moreover, a (graded) submanifold $\cC$ of a (graded) symplectic manifold $({\mathcal N},\symp)$ is said to be \textit{Lagrangian} if it is maximally isotropic, i.e.,
\begin{enumerate}\itemsep=0pt
\item The restriction of the symplectic form $\Om$ to $\cC$ vanishes, i.e., $\Om|_\cC=0$.
\item The submanifold $\cC$ has maximal dimension $\dim \cC=\half \dim {\mathcal N}$.
\end{enumerate}
For graded Lagrangian submanifolds, it will be assumed that the underlying bosonic manifold of~$\cC$ (whose algebra of functions is coordinatised by coordinates of degree $0$) identifies with the one of~$\mN$.
Lagrangian submanifolds are coisotropic, so that the corresponding vanishing ideal is closed under the (graded) Poisson bracket associated to~$\Omega$ and hence is a Poisson subalgebra of the algebra of functions on ${\mathcal N}$.} As noted by D.~Roytenberg~\cite{Roytenberg2001}, both $E[1]$ and $E^*[1]$ are Lagrangian submanifolds of $\TE$, thus motivating to consider the Poisson subalgebras of functions vanishing on them. We start by defining the Poisson subalgebra~$\A^E_\text{\rm Lie-quasi}\subset\foncTE$ as the vanishing ideal of~$E[1]$.
In plain words, the subalgebra $\A^E_\text{\rm Lie-quasi}$ is obtained by discarding all functions of the form $\psi_{a_1\cdots a_m}(x) \xi^{a_1}\cdots\xi^{a_m}$, for all $m\geq0$.
The most general Hamiltonian function of $\A^E_\text{\rm Lie-quasi}$ reads as \eqref{Hamiltonian} with $\psi\equiv0$, where the functions $\pset{\rho,f,R,C,\varphi}$ satisfy \eqref{consHam1}--\eqref{consHam8} with $\psi\equiv0$. In particular, setting $\psi\equiv0$ in equations~\eqref{consHam1}--\eqref{consHam2} ensures that the pair $\pset{\rho,f}$ defines a~genuine Lie algebroid structure on $E$. The resulting equations reproduce the defining conditions of a \textit{Lie-quasi bialgebroid} on $(E,E^*)$. The pair $\big((\TE,\Om,\cH),E[1]\big)$ defines a~Manin pair,\footnote{A \textit{Manin pair} of degree $n$ is defined as a dg symplectic manifold $(\mN,\Om,\cH)$ of degree $n$ supplemented with a dg Lagrangian submanifold (also called a $\Lambda$-structure \cite{Severa2007}), i.e., a graded Lagrangian submanifold $\cC\subset\mN$ such that $\cH\in I_\cC$, with~$I_\cC$ the vanishing ideal of $\cC$. Manin pairs of degree $0$ correspond to a pair formed by a~(bosonic) symplectic manifold endowed with a Lagrangian submanifold. Manin pairs of degree $1$ are in bijective correspondence with pairs composed of a Poisson manifold together with a coisotropic submanifold while Manin pairs of degree $2$ identify with Courant algebroids endowed with a Dirac structure~\cite{Severa2007}.} in the terminology of Roytenberg~\cite{Roytenberg2001}.

Dually to the previous case, we define the graded Poisson subalgebra $\A^E_\text{\rm quasi-Lie}\subset\foncTE$ as the vanishing ideal of $E^*[1]$, i.e., as the subspace obtained by discarding all functions of the form $\varphi^{a_1\cdots a_n}(x) \zeta_{a_1}\cdots\zeta_{a_n}$ with $n\geq0$. The most general Hamiltonian function of $\A^E_\text{\rm quasi-Lie}$ reads as~\eqref{Hamiltonian} with $\varphi\equiv0$, where the functions $\pset{\rho,f,R,C,\psi}$ satisfy \eqref{consHam1}--\eqref{consHam7} and~\eqref{consHam9} with~$\varphi\equiv0$. Dually to the Lie-quasi case, setting~$\varphi\equiv0$ in equations~\eqref{consHam3}--\eqref{consHam4} ensures that the pair~$\pset{R,C}$ defines a genuine Lie algebroid structure on~$E^*$ (while the structure defined on~$E$ by~$\pset{\rho,f}$ is still ``quasi'' due to the presence of~$\psi$). The resulting equations reproduce the defining conditions of a dual structure on $(E,E^*)$, dubbed \textit{quasi-Lie bialgebroid} in~\cite{Kosmann-Schwarzbach}. Quasi-Lie bialgebroids are then equivalently characterised as Manin pairs of the form $\big((\TE,\Om,\cH),E^*[1]\big)$.

 We conclude by defining the subspace $\A^E_{\rm Lie}:= \A^E_\text{\rm Lie-quasi}\cap \A^E_\text{\rm quasi-Lie}$.
 The latter subspace can again be checked to be a graded Poisson subalgebra of $\foncTE$ obtained by discarding all functions of the form $\psi_{a_1\cdots a_m}(x)\, \xi^{a_1}\cdots\xi^{a_m}$ or $\varphi^{a_1\cdots a_n}(x) \zeta_{a_1}\cdots\zeta_{a_n}$, for all $m,n\geq0$. In particular, the most general Hamiltonian function of $\A^E_{\rm Lie}$ reads as \eqref{Hamiltonian} with $\varphi\equiv0$ and $\psi\equiv0$, where the functions $\pset{\rho,f,R,C}$ satisfy \eqref{consHam1}--\eqref{consHam7} with $\varphi\equiv0$ and $\psi\equiv0$. Constraints \eqref{consHam1}--\eqref{consHam2} (resp.\ \eqref{consHam3}--\eqref{consHam4}) ensure that the pair $\pset{\rho,f}$ (resp.~$\pset{R,C}$) defines a genuine Lie algebroid structure on $E$ (resp. $E^*$). These two Lie algebroids are furthermore compatible with each other due to \eqref{consHam4}--\eqref{consHam7} and hence define a \textit{Lie bialgebroid} structure on $(E,E^*)$ (cf.\ Appendix~\ref{section:Appendix}).
Lie bialgebroids are furthermore equivalently characterised\footnote{We refer to Appendix~\ref{section:Appendix} for analogues of these statements within the framework of Courant algebroids.} as Manin triples\footnote{A \textit{Manin triple} is a dg symplectic manifold supplemented with two {\it transverse} dg Lagrangian submanifolds $\cC,\dD\subset\mN$ such that $\cH\in I_\cC$ and $\cH\in I_\dD$.} of the form $((\TE,\Om,\cH),E[1], E^*[1])$.

As usual, gauge transformations for the Hamiltonian function \eqref{Hamiltonian} are generated by functions\footnote{Gauge transformations for Lie-quasi (resp. quasi-Lie) bialgebroids are generated by~\eqref{gauge} with $\om\equiv0$ (resp. $\Lambda\equiv0$), and with $\Lambda\equiv0\equiv\om$ for Lie bialgebroids.} of degree $2$ in $\foncTE$ reading explicitly:
\begin{gather}
\mathscr X=X^\mu p_\mu+\lambda^a{}_b \xi^b\zeta_a-\tfrac{1}{2}\Lambda^{ab} \zeta_a\zeta_b-\tfrac{1}{2}\om_{ab} \xi^a\xi^b,\label{gauge}
\end{gather}
so that the gauge algebra for proto-Lie bialgebroids is isomorphic to $\vf\oplus\End\big({\fE}\big)\oplus\Gamma\big({\w}^2E\big)\oplus\Gamma\big({\w}^2E^*\big)$ as a vector space.
We refer to \cite[Section~5]{Roytenberg2001} for a careful treatment of the group structure of these gauge transformations.
The various parameters in~\eqref{gauge} can be interpreted\footnote{We refer to Appendix~\ref{section:Appendix} for more details on this interpretation.} in terms of infinitesimal morphisms of proto-Lie bialgebroids as follows:
\begin{itemize}\itemsep=0pt
\item $X\in\vf$: diffeomorphism of $\M$,
\item $\lambda\in\End\big({\fE}\big)$: rotation of the fibers of $(E,E^*)$,
\item $\Lambda\in\Gamma\big({\w}^2E\big)$ twist of Lie-quasi bialgebroids,
\item $\om\in\Gamma\big({\w}^2E^*\big)$ twist of quasi-Lie bialgebroids,
\end{itemize}
where the denomination \textit{twist} refers to the following construction: given a proto-Lie bialgebroid on $(E,E^*)$ defined by the Hamiltonian function $\cH$, one can define a new proto-Lie bialgebroid structure on $(E,E^*)$ via twisting\footnote{The twisting procedure was introduced in \cite{Drinfeld1989} for quasi-Hopf algebras, in~\cite{Kosmann-Schwarzbach1992} for quasi-Lie bialgebras and in~\cite{Roytenberg2001} for quasi-Lie bialgebroids.} by an arbitrary bivector $\Lambda\in\Gamma\big({\w}^2E\big)$. Explicitly, the twisted proto-Lie bialgebroid $\cH_\Lambda$ is obtained by performing a canonical transformation generated by the flow of the Hamiltonian vector field $\mathscr X_\Lambda:=\pb{{-}\Lambda}{\cdot}^E_\symp$ where $\Lambda:=\half \Lambda^{ab}\zeta_a\zeta_b$ is a function of degree $2$ in $\foncTE$. Explicitly, a twist by $\Lambda$ corresponds to the canonical transformation\footnote{Dually, one can consider twisting by a 2-form field $\om\in\Gamma\big({\w}^2 E^*\big)$ corresponding to the canonical transformation generated by the flow of the Hamiltonian vector field $\mathscr X_\om:=\pb{-\om}{\cdot}$, where $\om:=\half \om_{ab}\, \xi^a\xi^b$:
\begin{gather}
x^\mu\overset{\om}{\longrightarrow} x^\mu ,\qquad p_\mu\overset{\om}{\longrightarrow} p_\mu-\half \p_\mu \om_{ab}\xi^a\xi^b ,\qquad\xi^a\overset{\om}{\longrightarrow} \xi^a ,\qquad \zeta_a\overset{\om}{\longrightarrow} \zeta_a-\om_{ab}\xi^b\label{canonicaltwistom}
\end{gather}
which amounts to a shift of the components of $\cH$ as in equation~\eqref{twistomexp}.}
\begin{gather*}
x^\mu\overset{\Lambda}{\longrightarrow} x^\mu ,\qquad p_\mu\overset{\Lambda}{\longrightarrow} p_\mu-\tfrac{1}{2} \p_\mu \Lambda^{ab}\zeta_a\zeta_b , \qquad\xi^a\overset{\Lambda}{\longrightarrow} \xi^a-\Lambda^{ab}\zeta_b ,\qquad \zeta_a\overset{\Lambda}{\longrightarrow} \zeta_a,
\end{gather*}
which amounts to a shift of the components of $\cH$ (cf.\ equation~\eqref{twistLambdaexp} in Appendix~\ref{section:Appendix} for explicit expressions). The previous canonical transformation maps Lie-quasi bialgebroids to Lie-quasi bialgebroids but generically fails to preserve quasi-Lie bialgebroids, as is consistent with the fact that $\Lambda\notin\A^E_\text{\rm quasi-Lie}$ so that $\Gamma\big({\w}^2E\big)$ is not part of the gauge algebra for quasi-Lie bialgebroids.\footnote{Dually, the canonical transformation \eqref{canonicaltwistom} preserves the space of quasi-Lie bialgebroids, but generically maps Lie-quasi bialgebroids to proto-Lie bialgebroids. This is consistent with the fact that $\om\notin\A^E_\text{\rm Lie-quasi}$ so that $\Gamma\big({\w}^2E^*\big)$ is not part of the gauge algebra for quasi-Lie bialgebroids.} We refer to Appendix~\ref{section:Appendix} for more details on the twisting procedure.

 We sum up the previous discussion in the following proposition, generalising Proposition~\ref{propbialgHam} to bialgebroids:
\begin{Proposition}
Let $E\overset{\pi}{\to} \M$ be a vector bundle. The following correspondences hold:
\begin{itemize}\itemsep=0pt
\item Hamiltonians in $\foncTE$ are in bijective correspondence with proto-Lie bialgebroid structures on $(E,E^*)$.
\item \hspace{0.5cm}``\hspace{0.5cm}$\A^E_\text{\rm Lie-quasi}$\hspace{1.5cm}``\hspace{1.3cm} Lie-quasi bialgebroid \hspace{0.35cm}``\hspace{0.35cm}.
\item \hspace{0.5cm}``\hspace{0.5cm}$\A^E_\text{\rm quasi-Lie}$ \hspace{1.4cm}``\hspace{1.4cm}quasi-Lie bialgebroid \hspace{0.35cm}``\hspace{0.35cm}.
\item \hspace{0.5cm}``\hspace{0.5cm}$\A^E_\text{\rm Lie}$ \hspace{2.2cm}``\hspace{1.9cm}Lie bialgebroid \hspace{0.85cm}``\hspace{0.28cm}.
\end{itemize}
\end{Proposition}
As noted previously, the graded geometric interpretation of algebro-geometric structures is instrumental to the construction of corresponding universal models, formulated in terms of graph complexes. The next section will introduce the relevant graph complexes which will be shown in Section \ref{section:Actions} to act on the various (sub)-algebras previously introduced.

\section{(Multi)-oriented graph complexes}\label{section:Multi)-oriented graph complexes}
The aim of this section is to review the definition and main results regarding multi-oriented graph complexes and their cohomology, as introduced and studied in \cite{Merkulov2017, Zivkovic2017a,Zivkovic2017} (cf.\ also~\cite{Merkulov2019} for a review). Graph complexes are most clearly defined as deformation complexes of a suitable morphism of operads~\cite{Merkulov2007}.\footnote{Or equivalently, as convolution Lie algebras constructed from suitable graph operads.} We start by introducing the relevant graph operads of multi-directed and multi-oriented graphs from a combinatorial point of view (Section~\ref{section:Operads}) before moving to the definition of the associated graph complexes. We conclude by discussing known results regarding the cohomology of multi-oriented graph complexes (Section \ref{section:Cohomology}) by putting the emphasis on some particular classes relevant for our purpose (cf.\ Section~\ref{section:Actions}).

\subsection{Directed, oriented and sourced graphs}\label{section:Operads}

We will denote $\ugra_{N,k}$ (resp.\ $\gra_{N,k}$) the set of multi(di)graphs\footnote{Recall that \textit{multi$($di$)$graphs} are undirected (resp. directed) graphs allowed to contain both loops and multiple edges.} with $N$ vertices and $k$ directed edges.

The set $\dgra{N,k}{c}$ of \textit{multi-directed graphs} with $c$ colors is defined as the set of ordered pairs $(\graphGam,\mathfrak c)$ where:
\begin{itemize}\itemsep=0pt
\item $\graphGam\in\gra_{N,k}$ is a multidigraph. We will denote $V_\graphGam$ (resp.~$E_\graphGam$) the set of vertices (resp. edges) of $\graphGam$.
 \item $\mathfrak c$ stands for a map $\mathfrak c\colon E_\graphGam\times[c-1]\to \pset{+,-}$ where $c\in\mathbb N$ stands for the total number of \textit{colors}\footnote{That is, including the underlying (black) arrow of the multidigraph.} and $[c-1]:=\pset{1,2,\dots,c-1}$.
\end{itemize}
 A pictorial representation of a multi-directed graph in $\dgra{N,k}{c}$ can be given by decorating each directed edge of the underlying multidigraph in $\gra_{N,k}$ with $c-1$ additional arrows of different colors (cf.\ Figure~\ref{figgraphexmult} for an example).
 \begin{figure}[ht]\centering
 \text{}\raisebox{-0.8ex}{\hbox{\begin{tikzpicture}[scale=0.5, >=stealth']
\tikzstyle{w}=[circle, draw, minimum size=4, inner sep=1]
\tikzstyle{b}=[circle, draw, fill, minimum size=2, inner sep=0.02]
\draw (0,0) node [ext] (b1) {\rm 1};
\draw (4,0) node [ext] (b2) {\rm 2};
\draw (4,-4) node [ext] (b3) {\rm 3};
\draw (0,-4) node [ext] (b4) {\rm 4};
\draw (2,0.4) node (b5) {\rm \scriptsize$i$};
\draw (2,-4.4) node (b6) {\rm \scriptsize$ii$};
\draw (2,-1.4) node (b7) {\rm \scriptsize$iii$};
\draw (4.6,-2) node (b8) {\rm \scriptsize$iv$};
\draw (-0.4,-2) node (b8) {\rm \scriptsize$v$};
\draw[black,->,>=latex,line width=0.2mm,scale=1.35][decoration={markings,mark=at position .3 with {\arrow[scale=1.35,linkcolor]{latex}}},
 postaction={decorate},
 shorten >=0.4pt][decoration={markings,mark=at position .45 with {\arrowreversed[scale=1.35,carrotorange]{latex}}},
 postaction={decorate},
 shorten >=0.4pt][decoration={markings,mark=at position .85 with {\arrow[scale=1.35,cerulean]{latex}}},
 postaction={decorate},
 shorten >=0.4pt] (b4) to (b2);
\draw[black,->,>=latex,line width=0.2mm,scale=1.35][decoration={markings,mark=at position .3 with {\arrow[scale=1.35,linkcolor]{latex}}},
 postaction={decorate},
 shorten >=0.4pt][decoration={markings,mark=at position .4 with {\arrowreversed[scale=1.35,carrotorange]{latex}}},
 postaction={decorate},
 shorten >=0.4pt][decoration={markings,mark=at position .66 with {\arrowreversed[scale=1.35,cerulean]{latex}}},
 postaction={decorate},
 shorten >=0.4pt] (b1) to (b2);
\draw[black,->,>=latex,line width=0.2mm,scale=1.35][decoration={markings,mark=at position .3 with {\arrow[scale=1.35,linkcolor]{latex}}},
 postaction={decorate},
 shorten >=0.4pt][decoration={markings,mark=at position .4 with {\arrowreversed[scale=1.35,carrotorange]{latex}}},
 postaction={decorate},
 shorten >=0.4pt][decoration={markings,mark=at position .85 with {\arrow[scale=1.35,cerulean]{latex}}},
 postaction={decorate},
 shorten >=0.4pt] (b2) to (b3);
\draw[black,<-,>=latex,line width=0.2mm,scale=1.35][decoration={markings,mark=at position .3 with {\arrow[scale=1.35,linkcolor]{latex}}},
 postaction={decorate},
 shorten >=0.4pt][decoration={markings,mark=at position .55 with {\arrow[scale=1.35,carrotorange]{latex}}},
 postaction={decorate},
 shorten >=0.4pt][decoration={markings,mark=at position .85 with {\arrow[scale=1.35,cerulean]{latex}}},
 postaction={decorate},
 shorten >=0.4pt] (b3) to (b4);
\draw[black,<-,>=latex,line width=0.2mm,scale=1.35][decoration={markings,mark=at position .3 with {\arrow[scale=1.35,linkcolor]{latex}}},
 postaction={decorate},
 shorten >=0.4pt][decoration={markings,mark=at position .4 with {\arrowreversed[scale=1.35,carrotorange]{latex}}},
 postaction={decorate},
 shorten >=0.4pt][decoration={markings,mark=at position .66 with {\arrowreversed[scale=1.35,cerulean]{latex}}},
 postaction={decorate},
 shorten >=0.4pt] (b4) to (b1);
\end{tikzpicture}}}
 \caption{Example of a multi-directed graph in $\dgra{4,5}{4}$.} \label{figgraphexmult}
\end{figure}

\noindent The direction of the arrow of color $i$ on the edge $e$ is aligned with the one of the black arrow if $\mathfrak c(e,i)=+$ and opposite to it if $\mathfrak c(e,i)=-$. There is a natural right-action of the permutation group $\mS_N$ (resp.~$\mS_k$) on elements of~$\dgra{N,k}{c}$ by permutation of the labeling of vertices (resp.\ edges).

{\bf Operads.} For all $d\in\mathbb N^*$, we define the collection $\pset{\dGra{\dime}{c}(N)}_{N\geq1}$ of $\mathbb S_N$-modules:
\begin{gather}
\dGra{\dime}{c}(N):=\bigoplus_{k\geq0}\big(\corps\dl\dgra{N,k}{c}\dr\otimes_{\mS_k}\sgn_k^{\otimes |d-1|}\big)[k(\dime-1)],\label{subop}
\end{gather}
where $\sgn_k$ denotes the $1$-dimensional sign representation of $\mathbb S_k$. The subscript stands for taking coinvariants with respect to the diagonal right action of $\mS_k$ and the term between brackets denotes degree suspension. In plain words, this means that edges carry an intrinsic degree $1-\dime$ and are bosonic for $d$ odd and fermionic for $d$ even. The $\mathbb S$-module $\pset{\dGra{\dime}{c}(N)}_{N\geq1}$ can further be given the structure of an operad\footnote{The identity element $\Id\in\dGra{\dime}{c}(1)$ is defined as the graph $\Id:=
\vcenter{\hbox{\begin{tikzpicture}[scale=0.5, >=stealth']
\node [ext] (b1) at (0,0) {\rm 1};
\end{tikzpicture}}}$.} by endowing it with the usual equivariant partial composition operations $\circ_i\colon\dGra{\dime}{c}(M)\otimes\dGra{\dime}{c}(N)\to\dGra{\dime}{c}(M+N-1)$, cf.\ Figure~\ref{figtypicaloddmulti} for an example and, e.g., of \cite[Section~4]{Morand2019a} for more details.
\begin{figure}[ht]\centering
\begin{gather*}\hspace*{10mm}
\begin{split}
&\raisebox{-4ex}{\hbox{
\begin{tikzpicture}[scale=0.5, >=stealth']
\tikzstyle{w}=[circle, draw, minimum size=4, inner sep=1]
\tikzstyle{b}=[circle, draw, fill, minimum size=4, inner sep=1]
\node [ext] (b1) at (0,0) {1};
\node [ext] (b2) at (2,0) {2};
\node [ext] (b3) at (1,-1.73) {3};
\draw[black,->,>=latex,>=latex][rightred] (b1) to (b2);
\draw (1,0.5) node[anchor=center] {{\small $i$}};
\draw (2,-0.85) node[anchor=center] {{\small $ii$}};
\draw (-0.08,-0.9) node[anchor=center] {{\small $iii$}};
\draw[black,->,>=latex][decoration={markings,mark=at position .2 with {\arrowreversed[scale=1.2,linkcolor]{latex}}},
 postaction={decorate}] (b2) to (b3);
\draw[black,->,>=latex,>=latex][rightred] (b3) to (b1);
\end{tikzpicture}}}
\quad\circ_2\quad
\raisebox{-0.8ex}{\hbox{
\begin{tikzpicture}[scale=0.5, >=stealth']
\tikzstyle{w}=[circle, draw, minimum size=4, inner sep=1]
\tikzstyle{b}=[circle, draw, fill, minimum size=4, inner sep=1]
\node [ext] (b4) at (2,0) {1};
\node [ext] (b5) at (4,0) {2};
\draw (3,0.5) node[anchor=center] {{\small $i$}};
\draw[black,->][decoration={markings,mark=at position .18 with {\arrowreversed[scale=1.2,linkcolor]{latex}}},
 postaction={decorate}] (b4) to (b5);
\end{tikzpicture}}}
\quad=\quad
\raisebox{-5ex}{\hbox{
 \begin{tikzpicture}[scale=0.7,yscale=-1]
 \node [ext] (b1) at (0,0) {1};
\draw (2,0) node[draw, color=gray, circle, minimum size=35, inner sep=1] (b2) {~};
\node [ext] (b3) at (1,1.73) {4};
\draw[black,->,>=latex,>=latex][rightred] (b1) to (b2);
\draw (0.8,-0.3) node[anchor=center] {{\small $i$}};
\draw (2,-0.3) node[anchor=center] {{\small $iv$}};
\draw (0.1,1) node[anchor=center] {{\small $iii$}};
\draw[black,->,>=latex,>=latex][leftred] (b2) to (b3);
\draw[black,->,>=latex,>=latex][rightred] (b3) to (b1);
 \node [ext] (b4) at (1.5,0) {2};
\node [ext] (b5) at (2.6,0) {3};
\draw (1.7,1) node[anchor=center] {{\small $ii$}};
\draw[black,->,>=latex,>=latex][decoration={markings,mark=at position .02 with {\arrowreversed[scale=1.2,linkcolor]{latex}}},
 postaction={decorate}] (b4) to (b5);
 \end{tikzpicture}}}\\
&\qquad =
\raisebox{-5ex}{\hbox{
\begin{tikzpicture}[scale=0.5, >=stealth']
\tikzstyle{w}=[circle, draw, minimum size=4, inner sep=1]
\tikzstyle{b}=[circle, draw, fill, minimum size=4, inner sep=1]
\node [ext] (b4) at (0,0) {4};
\node [ext] (b3) at (2,0) {3};
\node [ext] (b2) at (2,2) {2};
\node [ext] (b1) at (0,2) {1};
\draw (1,2.5) node[anchor=center] {{\small $i$}};
\draw (1,-0.5) node[anchor=center] {{\small $ii$}};
\draw (-0.5,1) node[anchor=center] {{\small $iii$}};
\draw (2.5,1) node[anchor=center] {{\small $iv$}};
\draw[black,->,>=latex,>=latex][rightred] (b1) to (b2);
\draw[black,->,>=latex,>=latex][leftred] (b2) to (b3);
\draw[black,->,>=latex,>=latex][leftred] (b3) to (b4);
\draw[black,->,>=latex,>=latex][rightred] (b4) to (b1);
\end{tikzpicture}}}
+
\raisebox{-5ex}{\hbox{
\begin{tikzpicture}[scale=0.5, >=stealth']
\tikzstyle{w}=[circle, draw, minimum size=4, inner sep=1]
\tikzstyle{b}=[circle, draw, fill, minimum size=4, inner sep=1]
\node [ext] (b4) at (0,0) {4};
\node [ext] (b3) at (2,0) {2};
\node [ext] (b2) at (2,2) {3};
\node [ext] (b1) at (0,2) {1};
\draw (1,2.5) node[anchor=center] {{\small $i$}};
\draw (1,-0.5) node[anchor=center] {{\small $ii$}};
\draw (-0.5,1) node[anchor=center] {{\small $iii$}};
\draw (2.5,1) node[anchor=center] {{\small $iv$}};
\draw[black,->,>=latex,>=latex][rightred] (b1) to (b2);
\draw[black,->,>=latex,>=latex][leftred] (b3) to (b2);
\draw[black,->,>=latex,>=latex][leftred] (b3) to (b4);
\draw[black,->,>=latex,>=latex][rightred] (b4) to (b1);
\end{tikzpicture}}}
+
\raisebox{-5ex}{\hbox{
\begin{tikzpicture}[scale=0.5, >=stealth']
\tikzstyle{w}=[circle, draw, minimum size=4, inner sep=1]
\tikzstyle{b}=[circle, draw, fill, minimum size=4, inner sep=1]
\node [ext] (b1) at (0,0) {1};
\node [ext] (b2) at (2,0) {2};
\node [ext] (b3) at (1,-1.73) {4};
\node [ext] (b4) at (3.5,1.3) {3};
\draw[black,->,>=latex,>=latex][rightred] (b1) to (b2);
\draw (1,0.5) node[anchor=center] {{\small $i$}};
\draw (2,-0.85) node[anchor=center] {{\small $ii$}};
\draw (-0.08,-0.9) node[anchor=center] {{\small $iii$}};
\draw (2.4,1.25) node[anchor=center] {{\small $iv$}};
\draw[black,->,>=latex,>=latex][leftred] (b2) to (b3);
\draw[black,->,>=latex,>=latex][rightred] (b3) to (b1);
\draw[black,->,>=latex,>=latex][decoration={markings,mark=at position .1 with {\arrowreversed[scale=1.2,linkcolor]{latex}}}, postaction={decorate}, shorten >=0.4pt] (b2) to (b4);
\end{tikzpicture}}}
+
\raisebox{-5ex}{\hbox{
\begin{tikzpicture}[scale=0.5, >=stealth']
\tikzstyle{w}=[circle, draw, minimum size=4, inner sep=1]
\tikzstyle{b}=[circle, draw, fill, minimum size=4, inner sep=1]
\node [ext] (b1) at (0,0) {1};
\node [ext] (b2) at (2,0) {3};
\node [ext] (b3) at (1,-1.73) {4};
\node [ext] (b4) at (3.5,1.3) {2};
\draw[black,->,>=latex,>=latex][rightred] (b1) to (b2);
\draw (1,0.5) node[anchor=center] {{\small $i$}};
\draw (2,-0.85) node[anchor=center] {{\small $ii$}};
\draw (-0.08,-0.9) node[anchor=center] {{\small $iii$}};
\draw (2.4,1.25) node[anchor=center] {{\small $iv$}};
\draw[black,->,>=latex,>=latex][leftred] (b2) to (b3);
\draw[black,->,>=latex,>=latex][rightred] (b3) to (b1);
\draw[black,->,>=latex,>=latex][decoration={markings,mark=at position .1 with {\arrowreversed[scale=1.2,linkcolor]{latex}}}, postaction={decorate}, shorten >=0.4pt] (b4) to (b2);
\end{tikzpicture}}}
\end{split}
\end{gather*}
 \caption{Example of partial composition $\circ_2\colon \dGra{d}{2}(3)\otimes\dGra{d}{2}(2)\to\dGra{d}{2}(4)$.}
 \label{figtypicaloddmulti}
\end{figure}

There is a natural sequence of embeddings of operads\footnote{Here, $\Graplain_{d}\equiv\dGra{d}{0}$ (resp.\ $\dGraplain_{d}\equiv\dGra{d}{1}$) stands for the operad of one-colored undirected (resp.\ directed) graphs.}
\begin{gather}
\Graplain_{d}\ \overset{\Or}{{\longhookrightarrow}}\ \dGraplain_{d}\ \overset{\Or}{{\longhookrightarrow}}\ \dGra{d}{2}\ \overset{\Or}{{\longhookrightarrow}}\ \dGra{d}{3}\ \overset{\Or}{{\longhookrightarrow}}\ \cdots\label{ormorph}
\end{gather}
given by mapping each graph in $\dGra{\dime}{c}$ to a sum of graphs in $\dGra{d}{c+1}$ where the summation runs over all the possible ways to orient the arrow of the additional direction. We call such mapping the \textit{orientation morphism} $\Or\colon \dGra{\dime}{c}\longhookrightarrow\dGra{d}{c+1}$ (cf.\ Figure~\ref{figorientmorph} for an example).
\begin{figure}[ht]\centering
$$
\Or\bigg(\raisebox{-4.8ex}{\hbox{\begin{tikzpicture}[scale=0.5, >=stealth']
\tikzstyle{b}=[circle, draw, fill, minimum size=2, inner sep=0.02]
\node [ext] (b2) at (0,0) {1};
\node [ext] (b3) at (2,0) {2};
\draw (1,1.2) node (b4) {$i$};
\draw (1,-1.2) node (b4) {$ii$};
\draw[black,->,>=latex] (b2) to[out=65, in=125, looseness=1.] (b3);
\draw[black,->,>=latex] (b2) to[out=-65, in=-125, looseness=1.] (b3);
\end{tikzpicture}}}\bigg)=
\raisebox{-4.8ex}{\hbox{\begin{tikzpicture}[scale=0.5, >=stealth']
\tikzstyle{b}=[circle, draw, fill, minimum size=2, inner sep=0.02]
\node [ext] (b2) at (0,0) {1};
\node [ext] (b3) at (2,0) {2};
\draw (1,1.2) node (b4) {$i$};
\draw (1,-1.2) node (b4) {$ii$};
\draw[black,->,>=latex][decoration={markings,mark=at position .55 with {\arrow[scale=1.35,linkcolor,rotate=10,yshift=0.017cm]{latex}}}, postaction={decorate}, shorten >=0.4pt] (b2) to[out=65, in=125, looseness=1.] (b3);
\draw[black,->,>=latex][decoration={markings,mark=at position .55 with {\arrow[scale=1.35,linkcolor,rotate=-10,yshift=-0.017cm]{latex}}}, postaction={decorate},shorten >=0.4pt] (b2) to[out=-65, in=-125, looseness=1.] (b3);
\end{tikzpicture}}}
+\raisebox{-4.8ex}{\hbox{\begin{tikzpicture}[scale=0.5, >=stealth']
\tikzstyle{b}=[circle, draw, fill, minimum size=2, inner sep=0.02]
\node [ext] (b2) at (0,0) {1};
\node [ext] (b3) at (2,0) {2};
\draw (1,1.2) node (b4) {$i$};
\draw (1,-1.2) node (b4) {$ii$};
\draw[black,->,>=latex][decoration={markings,mark=at position .55 with {\arrow[scale=1.35,linkcolor,rotate=10,yshift=0.017cm]{latex}}}, postaction={decorate}, shorten >=0.4pt] (b2) to[out=65, in=125, looseness=1.] (b3);
\draw[black,->,>=latex][decoration={markings,mark=at position .35 with {\arrowreversed[scale=1.35,linkcolor,rotate=5,yshift=-0.007cm]{latex}}}, postaction={decorate},shorten >=0.4pt] (b2) to[out=-65, in=-125, looseness=1.] (b3);
\end{tikzpicture}}}
+\raisebox{-4.8ex}{\hbox{\begin{tikzpicture}[scale=0.5, >=stealth']
\tikzstyle{b}=[circle, draw, fill, minimum size=2, inner sep=0.02]
\node [ext] (b2) at (0,0) {1};
\node [ext] (b3) at (2,0) {2};
\draw (1,1.2) node (b4) {$i$};
\draw (1,-1.2) node (b4) {$ii$};
\draw[black,->,>=latex][decoration={markings,mark=at position .4 with {\arrowreversed[scale=1.35,linkcolor,rotate=4]{latex}}},
 postaction={decorate},
 shorten >=0.4pt] (b2) to[out=65, in=125, looseness=1.] (b3);
\draw[black,->,>=latex][decoration={markings,mark=at position .55 with {\arrow[scale=1.35,linkcolor,rotate=-10,yshift=-0.017cm]{latex}}}, postaction={decorate},shorten >=0.4pt] (b2) to[out=-65, in=-125, looseness=1.] (b3);
\end{tikzpicture}}}
+\raisebox{-4.8ex}{\hbox{\begin{tikzpicture}[scale=0.5, >=stealth']
\tikzstyle{b}=[circle, draw, fill, minimum size=2, inner sep=0.02]
\node [ext] (b2) at (0,0) {1};
\node [ext] (b3) at (2,0) {2};
\draw (1,1.2) node (b4) {$i$};
\draw (1,-1.2) node (b4) {$ii$};
\draw[black,->,>=latex][decoration={markings,mark=at position .4 with {\arrowreversed[scale=1.35,linkcolor,rotate=4]{latex}}},
 postaction={decorate},
 shorten >=0.4pt] (b2) to[out=65, in=125, looseness=1.] (b3);
\draw[black,->,>=latex][decoration={markings,mark=at position .4 with {\arrowreversed[scale=1.35,linkcolor,rotate=0]{latex}}},
 postaction={decorate},
 shorten >=0.4pt] (b2) to[out=-65, in=-125, looseness=1.] (b3);
\end{tikzpicture}}}
$$
 \caption{The orientation morphism $\Or\colon \dGraplain_{d}\longhookrightarrow\dGra{d}{2}$.} \label{figorientmorph}
\end{figure}

\noindent Denoting $\mLie\pset{1-d}$ the $(1-d)$-suspended Lie operad,\footnote{Recall that representations of $\mLie\pset{1-d}$ on a vector space $\alg$ are in bijective correspondence with Lie algebra structures on $\alg[1-d]$, hence the graded Lie bracket on $\alg$ has intrinsic degree $1-d$.} there is an operad morphism
\[
\gamma_0\colon \ \mLie\pset{1-d}\to\dGra{d}{c}
\]
 sending the generator
$\raisebox{-2ex}{\hbox{\begin{tikzpicture}[scale=0.8, >=stealth']
 \coordinate (1) at (0,0);
 \coordinate (2) at (135:0.5);
 \coordinate (3) at (45:0.5);
 \coordinate (4) at (-90:0.45);

 \draw (1) -- (2) ;
 \draw (1) -- (3) ;
 \draw (1) -- (4) ;
 \foreach \x in {(2), (3), (4)}{
 \fill \x circle[radius=3pt];
 }
\end{tikzpicture}}}\in\mLie\pset{1-d}(2)$ to the graph
\[
\raisebox{-0.5ex}{\hbox{\begin{tikzpicture}[scale=0.5, >=stealth']
\tikzstyle{b}=[circle, draw, fill, minimum size=2, inner sep=0.02]
\node [b] (b2) at (0,0) {1};
\node [b] (b3) at (2,0) {2};
\draw[black,->,>=latex][decoration={markings,mark=at position .6 with {\arrow[scale=1.35]{latex}}},
 postaction={decorate},
 shorten >=0.4pt] (b2) to (b3);
\end{tikzpicture}}}:=\raisebox{-0.7ex}{\hbox{\begin{tikzpicture}[scale=0.5, >=stealth']
\tikzstyle{b}=[circle, draw, fill, minimum size=2, inner sep=0.02]
\node [ext] (b2) at (0,0) {1};
\node [ext] (b3) at (2,0) {2};
\draw[black,->,>=latex][decoration={markings,mark=at position .6 with {\arrow[scale=1.35]{latex}}},
 postaction={decorate},
 shorten >=0.4pt] (b2) to (b3);
\end{tikzpicture}}}+(-1)^d\, \raisebox{-0.7ex}{\hbox{\begin{tikzpicture}[scale=0.5, >=stealth']
\tikzstyle{b}=[circle, draw, fill, minimum size=2, inner sep=0.02]
\node [ext] (b2) at (0,0) {2};
\node [ext] (b3) at (2,0) {1};
\draw[black,->,>=latex][decoration={markings,mark=at position .6 with {\arrow[scale=1.35]{latex}}},
 postaction={decorate},
 shorten >=0.4pt] (b2) to (b3);
\end{tikzpicture}}}
\]
 where the graph $\raisebox{-0.7ex}{\hbox{\begin{tikzpicture}[scale=0.5, >=stealth']
\tikzstyle{b}=[circle, draw, fill, minimum size=2, inner sep=0.02]
\node [ext] (b2) at (0,0) {1};
\node [ext] (b3) at (2,0) {2};
\draw[black,->,>=latex][decoration={markings,mark=at position .6 with {\arrow[scale=1.35]{latex}}},
 postaction={decorate},
 shorten >=0.4pt] (b2) to (b3);
\end{tikzpicture}}}\in\dGra{d}{c}(2)$ is obtained by decorating
$\raisebox{-0.7ex}{\hbox{\begin{tikzpicture}[scale=0.5, >=stealth']
\tikzstyle{b}=[circle, draw, fill, minimum size=2, inner sep=0.02]
\node [ext] (b2) at (0,0) {1};
\node [ext] (b3) at (1.8,0) {2};
\draw[black,->,>=latex] (b2) to (b3);
\end{tikzpicture}}}$ with $c-1$ additional colors and summing over all the possible orientations.\footnote{\label{footnotegraph}In other words, $\raisebox{-0.9ex}{\hbox{\begin{tikzpicture}[scale=0.5, >=stealth']
\tikzstyle{b}=[circle, draw, fill, minimum size=2, inner sep=0.02]
\node [ext] (b2) at (0,0) {1};
\node [ext] (b3) at (2,0) {2};
\draw[black,->,>=latex][decoration={markings,mark=at position .6 with {\arrow[scale=1.35]{latex}}},
 postaction={decorate},
 shorten >=0.4pt] (b2) to (b3);
\end{tikzpicture}}}:=(\Or)^{c-1}\big(\raisebox{-0.9ex}{\hbox{\begin{tikzpicture}[scale=0.5, >=stealth']
\tikzstyle{b}=[circle, draw, fill, minimum size=2, inner sep=0.02]
\node [ext] (b2) at (0,0) {1};
\node [ext] (b3) at (1.8,0) {2};
\draw[black,->,>=latex] (b2) to (b3);
\end{tikzpicture}}}\big)$. For example, $\raisebox{-0.9ex}{\hbox{\begin{tikzpicture}[scale=0.5, >=stealth']
\tikzstyle{b}=[circle, draw, fill, minimum size=2, inner sep=0.02]
\node [ext] (b2) at (0,0) {1};
\node [ext] (b3) at (2,0) {2};
\draw[black,->,>=latex][decoration={markings,mark=at position .6 with {\arrow[scale=1.35]{latex}}},
 postaction={decorate},
 shorten >=0.4pt] (b2) to (b3);
\end{tikzpicture}}}:=\raisebox{-0.9ex}{\hbox{\begin{tikzpicture}[scale=0.5, >=stealth']
\tikzstyle{b}=[circle, draw, fill, minimum size=2, inner sep=0.02]
\node [ext] (b2) at (0,0) {1};
\node [ext] (b3) at (2,0) {2};
\draw[black,->,>=latex][decoration={markings,mark=at position .6 with {\arrow[scale=1.35,linkcolor]{latex}}},
 postaction={decorate},
 shorten >=0.4pt] (b2) to (b3);
\end{tikzpicture}}}
+
\raisebox{-0.9ex}{\hbox{\begin{tikzpicture}[scale=0.5, >=stealth']
\tikzstyle{b}=[circle, draw, fill, minimum size=2, inner sep=0.02]
\node [ext] (b2) at (0,0) {1};
\node [ext] (b3) at (2,0) {2};
\draw[black,->,>=latex][decoration={markings,mark=at position .13 with {\arrowreversed[scale=1.35,linkcolor]{latex}}},
 postaction={decorate},
 shorten >=0.4pt] (b2) to (b3);
\end{tikzpicture}}}$ in $\dGra{d}{2}$.}

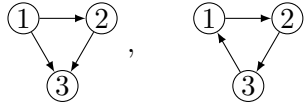
\begin{figure}[ht]\centering
\raisebox{-4ex}{\hbox{\begin{tikzpicture}[scale=0.4, >=stealth']
\tikzstyle{w}=[circle, draw, minimum size=4, inner sep=1]
\tikzstyle{b}=[circle, draw, fill, minimum size=2, inner sep=0.02]
\draw (30:1.5) node [ext] (b2) {2};
\draw (150:1.5) node [ext] (b1) {1};
\draw (-90:1.5) node [ext] (b3) {3};
\draw[black,->,>=latex] (b1) to (b2);
\draw[black,->,>=latex] (b1) to (b3);
\draw[black,->,>=latex] (b2) to (b3);
\end{tikzpicture}}}
,\qquad
\raisebox{-4ex}{\hbox{\begin{tikzpicture}[scale=0.4, >=stealth']
\tikzstyle{w}=[circle, draw, minimum size=4, inner sep=1]
\tikzstyle{b}=[circle, draw, fill, minimum size=2, inner sep=0.02]
\draw (30:1.5) node [ext] (b2) {2};
\draw (150:1.5) node [ext] (b1) {1};
\draw (-90:1.5) node [ext] (b3) {3};
\draw[black,->,>=latex] (b1) to (b2);
\draw[black,->,>=latex] (b2) to (b3);
\draw[black,->,>=latex] (b3) to (b1);
\end{tikzpicture}}}
 \caption{Example of oriented (left) and cyclic (right) graph in $\gra_{3,3}$.}
 \label{figgraphexori}
\end{figure}

A multidigraph in $\gra_{N,k}$ will be said \textit{oriented} (or acyclic) if it does not contain cycles.\footnote{Recall that a \textit{cycle} (or wheel) is a (non-trivial) directed path from a vertex to itself.} Contrariwise, it will be said non-oriented (or cyclic) if it contains at least one cycle, cf.\ Fi\-gure~\ref{figgraphexori}.
The subset of oriented multidigraphs will be denoted $\odgraplain_{N,k}\subset\gra_{N,k}$. This definition can be extended to multi-directed graphs by defining the subset\footnote{By convention, we identify $\odgra{0}{0}_{N,k}\equiv\ugra_{N,k}$, $\odgra{1}{0}_{N,k}\equiv\gra_{N,k}$ and $\odgra{0}{1}_{N,k}\equiv\odgraplain_{N,k}$.} $\odgra{j}{i}_{N,k}\subseteq\dgra{N,k}{c}$ of multi-directed graphs with $c=i+j$ directions for which there exists a subset of $i$ directions~-- black and/or colored~-- such that there are no cycles made of the corresponding arrows.\footnote{As an example, the graph of Figure \ref{figgraphexmult} belongs to $\odgra{2}{2}_{4,5}$ as it does not contain cycles for black and yellow arrows (although it does so for red and blue arrows). In the right-hand side of Figure \ref{figorientmorph}, the first and fourth graphs belong to $\odgra{0}{2}_{2,2}$ while the second and third graphs belong to $\odgra{1}{1}_{2,2}$. } Substituting $\odgra{j}{i}_{N,k}$ in place of $\dgra{N,k}{c}$ in \eqref{subop} allows to define the collection of $\mathbb S_N$-modules $\pset{\odGra{d}{j}{i}(N)}_{N\geq1}$ for all $i,j,d\geq0$. It is easy to check that the latter is closed under partial compositions and hence defines a suboperad $\odGra{d}{j}{i}\subseteq\dGra{d}{c}$ of multi-oriented graphs. Note that the graph $\raisebox{-0.5ex}{\hbox{\begin{tikzpicture}[scale=0.5, >=stealth']
\tikzstyle{b}=[circle, draw, fill, minimum size=2, inner sep=0.02]
\node [b] (b2) at (0,0) {1};
\node [b] (b3) at (2,0) {2};
\draw[black,->,>=latex][decoration={markings,mark=at position .6 with {\arrow[scale=1.35]{latex}}},
 postaction={decorate},
 shorten >=0.4pt] (b2) to (b3);
\end{tikzpicture}}}$ is (trivially) multi-oriented and hence defines a morphism of operads
\[
\gamma_0\colon \mLie\pset{1-d}\to\odGra{d}{j}{i} \qquad \text{for all $i,j,d\geq0$.}
\]

Oriented graphs belong to the larger subset of multidigraphs possessing a \textit{source}, i.e., a vertex admitting only outgoing arrows. More generally, the suboperad of multi-directed graphs with $c=|k|+j$ directions such that~$|k|$ directions are sourced will be denoted $\sodGra{d}{j}{k}\subseteq\dGra{d}{c}$ for all $k\in\mathbb Z$ and $d,j\geq0$, where negative values of~$k$ correspond to directions admitting a \textit{sink}, i.e., a vertex with only ingoing arrows. Finally, the suboperad of graphs such that $j$ directions admit at least one source and one sink will be denoted $\sodGra{d}{j}{i,-i}$, for all $i,j,d\geq0$. It is a~well-known result in graph theory that oriented graphs admit at least one source and one sink (see, e.g.,~\cite{Zivkovic2017} for a statement) so that we have a sequence of inclusions
\[
\odGra{d}{j}{|k|}\subseteq\sodGra{d}{j}{|k|,-|k|}\subseteq\sodGra{d}{j}{k} \qquad \text{for all $k\in\mathbb Z$ and $d,j\geq0$}.
\]

{\bf Graph complexes.} Given the multi-oriented graph operad $\odGra{d}{j}{i}$, one defines the dg Lie algebra of multi-oriented graphs $\odfGC{d}{j}{i}$ as the deformation complex
\[
\odfGC{d}{j}{i}:=\Def\big(\mLie\pset{1-d}\overset{\gamma_0}{\to}\odGra{d}{j}{i}\big)
\] of the morphism of operads~$\gamma_0$.\footnote{We refer to references \cite{Loday2012,Merkulov2007} for details.} As a graded vector space, $\odfGC{d}{j}{i}$ is defined as
\begin{gather}
\bullet \ \text{$\dime$ even:} \ \odfGC{d}{j}{i}:=\bigoplus_{N\geq1}\big(\odGra{d}{j}{i}(N)[\dime(1-N)]\big)^{\mS_N},\label{defGC1}\\
\bullet \ \text{$\dime$ odd:} \ \odfGC{d}{j}{i}:=\bigoplus_{N\geq1}\big(\odGra{d}{j}{i}(N)\otimes \sgn_N[\dime(1-N)]\big)^{\mS_N},\label{defGC2}
\end{gather}
where the terms between brackets denote degree suspension while the superscript stands for taking invariants with respect to the right action of $\mS_N$, with $\sgn_N$ the 1-dimensional signature representation of $\mS_N$. In other words, vertices are bosonic for $d$ even and fermionic for $d$ odd. According to the degree suspension in \eqref{subop} and \eqref{defGC1}--\eqref{defGC2}, the degree of an element\footnote{Note that the graph degree is insensitive to the number $i$ of oriented directions. } $\Gamgraph\in\odfGC{d}{j}{i}$ with $N$ vertices and $k$ edges is given by $|\Gamgraph|=\dime(N-1)+k(1-\dime)$.

 The graded Lie bracket $\brdot$ of degree 0 on $\odfGC{d}{j}{i}$ is defined as usual in terms of the partial composition operations (cf., e.g., \cite[Section~4.2]{Morand2019a} for details). The differential $\delta:=\br{\Upsilon_\mathsf{S}}{\cdot}$ is defined by taking the adjoint action with respect to the Maurer--Cartan element\footnote{For example, the graph
$\Upsilon_\mathsf{S}:=\raisebox{-0.5ex}{\hbox{\begin{tikzpicture}[scale=0.5, >=stealth']
\tikzstyle{b}=[circle, draw, fill, minimum size=2, inner sep=0.02]
\node [b] (b2) at (0,0) {2};
\node [b] (b3) at (1.8,0) {3};
\draw[black,->,>=latex][decoration={markings,mark=at position .6 with {\arrow[scale=1.35,linkcolor]{latex}}},
 postaction={decorate},
 shorten >=0.4pt] (b2) to (b3);
\end{tikzpicture}}}=\raisebox{-0.7ex}{\hbox{\begin{tikzpicture}[scale=0.5, >=stealth']
\tikzstyle{b}=[circle, draw, fill, minimum size=2, inner sep=0.02]
\node [ext] (b2) at (0,0) {1};
\node [ext] (b3) at (2,0) {2};
\draw[black,->,>=latex][decoration={markings,mark=at position .6 with {\arrow[scale=1.35,linkcolor]{latex}}},
 postaction={decorate},
 shorten >=0.4pt] (b2) to (b3);
\end{tikzpicture}}}
+
\raisebox{-0.7ex}{\hbox{\begin{tikzpicture}[scale=0.5, >=stealth']
\tikzstyle{b}=[circle, draw, fill, minimum size=2, inner sep=0.02]
\node [ext] (b2) at (0,0) {1};
\node [ext] (b3) at (2,0) {2};
\draw[black,->,>=latex][decoration={markings,mark=at position .13 with {\arrowreversed[scale=1.35,linkcolor]{latex}}},
 postaction={decorate},
 shorten >=0.4pt] (b2) to (b3);
\end{tikzpicture}}}
+(-1)^d\, \big(
\raisebox{-0.9ex}{\hbox{\begin{tikzpicture}[scale=0.5, >=stealth']
\tikzstyle{b}=[circle, draw, fill, minimum size=2, inner sep=0.02]
\node [ext] (b2) at (0,0) {2};
\node [ext] (b3) at (2,0) {1};
\draw[black,->,>=latex][decoration={markings,mark=at position .6 with {\arrow[scale=1.35,linkcolor]{latex}}},
 postaction={decorate},
 shorten >=0.4pt] (b2) to (b3);
\end{tikzpicture}}}
+
\raisebox{-0.9ex}{\hbox{\begin{tikzpicture}[scale=0.5, >=stealth']
\tikzstyle{b}=[circle, draw, fill, minimum size=2, inner sep=0.02]
\node [ext] (b2) at (0,0) {2};
\node [ext] (b3) at (2,0) {1};
\draw[black,->,>=latex][decoration={markings,mark=at position .13 with {\arrowreversed[scale=1.35,linkcolor]{latex}}},
 postaction={decorate},
 shorten >=0.4pt] (b2) to (b3);
\end{tikzpicture}}}
\big)$ is a Maurer--Cartan element in $\odfGC{d}{j}{i}$ with $i,j,d\geq0$ and $i+j=2$.}
\[ \Upsilon_\mathsf{S}:=\raisebox{-0.5ex}{\hbox{\begin{tikzpicture}[scale=0.5, >=stealth']
\tikzstyle{b}=[circle, draw, fill, minimum size=2, inner sep=0.02]
\node [b] (b2) at (0,0) {1};
\node [b] (b3) at (1.8,0) {2};
\draw[black,->,>=latex][decoration={markings,mark=at position .6 with {\arrow[scale=1.35]{latex}}},
 postaction={decorate},
 shorten >=0.4pt] (b2) to (b3);
\end{tikzpicture}}}=\raisebox{-0.7ex}{\hbox{\begin{tikzpicture}[scale=0.5, >=stealth']
\tikzstyle{b}=[circle, draw, fill, minimum size=2, inner sep=0.02]
\node [ext] (b2) at (0,0) {1};
\node [ext] (b3) at (2,0) {2};
\draw[black,->,>=latex][decoration={markings,mark=at position .6 with {\arrow[scale=1.35]{latex}}},
 postaction={decorate},
 shorten >=0.4pt] (b2) to (b3);
\end{tikzpicture}}}+(-1)^d\, \raisebox{-0.7ex}{\hbox{\begin{tikzpicture}[scale=0.5, >=stealth']
\tikzstyle{b}=[circle, draw, fill, minimum size=2, inner sep=0.02]
\node [ext] (b2) at (0,0) {2};
\node [ext] (b3) at (2,0) {1};
\draw[black,->,>=latex][decoration={markings,mark=at position .6 with {\arrow[scale=1.35]{latex}}},
 postaction={decorate},
 shorten >=0.4pt] (b2) to (b3);
\end{tikzpicture}}}.
\]

Note to conclude that the dg Lie algebra $\odfGC{d}{j}{|k|}$ of multi-oriented graphs is a sub-dg Lie algebra of the dg Lie algebra of multi-sourced/sinked graphs $\sodfGC{d}{j}{k}$ defined for all $k\in\mathbb Z$ and $d,j\geq0$ by\footnote{Equivalently, one can define $\sodfGC{d}{j}{k}:=\Def\big(\mLie\pset{1-d}\overset{\gamma_0}{\to}\sodGra{d}{j}{k}\big)$.} substituting $\sodGra{d}{j}{k}(N)$ in place of $\odGra{d}{j}{k}(N)$ in \eqref{defGC1}--\eqref{defGC2}.

The various graph operads and complexes discussed in the present section are summarised in Table~\ref{Tabgraphs}.
\begin{table}[h]\centering\small
 \caption{Summary of graph operads and complexes (and their connected variants) in dimension $d$.} \label{Tabgraphs}
 \vspace{1mm}

\begin{tabular}[5pt]{ccccc}
 \toprule
 \multicolumn{5}{c}{\rule[-0.3cm]{0cm}{0.8cm}Undirected graphs with one color (Kontsevich's graphs)}\\
 \rule[-0.3cm]{0cm}{0.8cm} $\Graplain_d$&$\cGra_d$&$\fGC_d$&$\fcGC_d$&$d\geq0$\\
 \cline{1-5}
 \multicolumn{5}{c}{\rule[-0.3cm]{0cm}{0.8cm}Multi-directed graphs with $c$ colors}\\
 \rule[-0.3cm]{0cm}{0.8cm} $\dGra{d}{c}$&$\dcGra{d}{c}$&$\dfGC{d}{c}$&$\dfcGC{d}{c}$&$c,d\geq0$\\
 \cline{1-5}
 \multicolumn{5}{c}{\rule[-0.3cm]{0cm}{0.8cm}Multi-oriented graphs with $c=i+j$ colors and $i$ oriented directions}\\
 \rule[-0.3cm]{0cm}{0.8cm} $\odGra{d}{j}{i}$&$\odcGra{d}{j}{i}$&$\odfGC{d}{j}{i}$&$\odfcGC{d}{j}{i}$&$i,j,d\geq0$\\
 \cline{1-5}
 \multicolumn{5}{c}{\rule[-0.3cm]{0cm}{0.8cm}Multi-sourced/sinked graphs with $c=|k|+j$ colors and $|k|$ sourced/sinked directions}\\
 \rule[-0.3cm]{0cm}{0.8cm} $\sodGra{d}{j}{k}$&$\sodcGra{d}{j}{k}$&$\sodfGC{d}{j}{k}$&$\sodfcGC{d}{j}{k}$&$j,d\geq0$, $k\in\mathbb Z$\\
 \cline{1-5}
 \multicolumn{5}{c}{\rule[-0.3cm]{0cm}{0.8cm}Multi-sourced/sinked graphs with $c=i+j$ colors and $i$ directions being both sourced and sinked}\\
 \rule[-0.3cm]{0cm}{0.8cm} $\sodGra{d}{j}{i,-i}$&$\sodcGra{d}{j}{i,-i}$&$\sodfGC{d}{j}{i,-i}$&$\sodfcGC{d}{j}{i,-i}$&$i,j,d\geq0$\\
 \bottomrule
 \end{tabular}

\end{table}

\subsection{Cohomology}\label{section:Cohomology}
Having introduced the graph complex of multi-directed graphs as well as its sub-complexes of multi-sourced and multi-oriented graphs, we conclude this review section by collecting some known facts about their respective cohomology.\footnote{See \cite[Section~5]{Merkulov2017} and \cite[Section~7]{Merkulov2019} for reviews.} We start by introducing the suboperad $\odcGra{d}{j}{i}\subset\odGra{d}{j}{i}$ spanned by {\it connected} graphs, which in turns yields a sub-dg Lie algebra $\odfcGC{d}{j}{i}\subset\odfGC{d}{j}{i}$.\footnote{Similarly, one can introduce the suboperads of connected multi-directed graphs $\dcGra{d}{j}$ and connected multi-sourced graphs $\sodcGra{d}{j}{i}$ as well as their corresponding sub-dg Lie algebras $\dfcGC{d}{j}$ and $\sodfcGC{d}{j}{i}$. } This is justified by the fact that the cohomology of the graph complex $\odGra{d}{j}{i}$ is captured by its connected part,\footnote{We refer to~\cite{Willwacher2015} and \cite{Willwacher2015c} for the cases $i=0,1$ respectively and to \cite{Merkulov2017} for a general statement.} as follows from
\begin{displaymath}
H^\bullet(\odfGC{d}{j}{i})=\widehat{\mathcal O}\big(H^\bullet(\odfcGC{d}{j}{i})[-d]\big)[d],
\end{displaymath} where $\widehat{\mathcal O}(\alg)$ denotes the completed symmetric algebra associated with the graded vector space~$\alg$. As far as cohomology is concerned, one can therefore restrict the analysis to the connected part of the above complexes. Moreover, a further simplification comes from the fact that the cohomology of the various complexes previously introduced can be related to one another. This is embodied by the following important theorem, due to T.~Willwacher\footnote{See \cite[Appendix~K]{Willwacher2015} for the first item (cf.\ also~\cite{Dolgushev2017}), \cite{Willwacher2015c} for the second and~\cite{Willwacher2015b} for the third.} for the case $i=0$ and to M.~{{\v{Z}}ivkovi{\'c}} \cite{Zivkovic2017a,Zivkovic2017} for its generalisation to arbitrary $i\geq0$.

\begin{Theorem}\label{thmWM}
For all integers $i,j,d\geq0$ and $k\in\mathbb Z$:
\begin{enumerate}\itemsep=0pt
\item[$1.$] The inclusion $\odfcGC{d}{j}{i}\longhookrightarrow \odfcGC{d}{j+1}{i}$ is a quasi-isomorphism.
\item[$2.$] There is a quasi-isomorphism $\odfcGC{d}{j}{i}\longrightarrow\odfcGC{d+1}{j}{i+1}$.
 \item[$3.$] The inclusion $\odfcGC{d}{j}{|k|}\longhookrightarrow\sodfcGC{d}{j}{k}$ is a quasi-isomorphism.
\end{enumerate}
The quasi-isomorphism of the second item preserves the additional grading provided by the first Betti number.\footnote{The \textit{first Betti number} is defined as $b:=k-N+1$ for a connected graph with $N$ vertices and $k$ edges.}
\end{Theorem}
A few comments are in order. We start by noting that the sequence of embeddings of operads~\eqref{ormorph} induces a sequence\footnote{Here, the dg Lie algebra $\fcGC_d\equiv \dfcGC{d}{0}$ (resp.~$\dfcGCplain_d\equiv \dfcGC{d}{1}$) stands for the usual Kontsevich graph complex of connected undirected (resp. directed) graphs.} of injective quasi-isomorphisms of complexes\footnote{See \cite{Dolgushev2017, Willwacher2015} for the first arrow and~\cite{Merkulov2017} for a general statement.}
\begin{displaymath}
\fcGC_d\underset{\Or}{\overset{\sim}{\longhookrightarrow}} \dfcGCplain_d\underset{\Or}{\overset{\sim}{\longhookrightarrow}} \dfcGC{d}{2}\ \underset{\Or}{\overset{\sim}{\longhookrightarrow}}\ \dfcGC{d}{3}\ \underset{\Or}{\overset{\sim}{\longhookrightarrow}}\ \cdots.
\end{displaymath}
Hence, adding extra colored direction does not affect the cohomology. The first item of Theo\-rem~\ref{thmWM} asserts that this result generalises to multi-oriented graphs where the number $i$ of oriented directions is kept fixed.

The situation gets more interesting when orienting extra directions since adding oriented directions \textit{does} change the cohomology. More precisely, the second item of the above theorem relates the cohomology of a given multi-oriented graph complex to the one of a (less oriented) complex in higher dimension. We will comment more on this important result in the next paragraph. Before doing so, let us note that the third item asserts that sourcing directions also affects the cohomology, but in a way that is completely captured by the cohomology of the (smaller) multi-oriented graph complex. Hence, the computation of the cohomology of the multi-sourced/sinked graph complex can always be reduced to computing the cohomology of the multi-oriented graph complex. In the next section, we will therefore express our results in terms of the smaller multi-oriented graph complex when possible.

For later use, Table~\ref{Tabcohom} summarises the cohomology in degrees $0$, $1$ and $2$ of the (undirected) graph complex in dimensions $d=1,2,3$.\footnote{The cocycle $L_3$ stands for the triangle loop $\raisebox{-3ex}{\hbox{\begin{tikzpicture}[scale=0.3, >=stealth']
\tikzstyle{w}=[circle, draw, minimum size=4, inner sep=1]
\tikzstyle{b}=[circle, draw, fill, minimum size=2, inner sep=0.02]
\draw (30:1.5) node [b] (b2) {2};
\draw (150:1.5) node [b] (b1) {1};
\draw (-90:1.5) node [b] (b3) {3};
\draw[black,line width=0.2mm,scale=1.35] (b1) to (b2);
\draw[black,line width=0.2mm,scale=1.35] (b1) to (b3);
\draw[black,line width=0.2mm,scale=1.35] (b2) to (b3);
\end{tikzpicture}}}$. Regarding the $\Theta$-cocycle and the Grothendieck--Teichm\"uller algebra $\grt_1$, see below.}
\begin{table}[h]\centering\small
 \caption{Cohomology in low dimension and degree.} \label{Tabcohom}

 \vspace{1mm}

 \begin{tabular}[5pt]{cccc}
 \toprule
 \rule[-0.2cm]{0cm}{0.6cm} &$H^0(\fcGC_d)$&$H^1(\fcGC_d)$&$H^2(\fcGC_d)$\\
 \bottomrule
 \rule[-0.2cm]{0cm}{0.6cm} $d=1$&\textbf{0}&$\corps\dl \Theta\dr$&$\corps\dl L_3\dr$\\
 \rule[-0.2cm]{0cm}{0.6cm} $d=2$&$\grt_1$&\hspace{0.4mm} \textbf{0}?&?\\
 \rule[-0.2cm]{0cm}{0.6cm} $d=3$&$\corps\dl L_3\dr$&\textbf{0}&\textbf{0}\\
 \bottomrule
 \end{tabular}
\end{table}

{\bf Climbing the dimension ladder.} Informally, the second item of Theorem~\ref{thmWM} allows to map familiar structures in low dimensions\footnote{See footnote \ref{footdim} for the interpretation of the dimension~$d$.} to novel incarnations thereof in higher dimensions.\footnote{Since the quasi-isomorphism of the second item of Theorem~\ref{thmWM} preserves both the graph degree and the first Betti number, a connected cocycle $\graphGam\in\gra_{N,k}$ in dimension $d$ is mapped to a cocycle $\graphGam'\in\gra_{N',k'}$ in dimension $d+1$ with $N'=k+1$ and $k'=2 k-N+1$. More generally, the incarnation of a cocycle $\graphGam\in\gra_{N,k}$ of $H^\bullet(\odfcGC{d}{j}{i})$ in dimension $d'>d$ is a graph $\graphGam'\in\gra_{N',k'}$ with $N'=N+(d'-d) b$ and $k'=k+(d'-d) b$ with $b$ the first Betti number.\label{footiso}} The most striking example of such hierarchy of structures stems from another important theorem of T.~Willwacher~\cite{Willwacher2015} showing the existence of an isomorphism of Lie algebras $H^0(\fcGC_2)\simeq\grt_1$ where $\grt_1$ denotes the infinite-dimensional Grothendieck--Teichm\"uller algebra. Combined with the second item of Theorem~\ref{thmWM}, we obtain a sequence $H^0(\fcGC_2)\simeq H^0(\odfcGC{3}{\bullet}{1})\simeq H^0(\odfcGC{4}{\bullet}{2})\simeq\cdots\simeq\grt_1$ of incarnations of $\grt_1$ in arbitrary dimensions $d\geq 2$.\footnote{For example, the tetrahedron graph $t_2\in H^0(\fcGC_2)$ is mapped to an oriented cocycle $t_3\in H^0(\odfcGC{3}{\bullet}{1})$ with $N=7$ vertices and $k=9$ edges (see footnote \ref{footiso}).} Different incarnations yield different actions of the Grothendieck--Teichm\"uller group\footnote{Recall that the Gro\-then\-dieck--Teichm\"uller group is defined by exponentiation of the pro-nilpotent Gro\-then\-dieck--Teichm\"uller algebra, i.e., $\GRT_1\equiv\exp(\grt_1)$, see~\cite{WillwachernotesGRT} for a review.} on various algebro-geometric structures. In its original incarnation via directed cocycles in $H^0(\dfGC{d}{1})$, the Grothendieck--Teichm\"uller group naturally acts via \Lieinf-automorphisms on the Schouten Lie algebra of polyvector fields $\Tpoly$ on (finite-dimensional) manifolds. This yields an action on the space of universal formality maps hence on the space of universal quantization maps for finite-dimensional Poisson manifolds \cite{Dolgushev2011,Jost2012, Kontsevich1997,Kontsevich:1997vb,Merkulov2008,Willwacher2015}. In dimension~3, there is a~natural action of $\GRT_1$ via oriented cocycles in $ H^0(\odfcGC{3}{0}{1})$ on the deformation complex of Lie bialgebras\footnote{As reviewed in Section~\ref{section:Lie bialgebras2}.} hence on the space of universal formality maps related to the deformation quantization of Lie bialgebras \cite{Merkulov2016a,Merkulov2016, Willwacher2015c}.\footnote{The fact that quantization of Lie bialgebras involves oriented graphs was already recognised in~\cite{Drinfeld1992}.}

Another important result for our story concerns the manifold incarnations of the $\Theta$-cocycle $\raisebox{-2.ex}{\hbox{\begin{tikzpicture}[
my angle/.style={draw, <->, angle eccentricity=1.3, angle radius=9mm},scale=0.4
 ]
\coordinate (O) at (0,0);
\coordinate[] (A) at (180:1);
\coordinate[] (B) at (0:1);
\draw[thick] (A) -- (B);
\draw[thick] (0,0) circle (1cm);
 \end{tikzpicture}}}\in\gra_{2,3}$ spanning the first cohomology class of $\fcGC_{1}$, i.e., $H^1(\fcGC_{1})\simeq\corps\dl\Theta\dr$. Again, applying the second item of Theorem \ref{thmWM} results in a sequence of isomorphisms
 \[
 H^1(\fcGC_{1})\simeq H^1(\odfcGC{2}{\bullet}{1})\simeq H^1(\odfcGC{3}{\bullet}{2})\simeq\cdots\simeq\corps.
 \]
 In its original incarnation as the cocycle spanning $H^1(\fcGC_{1})$, the $\Theta$-graph can be recursively extended\footnote{\label{footnote:recur}Potential obstructions in promoting the $\Theta$-graph to a full Maurer--Cartan element lie in $H^2(\fcGC_{1})\simeq\corps\dl L_{3}\dr$, cf.\ Table~\ref{Tabcohom}. Since the obstruction to the prolongation of the $\Theta$-graph at order $k\geq2$ has Betti number $k+2$, it never hits the loop graph $L_3$ (of Betti number~$1$) so that the prolongation is unobstructed at all orders and can be performed recursively. The argument carries identically to incarnations of the $\Theta$-graph in higher dimensions. } to a non-trivial Maurer--Cartan element $\Upsilon_\Theta:=\raisebox{-0.35ex}{\hbox{\begin{tikzpicture}[scale=0.4, >=stealth']
\tikzstyle{b}=[circle, draw, fill, minimum size=2, inner sep=0.02]
\node [b] (b2) at (0,0) {$\bullet$};
\node [b] (b3) at (2,0) {$\bullet$};
\draw[black,->,>=latex] (b2) to (b3);
\end{tikzpicture}}}+\raisebox{-2.ex}{\hbox{\begin{tikzpicture}[
my angle/.style={draw, <->, angle eccentricity=1.3, angle radius=9mm},scale=0.4
 ]
\coordinate (O) at (0,0);
\coordinate[] (A) at (180:1);
\coordinate[] (B) at (0:1);
\draw[thick] (A) -- (B);
\draw[thick] (0,0) circle (1cm);
 \end{tikzpicture}}}+\cdots$ in the graded Lie algebra $\big(\fcGC_{1},\brdot\big)$, see, e.g., \cite{Khoroshkin2014}. This Maurer--Cartan element is mapped to the Moyal star-commutator \cite{Groenewold1946,Moyal1949} via the natural action of $\fcGC_{1}$ on the algebra of functions of a (bosonic) symplectic manifold. Considering the case $d=2$ yields another important incarnation of the $\Theta$-graph as the oriented Kontsevich--Shoikhet cocycle~-- denoted $\Theta_2$ in the following, cf.\ Figure~\ref{figKS}~-- and spanning $H^1(\odfcGC{2}{\bullet}{1})$. The latter first appeared implicitly in \cite{Penkava1998} (see also \cite{Dito2013,Willwacher2013a}) as the obstruction to the existence of a cycle-less universal quantization of Poisson manifolds beyond order $\hbar^3$. It then appeared explicitly in \cite{Shoikhet2008a} as the obstruction to formality in infinite dimension while its graph theoretical interpretation~-- as the avatar of the $\Theta$-cocycle in dimension~2~-- has been elucidated in~\cite{Willwacher2015c}. The corresponding Maurer--Cartan element $\Upsilon_{\Theta_2}$ induces an exotic (and essentially unique) universal \Lieinf-structure on polyvector fields deforming non-trivially\footnote{For finite-dimensional manifolds, this exotic \Lieinf-structure can be shown to be isomorphic to the standard Schouten bracket, although in a highly non-trivial way, cf.~\cite{Merkulov2016} for explicit transcendental formulae.} the Schouten algebra on infinite-dimensional manifolds~\cite{Shoikhet2008a}. The latter can then be considered as the avatar in $d=2$ of the Moyal star-commutator in $d=1$.
 \begin{figure}[ht]\centering
$\raisebox{-5ex}{\hbox{\begin{tikzpicture}[scale=0.5, >=stealth']
\tikzstyle{w}=[circle, draw, minimum size=4, inner sep=1]
\tikzstyle{b}=[circle, draw, fill, minimum size=2, inner sep=0.02]
\node [b] (b2) at (30:2) {2};
\node [b] (b1) at (150:2) {1};
\node [b] (b3) at (-90:2) {3};
\node [b] (b4) at (0,0) {4};
\draw[black,->,>=latex] (b3) to (b1);
\draw[black,->,>=latex] (b4) to (b1);
\draw[black,->,>=latex] (b3) to (b2);
\draw[black,->,>=latex] (b4) to (b2);
\draw[black,->,>=latex] (b3) to (b4);
\end{tikzpicture}}}-2\, \raisebox{-5ex}{\hbox{\begin{tikzpicture}[scale=0.5, >=stealth']
\tikzstyle{w}=[circle, draw, minimum size=4, inner sep=1]
\tikzstyle{b}=[circle, draw, fill, minimum size=2, inner sep=0.02]
\node [b] (b2) at (30:2) {2};
\node [b] (b1) at (150:2) {1};
\node [b] (b3) at (-90:2) {3};
\node [b] (b4) at (0,0) {4};
\draw[black,->,>=latex] (b3) to (b1);
\draw[black,->,>=latex] (b4) to (b1);
\draw[black,->,>=latex] (b2) to (b3);
\draw[black,->,>=latex] (b2) to (b4);
\draw[black,->,>=latex] (b3) to (b4);
\end{tikzpicture}}}+\raisebox{-5ex}{\hbox{\begin{tikzpicture}[scale=0.5, >=stealth']
\tikzstyle{w}=[circle, draw, minimum size=4, inner sep=1]
\tikzstyle{b}=[circle, draw, fill, minimum size=2, inner sep=0.02]
\node [b] (b2) at (30:2) {2};
\node [b] (b1) at (150:2) {1};
\node [b] (b3) at (-90:2) {3};
\node [b] (b4) at (0,0) {4};
\draw[black,->,>=latex] (b1) to (b3);
\draw[black,->,>=latex] (b1) to (b4);
\draw[black,->,>=latex] (b2) to (b3);
\draw[black,->,>=latex] (b2) to (b4);
\draw[black,->,>=latex] (b4) to (b3);
\end{tikzpicture}}}$
 \caption{Kontsevich--Shoikhet cocycle $\Theta_2\in H^1(\odfcGC{2}{0}{1})$.} \label{figKS}
\end{figure}

Of special interest for our purpose is the incarnation of the $\Theta$-cocycle in dimension $3$, dubbed~$\Theta_3$ in the following. The latter is a combination of bi-oriented graphs with $N=6$ and $k=7$\footnote{\label{footMC}More generally, the Maurer--Cartan element $\Upsilon_{\Theta_d}$ associated with the incarnation $\Theta_d$ of the $\Theta$-graph in dimension $d$ is a sum $\Upsilon_{\Theta_d}=\raisebox{-0.4ex}{\hbox{\begin{tikzpicture}[scale=0.5, >=stealth']
\tikzstyle{b}=[circle, draw, fill, minimum size=2, inner sep=0.02]
\node [b] (b2) at (0,0) {1};
\node [b] (b3) at (1.8,0) {2};
\draw[black,->,>=latex][decoration={markings,mark=at position .6 with {\arrow[scale=1.35]{latex}}},
 postaction={decorate},
 shorten >=0.4pt] (b2) to (b3);
\end{tikzpicture}}}+\Theta_d+\cdots=\sum_{p\geq0}\Upsilon^p_{\Theta_d}$ where the graph $\Upsilon^p_{\Theta_d}$ has $N=2p\, (d-1)+2$ vertices and $k=2pd+1$ edges.

For $d=1$, we recover the sum of graphs $\Upsilon_{\Theta_1}:=\sum_{p\geq0}
 \frac{1}{(2p+1)!}
\raisebox{-4.5ex}{\hbox{\begin{tikzpicture}[scale=0.5, >=stealth']
\tikzstyle{b}=[circle, draw, fill, minimum size=2, inner sep=0.02]
\node [b] (b2) at (0,0) {1};
\node [] (b1) at (1,0.25) {$\vdots$};
\node [] (b4) at (1,-1) {\scriptsize{$2p+1$ edges}};
\node [b] (b3) at (2,0) {2};
\draw[black] (b2) to[out=55, in=125, looseness=1.] (b3);
\draw[black] (b2) to[out=-55, in=-125, looseness=1.] (b3);
\end{tikzpicture}}}$ yielding the Moyal star-commutator when represented on the algebra of functions on a symplectic manifold \cite{Khoroshkin2014}.
} (see Appendix \ref{section:Incarnation}) that will be argued to provide an obstruction to the universal quantization of Lie bialgebroids in Section \ref{section:Actions}.

\section{Universal models}
\label{section:Actions}
The aim of the present section is to introduce universal models of multi-oriented graphs (cf.\ Section~\ref{section:Multi)-oriented graph complexes}) for Lie bialgebroids~-- and variations thereof~-- using the graded geometric picture reviewed in Section~\ref{section:Graded geometry}. We start by providing an abstract characterisation of universal models and emphasise their relevance to address questions related to formality theory and deformation quantization.
 We then review universal models for Lie bialgebras (and their ``quasi'' versions) before moving on to the Lie bialgebroid case. We conclude by discussing the implications of our results regarding the deformation quantization problem for Lie bialgebroids.

\subsection{Armchair formality theory}\label{Armchair formality theory}
Let $(\alg,\delta,\brdot_\alg)$ be a dg Lie algebra and denote $(H(\alg),0,\brdot_{H(\alg)})$ the associated cohomology endowed with the canonical dg Lie algebra structure inherited from $\alg$ (with trivial differential).
 Let furthermore $\Phi\colon (H(\alg),0)\iso (\alg,\delta)$ be a quasi-isomorphism of complexes. Generically, $\Phi$~fails to preserve the additional graded Lie structures~-- i.e., $\Phi(\br{x}{y}_{H(\alg)})\neq\br{\Phi(x)}{\Phi(y)}_\alg$ with $x,y\in H(\alg)$~-- so that $\Phi$ is {\it not} a morphism of dg Lie algebras. According to the \textit{homotopy transfer theorem},\footnote{See, e.g., \cite[Theorem~10.3.1]{Loday2012} for a statement as well as Chapter~10 therein for details and history.} any dg Lie algebra $\alg$ is quasi-isomorphic (as a \Lieinf-algebra) to its cohomology $H(\alg)$ endowed with a certain \Lieinf-structure $(H(\alg),l)$ {deforming} the canonical dg Lie algebra structure $(H(\alg),0,\brdot_{H(\alg)})$, i.e., $l_1=0$, $l_2=\brdot_{H(\alg)}$ and the higher order brackets $l_{>2}$ are transferred from the dg Lie algebra structure on $\alg$. In other words, any quasi-isomorphism of complexes $\Phi\colon (H(\alg),0)\iso (\alg,\delta)$ can be upgraded to a quasi-isomorphism of \Lieinf-algebras $\U\colon (H(\alg),l)\iso(\alg,\delta,\brdot_\alg)$ with $\U_1=\Phi$. If the higher brackets $l_{>2}$ vanish, then $(H(\alg),0,\brdot_{H(\alg)})$ and $(\alg,\delta,\brdot_\alg)$ are quasi-isomorphic as \Lieinf-algebras and $\alg$ is said to be \textit{formal}. The homotopy transfer theorem thus allows to reduce questions regarding formality (such as existence of formality maps and their classification) to the study of the space of \Lieinf-structures deforming the canonical dg Lie structure $(H(\alg),0,\brdot_{H(\alg)})$. The relevant deformation theory is controlled by the {Chevalley--Eilenberg} dg Lie algebra $\CE\big(H(\alg)\big)$ endowed with the {Nijenhuis--Richardson bracket} $\brdot_{\NR}$ and the differential $\deltaS:=\big[\brdot_{H(\alg)},\cdot\big]_{\NR}$. Since we are interested in formality maps given by universal formulae, our aim is to introduce~-- for each deformation quantization problem at hand~-- a~\textit{universal model} for the deformation theory of $H(\alg)$ in the guise of a dg Lie algebra of graphs (collectively denoted $\GC$) together with a~morphism of dg Lie algebras $\GC\to \CE\big(H(\alg)\big)$.

\begin{Example}[universal model for polyvector fields]
The paradigmatic example of the above construction is due to M.~Kontsevich \cite{Kontsevich1997} in the context of the deformation quantization problem for Poisson manifolds.

In this context, the quasi-isomorphism of complexes is provided by the HKR map $\Phi_\mathsf{HKR}\colon\Tpoly\allowbreak \iso \Dpoly$ between:
\begin{itemize}\itemsep=0pt
\item $\Tpoly$: the Schouten graded Lie algebra of polyvector fields on the affine space $\mathbb R^m$,
\item $\Dpoly$: the Hochschild dg Lie algebra of multidifferential operators on $\mathbb R^m$.
\end{itemize}
According to the previous reasoning, the existence of a formality map $\U\colon \Tpoly\iso \Dpoly$ can be probed by studying the deformation theory of the Schouten algebra $(\Tpoly,0,\brdot_{\mathsf S})$, controlled by the Chevalley--Eilenberg dg Lie algebra $\CE(\Tpoly)$. In the formulation of his Formality conjecture~\cite{Kontsevich1997}, M.~Kontsevich introduced a~dg Lie algebra of {graphs}~-- denoted~$\fGC_2$~-- together with a morphism of dg Lie algebras $\fGC_2\to \CE(\Tpoly)$ given by explicit local formulae.\footnote{The explicit formulae defining the morphism $\fGC_2\to \CE(\Tpoly)$ take advantage of the graded geometric formulation of Poisson manifolds as dg symplectic manifolds of degree~1, see~\cite{Willwacher2015} for the affine space case (cf.\ also the earlier work \cite{Merkulov2008} as well as \cite{Morand2019a} for a generalisation to dg symplectic manifolds of arbitrary degree), \cite{Jost2012}~for a~globalisation to any smooth manifolds and~\cite{Dolgushev2012b} for a generalisation to the sheaf of polyvector fields on any smooth algebraic variety.} The dg Lie algebra $\fGC_2$ can therefore be interpreted as a {universal model} for $\CE(\Tpoly)$ allowing to reduce important questions related to formality theory to the cohomology of $\fGC_2$:
\begin{itemize}\itemsep=0pt
\item {\textit{Existence}}: Obstructions to the existence of universal formality maps live in $H^1(\fGC_2)$.
\item {\textit{Classification}}: The space of universal formality maps is classified by $H^0(\fGC_2)$.
\end{itemize}
This characterisation of the universal solutions to deformation quantization problems via the cohomology of suitable graph complexes can be generalised to other algebro-geometric structures.
In the next sections, we will mimic the Kontsevich construction for Lie bialgebras and Lie bialgebroids by resorting to models of (multi)-oriented graphs.
\end{Example}

\subsection{Lie bialgebras}
\label{section:Lie bialgebras2}
We start by reviewing some known results regarding universal models on Lie bialgebras, see, e.g., \cite{Merkulov2016,Merkulov2016a, Willwacher2015c}. Let us come back to the graded symplectic manifold $T^*(\alg[1])\simeq(\alg\oplus\alg^*)[1]$ of Section~\ref{section:Lie bialgebras} endowed with a set of homogeneous local coordinates $\big\{\underset{1}{\xi^a},\underset{1}{\zeta_a}\big\}$, with $a\in\pset{1,\dots,\dim\alg}$. Using this set of coordinates, one can (locally) endow the graded algebra of functions~$\foncg$ with a natural structure of $\dGraplain_{3}$-algebra, with $\dGraplain_{3}$ the operad of 1-directed graphs in dimension~3. Explicitly, we define a morphism of operads\footnote{Here $\End_V$ stands for the endomorphism operad associated with the (graded) vector space $V$.} $\Rep^\alg\colon\dGraplain_{3}\rightarrow \End_{\foncg}$ via the following sequence of morphisms of graded\footnote{Recall that the grading of a graph in $\odGra{d}{j}{i}$ is given by $|\gamma|=k(1-d)$ with $k$ the number of edges.} vector spaces, for $N\geq1$:
\begin{align}
\Rep^\alg_N\colon \ \dGraplain_{3}(N)\otimes\foncg^{\otimes N}& \to\foncg,\nonumber\\
\Rep^\alg_N(\graphGam)(f_1\otimes \cdots \otimes f_N)& =\mu_N\bigg( \underset{e\in E_\graphGam}{\prod}\Delta_{e}(f_1\otimes \cdots \otimes f_N)\bigg),
\label{Orientationmorphism}
\end{align}
where
\begin{itemize}\itemsep=0pt
\item The $f_i$'s are functions on $(\alg\oplus\alg^*)[1]$.
\item The symbol $\mu_N$ denotes the multiplication map on $N$ elements:
\begin{align*}
\mu_N\colon\ \foncg^{\otimes N}& \to \foncg,\\
 f_1\otimes \cdots \otimes f_N& \mapsto f_1 \cdots f_N.
\end{align*}
\item The product is performed over the set $E_\graphGam$ of edges of the graph $\graphGam\in\dGraplain_{3}(N)$.
\item For each edge $e\in E_\graphGam$ connecting vertices labeled by integers $i$ and $j$, the operator $\Delta_{e}$ is defined as
\begin{displaymath}
\Delta_{\raisebox{-0ex}{\hbox{\begin{tikzpicture}[scale=0.5, >=stealth']
\tikzstyle{w}=[circle, draw, minimum size=4, inner sep=1]
\tikzstyle{b}=[circle, draw, fill, minimum size=2, inner sep=0.02]
\draw (0,0) node [ext] (b1) {\rm \tiny{$i$}};
\draw (1.3,0) node [exttiny] (b2) {\rm \tiny{$j$}};
\draw(0.6,0.3) node {};
\draw[black,->,>=latex,line width=0.2mm]
 (b1) to (b2);
\end{tikzpicture}}}}=\frac{\p }{\p \xi^a_{(i)}}\frac{\p }{\p \zeta_{a}^{(j)}},
\end{displaymath}
where the sub(super)scripts~$(i)$ or~$(j)$ indicate that the derivative acts on the $i$-th or $j$-th factor in the tensor product. Note that $|\Delta_{e}|=-2$, consistently with the grading of edges in $\dGraplain_{3}$, and that $\Delta_{e}\Delta_{e'}=\Delta_{e'}\Delta_{e}$ as is consistent with the fact that edges are bosonic for~$d$ odd.
\end{itemize}
We refer to \cite[Proposition 5.1]{Morand2019a} for a proof (in a slightly more general context) that the map $\Rep^\alg$ is well-defined and satisfies the axioms of a morphism of operads. The representation $\Rep^\alg$ yields a sequence\footnote{More generally, the action of the 3-Gerstenhaber operad $\mGer_3$ on the algebra of functions factors through $\dGraplain_{3}$ as $\mGer_3\longhookrightarrow\dGraplain_{3}\overset{\Rep^\alg}{\longrightarrow}\End_{\foncg}$, see footnote \ref{footnoteGer} for more details.} of morphisms of operads $\mLie\pset{-2}\overset{\gamma_0}{\longhookrightarrow}\dGraplain_{3}\overset{\Rep^\alg}{\longrightarrow}\End_{\foncg}$ mapping the generator of $\mLie\pset{-2}$ to the graded Poisson bracket \eqref{Poissonbi} via the graph $\raisebox{-0.5ex}{\hbox{\begin{tikzpicture}[scale=0.5, >=stealth']
\tikzstyle{b}=[circle, draw, fill, minimum size=2, inner sep=0.02]
\node [b] (b2) at (0,0) {1};
\node [b] (b3) at (1.8,0) {2};
\draw[black,->,>=latex]
 (b2) to (b3);
\end{tikzpicture}}}:=\raisebox{-0.7ex}{\hbox{\begin{tikzpicture}[scale=0.5, >=stealth']
\tikzstyle{b}=[circle, draw, fill, minimum size=2, inner sep=0.02]
\node [ext] (b2) at (0,0) {1};
\node [ext] (b3) at (1.8,0) {2};
\draw[black,->,>=latex]
 (b2) to (b3);
\end{tikzpicture}}}-\raisebox{-0.7ex}{\hbox{\begin{tikzpicture}[scale=0.5, >=stealth']
\tikzstyle{b}=[circle, draw, fill, minimum size=2, inner sep=0.02]
\node [ext] (b2) at (0,0) {2};
\node [ext] (b3) at (1.8,0) {1};
\draw[black,->,>=latex]
 (b2) to (b3);
\end{tikzpicture}}}$.

Although this action of $\dGraplain_{3}$ on $\foncg$ is well-defined, it generically fails to preserve the various subalgebras introduced in Section~\ref{section:Lie bialgebras}.

\begin{Example}
Consider the following action of a cycle graph
\[
\Rep^\alg_2\big(\raisebox{-1.5ex}{\hbox{\begin{tikzpicture}[scale=0.35, >=stealth']
\tikzstyle{b}=[circle, draw, fill, minimum size=2, inner sep=0.02]
\node [ext] (b2) at (0,0) {1};
\node [ext] (b3) at (2,0) {2};
\draw[black,->,>=latex] (b2) to[out=45, in=140, looseness=1.] (b3);
\draw[black,<-,>=latex] (b2) to[out=-45, in=-125, looseness=1.] (b3);
\end{tikzpicture}}}\big)(f_1\otimes f_2)=\frac{\p^2 f_1}{\p\xi^a\p\zeta_b}\frac{\p^2 f_2}{\p\zeta_a\p\xi^b},
\]
\begin{itemize}
 \item on $f_1=f_2\sim \xi\xi\zeta\in\A^\alg_\text{\rm Lie-quasi}$ yields $\Rep^\alg_2\big(\raisebox{-1.5ex}{\hbox{\begin{tikzpicture}[scale=0.35, >=stealth']
\tikzstyle{b}=[circle, draw, fill, minimum size=2, inner sep=0.02]
\node [ext] (b2) at (0,0) {1};
\node [ext] (b3) at (2,0) {2};
\draw[black,->,>=latex] (b2) to[out=45, in=140, looseness=1.] (b3);
\draw[black,<-,>=latex] (b2) to[out=-45, in=-125, looseness=1.] (b3);
\end{tikzpicture}}}\big)(f_1\otimes f_2)\sim\xi\xi\notin\A^\alg_\text{\rm Lie-quasi}$,
 \item on $f_1=f_2\sim \xi\zeta\zeta\in\A^\alg_\text{\rm quasi-Lie}$ yields $\Rep^\alg_2\big(\raisebox{-1.5ex}{\hbox{\begin{tikzpicture}[scale=0.35, >=stealth']
\tikzstyle{b}=[circle, draw, fill, minimum size=2, inner sep=0.02]
\node [ext] (b2) at (0,0) {1};
\node [ext] (b3) at (2,0) {2};
\draw[black,->,>=latex] (b2) to[out=45, in=140, looseness=1.] (b3);
\draw[black,<-,>=latex] (b2) to[out=-45, in=-125, looseness=1.] (b3);
\end{tikzpicture}}}\big)(f_1\otimes f_2)\sim\zeta\zeta\notin\A^\alg_\text{\rm quasi-Lie}$,
\item {\it A fortiori}, the action of $\raisebox{-1.5ex}{\hbox{\begin{tikzpicture}[scale=0.35, >=stealth']
\tikzstyle{b}=[circle, draw, fill, minimum size=2, inner sep=0.02]
\node [ext] (b2) at (0,0) {1};
\node [ext] (b3) at (2,0) {2};
\draw[black,->,>=latex] (b2) to[out=45, in=140, looseness=1.] (b3);
\draw[black,<-,>=latex] (b2) to[out=-45, in=-125, looseness=1.] (b3);
\end{tikzpicture}}}$ fails to preserve $\A^\alg_{\rm Lie}:= \A^\alg_\text{\rm Lie-quasi}\cap \A^\alg_\text{\rm quasi-Lie}$.
\end{itemize}
\end{Example}

This defect can be cured by carefully restricting the space of graphs, as embodied in the following lemma.
\begin{Lemma}\label{PropalgLie}Let $\alg$ be a vector space.
\begin{itemize}\itemsep=0pt
\item The graded algebra of functions $\foncg$ is endowed with a structure of $\dGraplain_{3}$-algebra via the action of $1$-directed graphs.
\item The graded subalgebra $\A^\alg_\text{\rm Lie-quasi}$ is endowed with a structure of $\sodGra{3}{0}{1}$-algebra via the action of sourced $1$-directed graphs.
\item The graded subalgebra $\A^\alg_\text{\rm quasi-Lie}$ is endowed with a structure of $\sodGra{3}{0}{-1}$-algebra via the action of sinked $1$-directed graphs.
\item The graded subalgebra $\A^\alg_\text{\rm Lie}$ is endowed with a structure of $\sodGra{3}{0}{1,-1}$-algebra via the action of $1$-directed graphs being both sourced and sinked.
\end{itemize}
\end{Lemma}
\begin{proof}
Let us first consider the action of a graph containing a source vertex and let $f$ denote the function decorating the source upon the action of $\Rep^\alg$. Then the differential operator acting on~$f$ is of the form $\frac{\p}{\p\xi}\cdots\frac{\p}{\p\xi}f$. Assuming that $f$ belongs to the subalgebra $\A^\alg_\text{\rm Lie-quasi}$ ensures that~$f$ is at least linear in~$\zeta$ (by definition of $\A^\alg_\text{\rm Lie-quasi}$, cf.\ Section~\ref{section:Lie bialgebras}). Hence the function obtained as the outcome of the action of a sourced graph on $\A^\alg_\text{\rm Lie-quasi}$ either vanishes or is at least linear in~$\zeta$, so that it belongs to~$\A^\alg_\text{\rm Lie-quasi}$. The subalgebra $\A^\alg_\text{\rm Lie-quasi}$ is therefore closed under the action of sourced graphs. A similar reasoning shows that $\A^\alg_\text{\rm quasi-Lie}$ is closed under the action of sinked graphs. It follows that the intersection $\A^\alg_{\rm Lie}:= \A^\alg_\text{\rm Lie-quasi}\cap \A^\alg_\text{\rm quasi-Lie}$ is closed under the action of graphs possessing at least one source and one sink.
\end{proof}

As noted previously (see Section~\ref{section:Operads}) oriented graphs necessarily possess at least one source and one sink. Denoting collectively by $\A^\alg_\text{\rm sub}$ the three subalgebras $\A^\alg_\text{\rm Lie-quasi}$, $\A^\alg_\text{\rm quasi-Lie}$ and $\A^\alg_\text{\rm Lie}\subset\foncg$, the previous lemma thus ensures that $\A^\alg_\text{\rm sub}$ is endowed with a structure of $\odGra{3}{0}{1}$-algebra via the action of oriented 1-directed graphs, i.e., the representation $\Rep^\alg$ induces morphisms of operads $\Rep^\alg\colon \odGra{3}{0}{1}\rightarrow \End_{\A^\alg_\text{\rm sub}}$.
Since we are primarily interested in the cohomology of the associated graph complexes, it is sufficient for our purpose to resort to the smaller operad of oriented graphs as this is where the cohomology lies (cf.\ the third item of Theorem~\ref{thmWM}).
Applying $\Def(\mLie\pset{-2}\to\cdot)$ on both sides of the morphism $\Rep^\alg$ yields the following proposition:
\begin{Proposition}\label{propmorphism}\quad
\begin{itemize}\itemsep=0pt
\item The morphism of operads $\Rep^\alg\colon \dGraplain_{3}\rightarrow \End_{\foncg}$ induces a morphism of dg Lie algebras
\begin{gather*}
\big(\dfGC{3}{1},\deltabr,\brdot\big)\to \big(\CE\big(\foncg[2]\big),\deltaS,\brdot_{\NR}\big).
\end{gather*}
\item The morphisms of operads $\Rep^\alg:\odGra{3}{0}{1}\rightarrow \End_{\A^\alg_\text{\rm sub}}$ induce morphisms of dg Lie algebras
\begin{gather}
\big(\odfGC{3}{0}{1},\deltabr,\brdot\big)\to \big(\CE(\A^\alg_\text{\rm sub}[2]),\deltaS,\brdot_{\NR}\big),\label{univmordg Lie algebra2}
\end{gather}
where $\A^\alg_\text{\rm sub}$ collectively denotes the subalgebras $\A^\alg_\text{\rm Lie-quasi}, \A^\alg_\text{\rm quasi-Lie},\A^\alg_\text{\rm Lie}\!\subset\!\foncg$.
\end{itemize}
\end{Proposition}

\begin{proof}
The proposition follows straightforwardly from Lemma \ref{PropalgLie} together with the equivariance of $\Rep^\alg$.
\end{proof}

We denoted $\CE(\alg)$ the Chevalley--Eilenberg cochain space (in the adjoint representation) associated with the vector space $\alg$ while $\brdot_{\NR}$ stands for the {Nijenhuis--Richardson bracket} and $\deltaS:=\br{\pbdot^\alg_\symp}{\cdot}_{\NR}$ for the Chevalley--Eilenberg differential associated with the graded Poisson bracket \eqref{Poissonbi}. In plain words, \Prop{propmorphism} states that the graph complex $\dfcGC{3}{1}$ provides a universal model for the deformation theory of proto-Lie bialgebras while $\odfcGC{3}{0}{1}$ can be seen as a universal model for the deformation theory of Lie-quasi, quasi-Lie and Lie bialgebras. This last fact combined with the cohomology computations reviewed in Section~\ref{section:Cohomology} yields the following well-known result:
\begin{Corollary}\label{corGroalg}
The Grothendieck--Teichm\"uller group acts via \Lieinf-automorphisms on the deformation complexes of Lie-quasi, quasi-Lie and Lie bialgebras.
\end{Corollary}
\begin{proof}
Going to the zeroth cohomology in \eqref{univmordg Lie algebra2} yields a map
\[
H^0\big(\odfGC{3}{0}{1}|\deltabr\big)\to H^0\big(\CE(\A^\alg_\text{\rm sub}[2])|\deltaS\big).
\]
Hence any non-trivial zero degree cocycle $\gamma$ in $\odfGC{3}{0}{1}$ is mapped to a \Lieinf-derivation $\Rep^\alg(\gamma)$ of the dg Lie algebras $\big(\A^\alg_\text{\rm sub},\pbdot^\alg_\symp\big)$, thus ensuring that $\exp\big(\Rep^\alg(\gamma)\big)$ is a \Lieinf-automorphism thereof. Therefore the pro-unipotent group $\exp\big(H^0\big(\odfGC{3}{0}{1}\big)\big)$ acts on the deformation complexes $\big(\A^\alg_\text{\rm sub},\pbdot^\alg_\symp\big)$ via \Lieinf-automorphisms. Focusing on connected graphs, the second item of Theorem~\ref{thmWM} ensures that $H^0(\odfcGC{3}{0}{1})\simeq H^0(\fcGC_2)\simeq \grt_1$, where the last equivalence is T.~Willwacher's theorem~\cite{Willwacher2015}. Exponentienting thus yields an action of the pro-unipotent Grothendieck--Teichm\"uller group $\GRT_1\equiv\exp(\grt_1)$ on the deformation complexes $\big(\A^\alg_\text{\rm sub},\pbdot^\alg_\symp\big)$ of Lie-quasi, quasi-Lie and Lie bialgebras via \Lieinf-automorphisms.
\end{proof}

This result is consistent with the known action of $\GRT_1$ (via the $\GRT_1$-torsor of Drinfeld associators) on quantization maps for Lie bialgebras \cite{Etingof1995,Ginot2016,Hinich2014,Merkulov2016,Merkulov2016a, Tamarkin2007b} and Lie-quasi bialgebras~\cite{Enriquez2008,Sakalos2013}.
In comparison, there is no such action of $\GRT_1$ on the deformation complex of proto-Lie bialgebras as the above \Lieinf-automorphisms become trivial when acting on $\foncg$ (consistently with the fact that $H^0(\dfcGC{3}{1})\simeq H^0(\fcGC_3) \simeq \corps$, cf.\ Table~\ref{Tabcohom}). In that sense, there seems to be no true ``intermediate case'' between the proto-Lie bialgebra case and the Lie bialgebra case since restricting to oriented graphs allows to preserve all three subalgebras at once. As we will see, this fact is in contradistinction with the ``bialgebroid'' case in which Lie-quasi and quasi-Lie bialgebroids provide a true intermediate case between proto-Lie and Lie bialgebroids.
\subsection{Lie bialgebroids}
\label{section:Lie bialgebroids2}
We now move on to the main result of this note by introducing novel universal models for the deformation complexes of the family of variations on Lie bialgebroids reviewed in Section~\ref{section:Lie bialgebroids} and Appendix~\ref{section:Appendix}. Let $E\overset{\pi}{\to} \M$ be a vector bundle and consider the graded symplectic manifold $\big(\TE,\symp\big)$, with $\TE\equiv T^*[2]E[1]$, coordinatised by the set of homogeneous coordinates $\big\{\underset{0}{x^\mu},\underset{1}{\xi^a},\underset{1}{\zeta_a},\underset{2}{p_\mu}\big\}$. The algebra of functions on $\TE$ carries a natural action $\Rep^E\colon \dGra{3}{2}\rightarrow \End_{\foncTE}$ of the operad $\dGra{3}{2}$ of bi-directed graphs containing both black and red directions. The action $\Rep^E$ is defined similarly as the action~\eqref{Orientationmorphism} where for each edge $e\in E_\graphGam$ connecting vertices labeled by integers~$i$ and~$j$, the operator $\Delta_{e}$ is defined as
\begin{displaymath}
\Delta_{\raisebox{-0ex}{\hbox{\begin{tikzpicture}[scale=0.5, >=stealth']
\tikzstyle{w}=[circle, draw, minimum size=4, inner sep=1]
\tikzstyle{b}=[circle, draw, fill, minimum size=2, inner sep=0.02]
\draw (0,0) node [ext] (b1) {\rm \tiny{$i$}};
\draw (1.3,0) node [exttiny] (b2) {\rm \tiny{$j$}};
\draw(0.6,0.3) node {};
\draw[black,->,>=latex,line width=0.2mm][decoration={markings,mark=at position .55 with {\arrow[scale=1,linkcolor]{latex}}},
 postaction={decorate},
 shorten >=0.4pt] (b1) to (b2);
\end{tikzpicture}}}}=\frac{\p }{\p x^\mu_{(i)}}\, \frac{\p }{\p p_\mu^{(j)}}
 ,\qquad
\Delta_{\raisebox{-0ex}{\hbox{\begin{tikzpicture}[scale=0.5, >=stealth']
\tikzstyle{w}=[circle, draw, minimum size=4, inner sep=1]
\tikzstyle{b}=[circle, draw, fill, minimum size=2, inner sep=0.02]
\draw (0,0) node [ext] (b1) {\rm \tiny{$i$}};
\draw (1.3,0) node [exttiny] (b2) {\rm \tiny{$j$}};
\draw(0.6,0.3) node {};
\draw[black,->,>=latex,line width=0.2mm][decoration={markings,mark=at position 0 with {\arrowreversed[scale=1,linkcolor]{latex}}},
 postaction={decorate},
 shorten >=0.4pt] (b1) to (b2);
\end{tikzpicture}}}}=\frac{\p }{\p \xi^a_{(i)}}\frac{\p }{\p \zeta_{a}^{(j)}},
\end{displaymath}
where the sub(super)scripts $(i)$ or $(j)$ indicate that the derivative acts on the $i$-th or $j$-th factor in the tensor product.
As in the Lie bialgebra case, we note that the operator $\Delta_e$ has grading $|\Delta_{e}|=-2$, consistently with the grading of edges in~$\dGraplain_{3}$, and that $\Delta_{e}\Delta_{e'}=\Delta_{e'}\Delta_{e}$ as is consistent with the fact that edges are bosonic for~$d$ odd.

The representation $\Rep^E$ maps the graph
\[
\raisebox{-0.5ex}{\hbox{\begin{tikzpicture}[scale=0.5, >=stealth']
\tikzstyle{b}=[circle, draw, fill, minimum size=2, inner sep=0.02]
\node [b] (b2) at (0,0) {2};
\node [b] (b3) at (1.8,0) {3};
\draw[black,->,>=latex][decoration={markings,mark=at position .6 with {\arrow[scale=1.35,linkcolor]{latex}}},
 postaction={decorate},
 shorten >=0.4pt] (b2) to (b3);
\end{tikzpicture}}}=\raisebox{-0.7ex}{\hbox{\begin{tikzpicture}[scale=0.5, >=stealth']
\tikzstyle{b}=[circle, draw, fill, minimum size=2, inner sep=0.02]
\node [ext] (b2) at (0,0) {1};
\node [ext] (b3) at (2,0) {2};
\draw[black,->,>=latex][decoration={markings,mark=at position .6 with {\arrow[scale=1.35,linkcolor]{latex}}},
 postaction={decorate},
 shorten >=0.4pt] (b2) to (b3);
\end{tikzpicture}}}
+
\raisebox{-0.7ex}{\hbox{\begin{tikzpicture}[scale=0.5, >=stealth']
\tikzstyle{b}=[circle, draw, fill, minimum size=2, inner sep=0.02]
\node [ext] (b2) at (0,0) {1};
\node [ext] (b3) at (2,0) {2};
\draw[black,->,>=latex][decoration={markings,mark=at position .13 with {\arrowreversed[scale=1.35,linkcolor]{latex}}},
 postaction={decorate},
 shorten >=0.4pt] (b2) to (b3);
\end{tikzpicture}}}
\, -\, \big(
\raisebox{-0.7ex}{\hbox{\begin{tikzpicture}[scale=0.5, >=stealth']
\tikzstyle{b}=[circle, draw, fill, minimum size=2, inner sep=0.02]
\node [ext] (b2) at (0,0) {2};
\node [ext] (b3) at (2,0) {1};
\draw[black,->,>=latex][decoration={markings,mark=at position .6 with {\arrow[scale=1.35,linkcolor]{latex}}},
 postaction={decorate},
 shorten >=0.4pt] (b2) to (b3);
\end{tikzpicture}}}
+
\raisebox{-0.7ex}{\hbox{\begin{tikzpicture}[scale=0.5, >=stealth']
\tikzstyle{b}=[circle, draw, fill, minimum size=2, inner sep=0.02]
\node [ext] (b2) at (0,0) {2};
\node [ext] (b3) at (2,0) {1};
\draw[black,->,>=latex][decoration={markings,mark=at position .13 with {\arrowreversed[scale=1.35,linkcolor]{latex}}},
 postaction={decorate},
 shorten >=0.4pt] (b2) to (b3);
\end{tikzpicture}}}
\big)
\]
 towards the graded Poisson bracket \eqref{eqPois} and furthermore yields a sequence\footnote{\label{footnoteGer}Similarly as in the Lie bialgebra case, the action of the 3-Gerstenhaber operad $\mGer_3$ (also called $\mathsf{e}_3$) on the algebra of functions $\foncTE$ factors through $\dGra{3}{2}$ as $\mGer_3\overset{i_3}{\longhookrightarrow}\dGra{3}{2}\overset{\Rep^E}{\longrightarrow}\End_{\foncTE}$ where the embedding is explicitly given by the following action on generators $a_1\w a_2,\pb{a_1}{a_2}\in\mGer_3(2)$:
\begin{itemize}\itemsep=0pt
\item $i_3(a_1\w a_2)=\raisebox{-0.8ex}{\hbox{
\begin{tikzpicture}[scale=0.5, >=stealth']
\tikzstyle{w}=[circle, draw, minimum size=4, inner sep=1]
\tikzstyle{b}=[circle, draw, fill, minimum size=4, inner sep=1]
\node [ext] (b4) at (2,0) {\rm 1};
\node [ext] (b5) at (3.8,0) {\rm 2};
\end{tikzpicture}}}$ with $\w$ the graded commutative associative product of degree $0$
\item $i_3(\pb{a_1}{a_2})=\raisebox{-0.5ex}{\hbox{\begin{tikzpicture}[scale=0.5, >=stealth']
\tikzstyle{b}=[circle, draw, fill, minimum size=2, inner sep=0.02]
\node [b] (b2) at (0,0) {2};
\node [b] (b3) at (1.8,0) {3};
\draw[black,->,>=latex][decoration={markings,mark=at position .6 with {\arrow[scale=1.35,linkcolor]{latex}}},
 postaction={decorate},
 shorten >=0.4pt] (b2) to (b3);
\end{tikzpicture}}}$ with $\pbdot$ the graded Lie bracket of degree $-2$.
\end{itemize}
 } of morphisms of operads
 \begin{displaymath}
 \mLie\pset{-2}\overset{\gamma_0}{\longhookrightarrow}\dGra{3}{2}\overset{\Rep^E}{\longrightarrow}\End_{\foncTE}.
 \end{displaymath}
Although the action of $\dGra{3}{2}$ on $\foncTE$ is well-defined and satisfies the axioms of a morphism of operads (cf.\ the discussion around~\eqref{Orientationmorphism}), it generically fails to preserve the various subalgebras introduced in Section~\ref{section:Lie bialgebroids}.
\begin{Example}
Considering the following action of a cycle graph
\begin{displaymath}
\displaystyle\Rep^E_2\big(
\raisebox{-1.7ex}{\hbox{\begin{tikzpicture}[scale=0.4, >=stealth']
\tikzstyle{b}=[circle, draw, fill, minimum size=2, inner sep=0.02]
\node [ext] (b2) at (0,0) {1};
\node [ext] (b3) at (2,0) {2};
\draw[black,->,>=latex][decoration={markings,mark=at position .55 with {\arrow[scale=1,linkcolor,rotate=4,yshift=0.017cm,xshift=0.08cm]{latex}}}, postaction={decorate}, shorten >=0.4pt] (b2) to[out=65, in=125, looseness=1.] (b3);
\draw[black,->,>=latex][decoration={markings,mark=at position .35 with {\arrow[scale=1,linkcolor,rotate=-3.5,yshift=0.cm,xshift=0.13cm]{latex}}}, postaction={decorate},shorten >=0.4pt] (b3) to[out=-125, in=-45, looseness=1.] (b2);
\end{tikzpicture}}}
\big)(f_1\otimes f_2)=\frac{\p^2 f_1}{\p x^\mu\p p_\nu}\frac{\p^2 f_2}{\p x^\nu\p p_\mu},
\end{displaymath}
 \begin{itemize}\itemsep=0pt
 \item on $f_1=f_1(x){}^\mu{}_a p_\mu \xi^a$ and $f_2=f_2(x){}^\mu{}_a p_\mu \xi^a\in\A^E_\text{\rm Lie-quasi}$ yielding
 \begin{displaymath}
\displaystyle\Rep^E_2\big(
\raisebox{-1.7ex}{\hbox{\begin{tikzpicture}[scale=0.4, >=stealth']
\tikzstyle{b}=[circle, draw, fill, minimum size=2, inner sep=0.02]
\node [ext] (b2) at (0,0) {1};
\node [ext] (b3) at (2,0) {2};
\draw[black,->,>=latex][decoration={markings,mark=at position .55 with {\arrow[scale=1,linkcolor,rotate=4,yshift=0.017cm,xshift=0.08cm]{latex}}}, postaction={decorate}, shorten >=0.4pt] (b2) to[out=65, in=125, looseness=1.] (b3);
\draw[black,->,>=latex][decoration={markings,mark=at position .35 with {\arrow[scale=1,linkcolor,rotate=-3.5,yshift=0.cm,xshift=0.13cm]{latex}}}, postaction={decorate},shorten >=0.4pt] (b3) to[out=-125, in=-45, looseness=1.] (b2);
\end{tikzpicture}}}
\big)(f_1\otimes f_2)=\p_\mu f_1(x){}^\nu{}_{[a}\p_\nu f_2(x){}^\mu{}_{b]} \xi^a\xi^b \notin\A^E_\text{\rm Lie-quasi},
\end{displaymath}
 \item on $f_1=f_1(x){}^{\mu|a} p_\mu \zeta_a$ and $f_2=f_2(x){}^{\mu|a} p_\mu \zeta_a\in\A^E_\text{\rm quasi-Lie}$ yielding
 \begin{displaymath}
\displaystyle\Rep^E_2\big(
\raisebox{-1.7ex}{\hbox{\begin{tikzpicture}[scale=0.4, >=stealth']
\tikzstyle{b}=[circle, draw, fill, minimum size=2, inner sep=0.02]
\node [ext] (b2) at (0,0) {1};
\node [ext] (b3) at (2,0) {2};
\draw[black,->,>=latex][decoration={markings,mark=at position .55 with {\arrow[scale=1,linkcolor,rotate=4,yshift=0.017cm,xshift=0.08cm]{latex}}}, postaction={decorate}, shorten >=0.4pt] (b2) to[out=65, in=125, looseness=1.] (b3);
\draw[black,->,>=latex][decoration={markings,mark=at position .35 with {\arrow[scale=1,linkcolor,rotate=-3.5,yshift=0.cm,xshift=0.13cm]{latex}}}, postaction={decorate},shorten >=0.4pt] (b3) to[out=-125, in=-45, looseness=1.] (b2);
\end{tikzpicture}}}
\big)(f_1\otimes f_2)=\p_\mu f_1(x){}^{\nu|[a}\p_\nu f_2(x){}^{\mu|b]}\, \zeta_a\zeta_b \notin\A^E_\text{\rm quasi-Lie}.
\end{displaymath}
\noindent A fortiori, the previous graph fails to preserve the intersection $\A^E_{\rm Lie}:= \A^E_\text{\rm Lie-quasi}\cap \A^E_\text{\rm quasi-Lie}$.
\end{itemize}
\end{Example}

Hence $\A^E_\text{\rm Lie-quasi}$, $\A^E_\text{\rm quasi-Lie}$ and $\A^E_\text{\rm Lie}$ generically fail to be $\dGra{3}{2}$-algebras. Similarly to the Lie bialgebra case, this apparent defect can be cured via a suitable restriction of the class of graphs, as performed in the following lemma:
\begin{Lemma}\label{lembialgebroid}
Let $E\overset{\pi}{\to} \M$ be a vector bundle.
\begin{itemize}\itemsep=0pt
\item The graded algebra of functions $\foncTE$ is endowed with a structure of $\dGra{3}{2}$-algebra via the action of bi-directed graphs.
\item The graded subalgebra $\A^E_\text{\rm Lie-quasi}$ is endowed with a structure of $\sodGrablack{3}{1}{1}$-algebra via the action of bi-directed graphs with a black source.
\item The graded subalgebra $\A^E_\text{\rm quasi-Lie}$ is endowed with a structure of $\sodGrared{3}{1}{1}$-algebra via the action of bi-directed graphs with a red source.
\item The graded subalgebra $\A^E_\text{\rm Lie}$ is endowed with a structure of $\sodGra{3}{0}{2}$-algebra via the action of bi-directed graphs with a black source and a red source.
\end{itemize}
\end{Lemma}
\begin{proof}
We mimic the proof in the Lie bialgebra case by considering first the action of a graph containing a black source vertex. Letting $f$ denote the function decorating the black source upon the action of $\Rep^E$, the differential operator acting on $f$ is of the form $\frac{\p}{\p x}\cdots\frac{\p}{\p x}\frac{\p}{\p\xi}\cdots\frac{\p}{\p\xi}f$. Assuming that $f$ belongs to the subalgebra $\A^E_\text{\rm Lie-quasi}$ ensures that $f$ is at least linear either in $p$ or $\zeta$ (by definition of $\A^E_\text{\rm Lie-quasi}$, cf.\ Section~\ref{section:Lie bialgebroids}). Hence the function obtained as the outcome of the action of a black sourced graph on $\A^E_\text{\rm Lie-quasi}$ cannot be of the form $\xi\cdots\xi$ so that it belongs to~$\A^E_\text{\rm Lie-quasi}$. The subalgebra $\A^E_\text{\rm Lie-quasi}$ is therefore closed under the action of black sourced graphs. A similar reasoning shows that $\A^E_\text{\rm quasi-Lie}$ is closed under the action of red sourced graphs.\footnote{Note that graphs with (black or red) sinks preserve neither $\A^E_\text{\rm Lie-quasi}$ nor $\A^E_\text{\rm quasi-Lie}$. } It follows that the intersection $\A^E_{\rm Lie}:= \A^E_\text{\rm Lie-quasi}\cap \A^E_\text{\rm quasi-Lie}$ is closed under the action of graphs possessing both a black and red source.
\end{proof}

As in the Lie bialgebra case, one can make use of the inclusion $\odGra{3}{\bullet}{\bullet}\subset\sodGra{3}{\bullet}{\bullet}$ (see Section \ref{section:Operads}) to extract from the previous lemma an action of (suitable) operads of multi-oriented graphs on the three subalgebras at hand. Explicitly, the subalgebra $\A^E_\text{\rm Lie-quasi}$ (resp. $\A^E_\text{\rm quasi-Lie}$) is acted upon by the operad of bi-directed graphs with oriented black (resp. red) arrows, denoted $\odGrablack{3}{1}{1}$ (resp. $\odGrared{3}{1}{1}$) in the following.\footnote{Alternatively, both actions could be written in terms of $\odGrablack{3}{1}{1}$ (say) by making use of $\Delta_{\raisebox{-0ex}{\hbox{\begin{tikzpicture}[scale=0.5, >=stealth']
\tikzstyle{w}=[circle, draw, minimum size=4, inner sep=1]
\tikzstyle{b}=[circle, draw, fill, minimum size=2, inner sep=0.02]
\draw (0,0) node [ext] (b1) {\rm \tiny{$i$}};
\draw (1.3,0) node [exttiny] (b2) {\rm \tiny{$j$}};
\draw(0.6,0.3) node {};
\draw[black,->,>=latex,line width=0.2mm][decoration={markings,mark=at position 0 with {\arrowreversed[scale=1,linkcolor]{latex}}},
 postaction={decorate},
 shorten >=0.4pt] (b1) to (b2);
\end{tikzpicture}}}}=\frac{\p }{\p \xi^a_{(i)}}\frac{\p }{\p \zeta_{a}^{(j)}}$ for $\A^E_\text{\rm Lie-quasi}$ and $\Delta'_{\raisebox{-0ex}{\hbox{\begin{tikzpicture}[scale=0.5, >=stealth']
\tikzstyle{w}=[circle, draw, minimum size=4, inner sep=1]
\tikzstyle{b}=[circle, draw, fill, minimum size=2, inner sep=0.02]
\draw (0,0) node [ext] (b1) {\rm \tiny{$i$}};
\draw (1.3,0) node [exttiny] (b2) {\rm \tiny{$j$}};
\draw(0.6,0.3) node {};
\draw[black,->,>=latex,line width=0.2mm][decoration={markings,mark=at position 0 with {\arrowreversed[scale=1,linkcolor]{latex}}},
 postaction={decorate},
 shorten >=0.4pt] (b1) to (b2);
\end{tikzpicture}}}}=\frac{\p }{\p \zeta_{a}^{(i)}}\frac{\p }{\p \xi^a_{(j)}}$ for $\A^E_\text{\rm quasi-Lie}$.}

Crucially, preserving the intersection $\A^E_{\rm Lie}:= \A^E_\text{\rm Lie-quasi}\cap \A^E_\text{\rm quasi-Lie}$ requires to orient {\it both} directions, thus yielding an action of $\odGra{3}{0}{2}$ on $\A^E_{\rm Lie}$. As pointed out in Section~\ref{section:Cohomology} (see Theorem~\ref{thmWM}), the number of oriented colors (in contradistinction with the number of directed colors) is the relevant factor to compute the respective cohomology. From this simple observation will then follow that the Lie bialgebroid case differs essentially from the dual cases of Lie-quasi and quasi-Lie bialgebroids.
Applying $\Def(\mLie\pset{-2}\to\cdot)$ on both sides of $\Rep^E\colon \odGra{3}{\bullet}{\bullet}\rightarrow \End_{\A^E}$ yields the following proposition:
\begin{Proposition}\label{propmorphismoid}\quad
\begin{itemize}\itemsep=0pt
\item The morphism of operads $\Rep^E\colon \dGra{3}{2}\rightarrow \End_{\foncTE}$ induces a morphism of dg Lie algebras
\begin{gather}
\big(\dfGC{3}{2},\deltabr,\brdot\big)\to \big(\CE\big(\foncTE[2]\big),\deltaS,\brdot_{\NR}\big).\label{univmordg Lie algebroid1}
\end{gather}
\item The morphisms of operads $\Rep^E\colon \odGra{3}{1}{1}\rightarrow \End_{\A^E_\text{\rm quasi}}$ induce morphisms of dg Lie algebras
\begin{gather}
\big(\odfGC{3}{1}{1},\deltabr,\brdot\big)\to \big(\CE\big(\A^E_\text{\rm quasi}[2]\big),\deltaS,\brdot_{\NR}\big),\label{univmordg Lie algebroid2}
\end{gather}
where $\A^E_\text{\rm quasi}$ stands for $\A^E_\text{\rm Lie-quasi}$, $\A^E_\text{\rm quasi-Lie}\subset\foncTE$.
\item The morphism of operads $\Rep^E\colon \odGra{3}{0}{2}\rightarrow \End_{\A^E_\text{\rm Lie}}$ induces a morphism of dg Lie algebras
\begin{gather*}
\big(\odfGC{3}{0}{2},\deltabr,\brdot\big)\to \big(\CE\big(\A^E_\text{\rm Lie}[2]\big),\deltaS,\brdot_{\NR}\big).
\end{gather*}
\end{itemize}
\end{Proposition}

\begin{proof}The proposition is a direct consequence of Lemma~\ref{lembialgebroid} and the equivariance of the morphism~$\Rep^E$.
\end{proof}

The conventions used are as in the bialgebra case and $\deltaS:=\big[\pbdot^E_\symp,\cdot\big]_{\NR}$ stands for the Chevalley--Eilenberg differential associated with the graded Poisson bracket~\eqref{eqPois}. Going into cohomology and using Theorem~\ref{thmWM} allows to compute the relevant cohomology groups, as summed up in Table~\ref{figgraphsym}.
\begin{table}[h]\centering\small
 \caption{Cohomology groups for Lie bialgebroid structures.} \label{figgraphsym}

 \vspace{1mm}

\begin{tabular}[5pt]{cccc}
 \toprule
 \rule[-0.2cm]{0cm}{0.6cm} Structure&Black&Red&Cohomology\\
 \bottomrule
 \rule[-0.2cm]{0cm}{0.6cm} proto-Lie bialgebroids&directed&directed&$H^\bullet({\dfGC{3}{2}})\simeq H^\bullet({\fGC_3})$\\
 \rule[-0.2cm]{0cm}{0.6cm} Lie-quasi bialgebroids&oriented&directed&$H^\bullet(\odfGCblack{3}{0}{1})\simeq H^\bullet({\fGC_2})$\\
 \rule[-0.2cm]{0cm}{0.6cm} quasi-Lie bialgebroids&directed&oriented&$H^\bullet(\odfGCred{3}{0}{1})\simeq H^\bullet({\fGC_2})$\\
 \rule[-0.2cm]{0cm}{0.6cm} Lie bialgebroids&oriented&oriented&$H^\bullet({\odfGC{3}{0}{2}})\simeq H^\bullet({\fGC_1})$\\

 \bottomrule
 \end{tabular}
\end{table}

The morphism \eqref{univmordg Lie algebroid1} on the deformation complex of proto-Lie bialgebroids can be seen as a~particular subcase of the action of $\fGC_3$ on the deformation complex of Courant algebroids \cite{Morand2019a} when restricted to the split case.\footnote{Here by split Courant algebroids we mean Courant algebroids whose underlying vector bundle is a Whitney sum $E\oplus E^*$.} The latter does not yield interesting structures in degrees~0 and~1 as the dominant level of the relevant cohomology $H^{\bullet}(\fcGC_3)$ is located in degree $-3$. The corresponding cohomology space $H^{-3}(\fcGC_3)$ is a unital commutative algebra spanned by trivalent graphs modulo the {\it IHX relation} where the r\^ole of the unit is played by the $\Theta$-graph (see, e.g.,~\cite{Bar-Natan1995-2}). Given a proto-Lie bialgebroid $E\overset{\pi}{\to} \M$ represented by the Hamiltonian function $\cH$ \big(see~\eqref{Hamiltonian}\big), each trivalent graph $\gamma\in H^{-3}(\fcGC_3)$ yields a cocycle function $\Om_\gamma\in\fonc{\M}$ \big(i.e., such that $\pb{\cH}{\Om_\gamma}_\symp^E=0$\big) thus yielding a conformal flow on the space of proto-Lie bialgebroids on $E$ (cf.~\cite{Morand2019a} for details).

Turning to Lie-quasi and quasi-Lie bialgebroids, the morphism \eqref{univmordg Lie algebroid2} yields a natural extension of the action of the Grothendieck--Teichm\"uller group $\exp\big(H^0(\odfcGC{3}{1}{1})\big)\simeq\GRT_1$ on the deformation complexes of Lie-quasi and quasi-Lie bialgebras to the ``bialgebroid'' case. This is the essence of Theorem~\ref{thmGro}.

\begin{proof}[Proof of Theorem \ref{thmGro}]
The proof is identical to the one of Corollary~\ref{corGroalg} upon substituting $\odfcGC{3}{0}{1}$ with $\odfcGC{3}{1}{1}$, both admitting the same cohomology thanks to the first item of Theorem~\ref{thmWM}.
\end{proof}

We refer to Section \ref{section:Application to quantization} for a discussion of the consequences of this action in the context of deformation quantization. Explicitly, the action of $\GRT_1$ is through graphs with one oriented color (either black or red) and as such generically {\it fails} to preserve the sub-deformation complex~$\A^E_{\rm Lie}$ for Lie bialgebroids. This is in contradistinction with the Lie bialgebra case whose deformation complex $\A^\alg_{\rm Lie}$ does carry a representation of~$\GRT_1$. Rather the action of $\odfcGC{3}{0}{2}$ endows~$\A^E_{\rm Lie}$ with a new \Lieinf-structure deforming non-trivially the big bracket~\eqref{eqPois}, as captured by Theorem~\ref{TheoremMain}, for which we are now in position to articulate the following proof:

\begin{proof}[Proof of Theorem~\ref{TheoremMain}]
It follows from the third item of Proposition~\ref{propmorphismoid} that Maurer--Cartan elements for the dg Lie algebra $\big(\odfGC{3}{0}{2},\deltabr,\brdot\big)$ are mapped via $\Rep^E$ to universal deformations of the graded Lie algebra $\big(\A^E_\text{\rm Lie},\pbdot^E_\symp\big)$ as a \Lieinf-algebra. As recalled in Section~\ref{section:Cohomology}, applying the second item of Theorem \ref{thmWM} results in a sequence of isomorphisms
$H^1(\odfcGC{3}{0}{2})\!\simeq\! H^1(\odfcGC{3}{0}{1})\!\simeq\! H^1(\fcGC_{1})\!\simeq\!\corps$,
 thus providing a non-trivial cocycle graph $\Theta_3\!\in\! H^1(\odfcGC{3}{0}{2})$ of degree $1$. The latter can be recursively extended to a non-trivial Maurer--Cartan element
\[
\Upsilon_{\Theta_3}:=\raisebox{-0.4ex}{\hbox{\begin{tikzpicture}[scale=0.5, >=stealth']
\tikzstyle{b}=[circle, draw, fill, minimum size=2, inner sep=0.02]
\node [b] (b2) at (0,0) {1};
\node [b] (b3) at (2,0) {2};
\draw[black,->,>=latex][decoration={markings,mark=at position .6 with {\arrow[scale=1.35,linkcolor]{latex}}},
 postaction={decorate},
 shorten >=0.4pt] (b2) to (b3);
\end{tikzpicture}}}+\Theta_3+\cdots=\sum_{p\geq0}\Upsilon^p_{\Theta_3}
\]
 in the graded Lie algebra $\big(\odfGC{3}{0}{2},\brdot\big)$~-- where the graph $\Upsilon^p_{\Theta_3}$ possesses $N=4p+2$ vertices and $k=6p+1$ edges~-- as the corresponding obstructions vanish, cf.\ footnote~\ref{footnote:recur}. Mapping the Maurer--Cartan element $\Upsilon_{\Theta_3}$ via the representation $\Rep^E$ yields an exotic \Lieinf-structure deforming non trivially the deformation complex $\big(\A^E_\text{\rm Lie},\pbdot^E_\symp\big)$ as a \Lieinf-algebra.
\end{proof}

The non-vanishing brackets of the (essentially unique\footnote{Up to gauge transformations and rescalings.}) exotic \Lieinf-structure of Theo\-rem~\ref{TheoremMain} -- denoted $\theta$ in the following~-- take the form\footnote{The intrinsic degree carried by each bracket is given by $|\theta_{4p+2}|=-12p-2$. Pulling back the brackets along the suspension map $s\colon\A^E_{\rm Lie}[2]\to\A^E_{\rm Lie}$ of degree 2 yields a series of brackets $\tilde \theta_{4p+2}$ on $\A^E_{\rm Lie}[2]$ with the usual degree~$-4p$.}
\begin{gather*}
\begin{split}
& \theta_2=\Rep^E_2(\Upsilon_\mathsf{S})=\pbdot^E_\symp ,\qquad \theta_6=\Rep^E_6(\Theta_3) ,\qquad \theta_{10}=\Rep^E_{10}\big(\Upsilon^2_{\Theta_3}\big) ,\qquad \dots ,\\ 
& \theta_{4p+2}=\Rep^E_{4p+2}\big(\Upsilon^p_{\Theta_3}\big) ,\qquad \dots.
\end{split}
\end{gather*}

The minimal\footnote{Recall that minimal \Lieinf-algebras are characterised by a vanishing differential $\theta_1\equiv0$.} \Lieinf-structure $(\A^E_{\rm Lie},\theta)$ can be interpreted as the avatar in dimension $d=3$ both of the Moyal star-commutator in $d=1$ and the Kontsevich--Shoikhet exotic \Lieinf-structure on infinite-dimensional manifolds in $d=2$. It relies on bi-oriented graphs and as such possesses no counterpart in the ``bialgebra'' realm where only one orientable direction is available. In fact, one can explicitly check that the first non-trivial deformed bracket $\theta_6=\Rep^E_6(\Theta_3)$ vanishes identically on the graded Poisson subalgebra\footnote{See footnote~\ref{footpoint}.} $\A^\alg_{\rm Lie}\subset\A^E_{\rm Lie}$ controlling deformations of Lie bialgebras, cf.\ Proposition~\ref{propob}.

\begin{Remark}\label{remarkQLB}\quad
\begin{itemize}\itemsep=0pt
\item In the next section, we will consider Maurer--Cartan elements in the formal extension $\A^E_{\rm Lie}\eps$ of $\A^E_{\rm Lie}$ by a formal parameter $\hbar$. By analogy with the $d=2$ case~\cite{Merkulov2016}, we will refer to Maurer--Cartan elements of $\big(\A^E_{\rm Lie}\eps,\theta\big)$ as \textit{formal ``quantizable Lie bialgebroids''}, to contrast with the Maurer--Cartan elements of $\big(\A^E_{\rm Lie}\eps,\pbdot^E_\symp\big)$ referred to simply as \textit{formal Lie bialgebroids}. Formal Lie bialgebroids being linear in $\hbar$ are just Lie bialgebroids and accordingly, we will refer to formal ``quantizable Lie bialgebroids'' linear in $\hbar$ as ``quantizable Lie bialgebroids''.
\item Note that ``quantizable Lie bialgebroids'' are in particular Lie bialgebroids as they satisfy $\pb{\cH}{\cH}_\symp^E=0$ on top of some higher consistency conditions \big($\theta_6\big(\cH^{\w 6}\big)=0$, etc.\big).
\end{itemize}
\end{Remark}

The distinction between Lie bialgebroids and ``quantizable Lie bialgebroids'' will become salient when applied to the quantization problem for Lie bialgebroids in Section~\ref{section:Application to quantization}.

\subsection{Application to quantization and future directions}
\label{section:Application to quantization}
The quantization problem for Lie-(quasi) bialgebras was formulated by V.~Drinfeld (cf.~\cite[Question~1.1]{Drinfeld1992} for Lie bialgebras and Section~5 for Lie-quasi bialgebras) and solved in \cite{Etingof1995} by Etingof--Kazhdan for the Lie bialgebra case\footnote{See also \cite{Ginot2016,Hinich2014,Merkulov2016,Merkulov2016a,Tamarkin2007b}.} and in \cite{Enriquez2008,Sakalos2013} for the Lie-quasi bialgebra case. In both cases, the solution is universal and involves the use of a Drinfeld associator, yielding an action of the Grothendieck--Teichm\"uller group $\GRT_1$ on the set of inequivalent universal quantization maps. The latter can be traced back to the action of $H^0(\odfcGC{3}{0}{1})\simeq\grt_1$ on the deformation complex of Lie and Lie-quasi bialgebras, as reviewed in Section \ref{section:Lie bialgebras2}. The situation for Lie bialgebras (and their quasi versions) is therefore much akin to the situation for finite-dimensional Poisson manifolds, in that both cases share the following important features (cf.\ the Introduction and Table~\ref{figgraphsymsum}):
\begin{enumerate}\itemsep=0pt
\item The Grothendieck--Teichm\"uller group plays a classifying r\^ole.
\item There is (conjecturally) no generic obstruction to the existence of universal quantizations.
\end{enumerate}
As for Lie bialgebroids, the corresponding quantization problem was formulated by P.~Xu in \cite{Xu1998a, Xu1999} as follows: Given a Lie bialgebroid structure on $(E,E^*)$, the associated quantum object is a topological deformation (called \textit{quantum groupoid}) of the standard (associative) bialgebroid structure on the universal enveloping algebra $\U_{R}(E)$ associated with the Lie algebroid structure on~$E$ \cite{Moerdijk2008}~-- with $R\equiv\fonc{\M}$~-- whose semi-classicalisation reproduces the original Lie bialgebroid structure. The quantization problem for Lie bialgebroids then consists in associating to each Lie bialgebroid a~quantum groupoid quantizing it. Although the quantization problem for a~generic Lie bialgebroid remains open, several explicit examples of quantizations for particular Lie bialgebroids have been exhibited in the literature. Apart from the above mentioned quantization of Lie bialgebras~\cite{Etingof1995}, it was shown in~\cite{Xu1999} that M.~Kontsevich's solution to the quantization problem for (finite-dimensional) Poisson manifolds \cite{Kontsevich:1997vb} ensures that Lie bialgebroids associated with Poisson manifolds (cf.\ Example~\ref{exLiebi}) constitute another example of Lie bialgebroids admitting a quantization.\footnote{More precisely, the Kontsevich star product $*$ quantizing the Poisson bivector $\pi$ provides a Drinfeld twistor $J_*\in \big(\U_R(E)\otimes_{R} \U_R(E)\big)\eps$ for the standard bialgebroid $\U_R(E)$, where $R\equiv\fonc{\M}$ and $E\equiv T\M$. Twisting $\U_R(E)$ by $J_*$ then provides a quantization of the Lie bialgebroid associated with $\pi$ on $(T\M,T^*\M)$, see~\cite{Xu1999}. } This result was shown to extend to regular triangular Lie bialgebroids in~\cite{Xu1999} (see also~\cite{Nest1996a}) using methods \`a la Fedosov~\cite{Fedosov1994} and to generic triangular Lie bialgebroids in \cite{Calaque2005} using a generalisation of Kontsevich's formality theorem for Lie algebroids. In this context, the following natural conjecture was formulated by P.~Xu:

\begin{Conjecture}[{Xu \cite[Section 6]{Xu1999}}] \label{conjXu}Every Lie bialgebroid admits a quantization as a quantum groupoid.
\end{Conjecture}

Although Conjecture~\ref{conjXu} might still hold true in the most general setting, we would like to argue for the non-existence of {\it universal}\footnote{Recall that a universal quantization admits formulae given by expansions in terms of graphs with universal coefficients.} quantizations of Lie bialgebroids, on the basis of the following results from Section~\ref{section:Lie bialgebroids2}:
\begin{enumerate}\itemsep=0pt
\item The Grothendieck--Teichm\"uller group plays no classifying r\^ole regarding the universal deformations (and hence quantizations) of Lie bialgebroids.
\item There exists a potential obstruction to the existence of universal quantizations of Lie bialgebroids.
\end{enumerate}
Contrasting these two features with their above mentioned counterparts for Lie bialgebras, one is led to conclude that the quantization problem for Lie bialgebroids differs essentially from its Lie bialgebra analogue and is in fact more akin to the quantization problem for infinite-dimensional manifolds. In view of this analogy, the obstruction appearing in Theorem \ref{TheoremMain} can be understood as the avatar in $d=3$ of the Kontsevich--Shoikhet obstruction in $d=2$. As shown in \cite{Dito2013,Penkava1998,Shoikhet2008a,Willwacher2013a}, the latter obstruction is hit in $d=2$ and thus prevents the existence of an {\it oriented} star product, thereby yielding a No-go result regarding the existence of universal quantizations for infinite-dimensional Poisson manifolds. Pursuing the analogy with the $d=2$ case, we conjecture the following:

\begin{Conjecture}[no-go]\label{congno}
There are no universal quantizations of Lie bialgebroids as quantum groupoids.
\end{Conjecture}

To explicitly show that the obstruction is hit would require a better understanding of the deformation theory of (associative) bialgebroids, which goes beyond the ambition of the present note.\footnote{While the deformation theory for associative algebras \cite{Gerstenhaber1963} and bialgebras \cite{Merkulov2007} are well understood using the frameworks of operads and properads respectively, the deformation theory for bialgebroids has for now been evading such (pr)operadic formulation. The underlying reason is that the relevant category to deal with bialgebroids is the one of $\big(R^e, R^e\big)$-bimodules~-- where $R$ is a ring and $R^e$ denotes its enveloping ring $R^e:=R\otimes R^\text{op}$~-- which is \textit{not} symmetric monoidal, as usually required to work with (pr)operads. Rather, the category of $\big(R^e, R^e\big)$-bimodules is naturally endowed with a structure of lax-oplax duoidal category (more precisely, a lax-strong duoidal category, that is the oplax structure is strong monoidal) whose bimonoids are bialgebroids, see, e.g.,~\cite{Basile2021}. As such, the deformation theory for bialgebroids cannot be described using the theory of properads (at least in its standard form). We are grateful to T.~Basile and D.~Lejay for clarifications regarding this fact. } We nevertheless conclude the present discussion by outlining a strategy of proof for Conjecture \ref{congno} by mimicking the two-steps procedure of \cite{Merkulov2016} for Poisson manifolds and Lie bialgebras and adapting it to the case at hand. Denoting $C^\bullet_{\mathsf{GS}}(\mathcal O_E,\mathcal O_E\big)$ the equivalent of the Gerstenhaber--Schack complex for the standard commutative co-commutative bialgebroid structure on the symmetric algebra $\mathcal O_E$ associated to the vector bundle $E$, the former should be endowed with a \Lieinf-algebra structure~$\mu$~-- generalising the one of \cite{Markl2004b,Merkulov2007} for the bialgebra case~-- whose corresponding Maurer--Cartan elements are quantum groupoids. By analogy with the bialgebra case, the cohomology of $\big(C^\bullet_{\mathsf{GS}}(\mathcal O_E,\mathcal O_E),\mu_1\big)$ should be isomorphic as a graded Lie algebra to the deformation complex $\big(\A^E_{\rm Lie},\pbdot^E_\symp\big)$ of Lie bialgebroids on $E$. Although these two \Lieinf-algebras should coincide in cohomology, we do not expect them to be quasi-isomorphic as \Lieinf-algebras, i.e., $\big(C^\bullet_{\mathsf{GS}}(\mathcal O_E,\mathcal O_E\big),\mu\big)$ is {\it not} formal. To show explicitly that the obstruction to formality is hit would require computing the \Lieinf-algebra structure obtained by transfer of $\mu$ on $H^\bullet\big(C^\bullet_{\mathsf{GS}}(\mathcal O_E,\mathcal O_E),\mu_1\big)$ and showing that the latter coincides with the exotic \Lieinf-structure~$\theta$ on~$\A^E_{\rm Lie}$,\footnote{This is only possible thanks to the fact that the exotic \Lieinf-structure~$\theta$ on~$\A^E_{\rm Lie}$ of Theorem~\ref{TheoremMain} is minimal and hence is a potential candidate for being the cohomology of another \Lieinf-structure.} as is the case in $d=2$ \cite{Merkulov2016, Shoikhet2008a}. This would provide a trivial (in the sense that no Drinfeld associator is needed) formality \Lieinf quasi-isomorphism $\big(\A^E_{\rm Lie},\theta\big)\iso\big(C^\bullet_{\mathsf{GS}}(\mathcal O_E,\mathcal O_E),\mu\big)$, yielding in turn a quantization map for (formal) ``quantizable Lie bialgebroids'' (the Maurer--Cartan elements of $\theta$ in $\A^E_{\rm Lie}\eps$, cf.\ Remark~\ref{remarkQLB}). Finally, the fact that $\theta$ is not \Lieinf-isomorphic to the big bracket in $\A^E_{\rm Lie}$ (cf.\ Theorem~\ref{TheoremMain}) would prevent the existence of a formality morphism for Lie bialgebroids. We sum up these (non)-formality conjectures for Lie bialgebroids in Figure~\ref{TableLie}.\looseness=-1

 \begin{figure}[ht]\centering\vspace{-7mm}
\begin{displaymath}\hspace*{20mm}
\begin{tikzcd}
{\big(\A^E_\text{\rm Lie},\pbdot^E_\symp\big)}\arrow[r, "\times", dashed]
 & \big(\A^E_\text{\rm Lie},\theta\big) \arrow[r, "\sim"] & {\big(C^\bullet_{\mathsf{GS}}(\mathcal O_E,\mathcal O_E),\mu\big)} \\[-3mm]
 {\substack{\text{Lie bialgebroids}}} \arrow[r, "\times", dashed] & {\substack{\text{``Quantizable} \\ \text{Lie bialgebroids''}}} \arrow[r, "\sim"] & {\substack{\text{Quantum} \\ \text{groupoids}}}
\end{tikzcd}
\end{displaymath}\vspace{-5mm}

\caption{Conjectural (non)-formality maps for Lie bialgebroids.} \label{TableLie}
\end{figure}

Note that the situation is markedly different in the Lie-quasi (and quasi-Lie) bialgebroid case.\footnote{For definiteness, we will focus on the Lie-quasi case, keeping in mind that the arguments apply similarly to the dual case.} Firstly, recall that the Maurer--Cartan element $\Upsilon_{\Theta_3}$ is {\it not} gauge-related to the Maurer--Cartan element $\Upsilon_\mathsf{S}$ in $\odfcGC{3}{0}{2}$~-- since $\Theta_3$ is a non-trivial cocycle in $\odfcGC{3}{0}{2}$~-- hence $\theta$ is indeed a non-trivial deformation of the big bracket in~$\A^E_{\rm Lie}$. However, the cocycle $\Theta_3$ {\it is} a coboundary in~$\odfcGC{3}{1}{1}$ (cf.\ Appendix~\ref{section:Incarnation}) so that there exists a combination of graphs $\vartheta_3\in\odfcGC{3}{1}{1}$ such that $\Theta_3=-\delta\vartheta_3\in\odfcGC{3}{0}{2}$. In order for $\Upsilon_{\Theta_3}$ and $\Upsilon_\mathsf{S}$ to be gauge-related in $\odfcGC{3}{1}{1}$, one needs to find a degree $0$ element\footnote{The graph $\vartheta_3^p$ has $N=4p+1$ vertices and $k=6p$ edges so that to have degree 0 in $d=3$.} $\vartheta=\vartheta_3+\vartheta_3^2+\cdots+\vartheta_3^p$ of $\odfcGC{3}{1}{1}$ such that $\Upsilon_{\Theta_3}=e^{\ad\, \vartheta}\Upsilon_\mathsf{S}$. Contrarily to the problem of prolongating the cocycle~$\Theta_3$ to the Maurer--Cartan element $\Upsilon_{\Theta_3}$~-- which can be solved by a trivial induction~-- to display an explicit gauge map $\vartheta$ is a highly non-trivial task\footnote{We refer to \cite{Merkulov2016} for an explicit construction in $d=2$. } as the higher obstructions\footnote{The first obstruction vanishes since $\Theta_3$ is exact in $\odfcGC{3}{1}{1}$. The second obstruction $\Upsilon^2_{\Theta_3}-\half\br{\vartheta_3}{\Theta_3}$ can be checked to live in $H^1(\odfcGC{3}{1}{1})$.} live in $H^1(\odfcGC{3}{1}{1})\simeq H^1(\fcGC_2)$, i.e., the recipient of the obstructions to the universal quantization of finite-dimensional Poisson manifolds. Although this cohomological space conjecturally vanishes\footnote{See Table~\ref{Tabcohom}.} (Drinfeld--Kontsevich), maps allowing to convert cocycles into coboundaries are highly non-trivial and necessarily involve the choice of a Drinfeld associator (consistently with the fact that two coboundaries differ by the choice of an element in $H^0(\fcGC_2)\simeq\grt_1$). Up to the Drinfeld--Kontsevich conjecture, it is nevertheless expected that, given a Drinfeld associator, one can define a \Lieinf-isomorphism $\big(\A^E_\text{\rm Lie-quasi},\pbdot^E_\symp\big)\iso \big(\A^E_\text{\rm Lie-quasi},\theta\big)$. Repeating the argument laid down in the Lie bialgebroid case, one needs to find a (trivial) \Lieinf quasi-isomorphism $\big(\A^E_\text{Lie-quasi},\theta\big)\iso\big(C^\bullet_{\text{\rm quasi}-\mathsf{GS}}(\mathcal O_E,\mathcal O_E),\mu\big)$, where the right-hand side stands for the deformation complex of the quantum object associated to Lie-quasi bialgebroids. The relevant category here is the one of \textit{quasi-bialgebroids}, a common generalisation of the notions of bialgebroids and quasi-bialgebras, allowing to define the associated notion of \textit{quasi-quantum groupoid} as topological deformation of the standard bialgebroid $\U_R(E)$ as a quasi-bialgebroid. Composing with the (non-trivial) \Lieinf-isomorphism $\big(\A^E_\text{\rm Lie-quasi},\pbdot^E_\symp\big)\iso \big(\A^E_\text{\rm Lie-quasi},\theta\big)$ would yield a formality morphism for Lie-quasi bialgebroids. We sum up the discussion by formulating the following conjecture and recap the corresponding conjectural formality maps in the Lie-quasi case in Figure~\ref{TablequasiLie}.
\begin{Conjecture}[yes-go]\label{congyes}
Given a Drinfeld associator, one can define a universal quantization of Lie-quasi bialgebroids as quasi-quantum groupoids.
\end{Conjecture}
 \begin{figure}[ht]\centering\vspace*{-7mm}
\begin{displaymath}
\hspace*{5mm}
\begin{tikzcd}
{\big(\A^E_\text{\rm Lie-quasi},\pbdot^E_\symp\big)}\arrow{r}{\sim}[swap]{{\circlearrowleft_{\text{\rm GRT}_1}}}
 & \big(\A^E_\text{\rm Lie-quasi},\theta\big) \arrow[r, "\sim"] & {\big(C^\bullet_{\text{\rm quasi}-\mathsf{GS}}(\mathcal O_E,\mathcal O_E),\mu\big)} \\[-3mm]
 {\substack{\text{Lie-quasi} \\ \text{bialgebroids}}} \arrow{r}{\sim}[swap]{{\circlearrowleft_{\text{\rm GRT}_1}}} & {\substack{\text{``Quantizable} \\ \text{Lie-quasi bialgebroids''}}} \arrow[r, "\sim"] & {\substack{\text{Quasi-quantum} \\ \text{groupoids}}}
\end{tikzcd}
\end{displaymath}\vspace{-5mm}

 \caption{Conjectural formality maps for Lie-quasi bialgebroids.} \label{TablequasiLie}
 \end{figure}

 A dual conjecture can be made about quasi-Lie bialgebroids. Of particular interest in this context would be to investigate the quantization of quasi-Lie bialgebroids associated to twisted Poisson manifolds (cf.\ Example~\ref{example:TWP}) by means of a Drinfeld twistor~\cite{Xu1999} for the corresponding non-associative Kontsevich star product~\cite{Cornalba2001}.

To conclude, let us note that particularising the conjectural quantization map of Conjecture~\ref{congyes} to Lie bialgebroids ensures that every Lie bialgebroid admits a quantization {\it as a quasi-quantum groupoid} (but generically not as a quantum groupoid as stated in Conjecture~\ref{conjXu}). However, for the particular subclass of ``quantizable Lie bialgebroids'' (which are in particular Lie bialgebroids, see Remark~\ref{remarkQLB})~-- such as Lie bialgebras and coboundary Lie bialgebroids (cf.\ Appendix~\ref{section:Incarnation})~-- the associated quantization should yield a quantum groupoid. According to the above picture, the exotic \Lieinf-structure $\theta$ of Theorem~\ref{TheoremMain} can therefore be seen as a concrete means to delineate the subclass of Lie bialgebroids susceptible to be quantized as quantum groupoids.

\appendix
\section{Geometry of Lie bialgebroids}\label{section:Appendix}
{\bf Lie bialgebras.}
A Lie bialgebra is a vector space $\alg$ endowed with a Lie algebra structure on both $\alg$ and its dual $\alg^*$ such that the cobracket $\Delta_\alg\colon \alg\to\w^2\alg$ is a cocycle for the Lie algebra $(\alg,\brdot_\alg)$, i.e., $\Delta_\alg(\br{x}{y})=\ad_x\Delta_\alg(y)-\ad_y\Delta_\alg(x)$ where the representation used is the extension $\ad:\alg\otimes(\w^2\alg)\to \w^2\alg$ of the adjoint action of $(\alg,\brdot_\alg)$ on $\w^2\alg$ as $\ad_x (y\w z)=\br{x}{y}_\alg\w z+y\w\br{x}{z}_\alg$. Letting $\pset{e_a}|_{a\in\pset{1,\dots, \dim \alg}}$ be a basis of $\alg$, one denotes $\br{e_a}{e_b}_\alg=f_{ab}{}^ce_c$ and $\Delta_\alg (e_c)=C_c{}^{a b}e_a\otimes e_b$. The three defining conditions of a Lie bialgebra read
\begin{gather}
f_{e[a}{}^{d}f_{bc]}{}^e=0 ,\qquad C_d{}^{e[a}C_e{}^{bc]}=0 ,\qquad f_{ab}{}^eC_e{}^{cd}-4 f_{e[a}{}^{[c}C_{b]}{}^{d]e}=0, \label{combi}
\end{gather}
where $f_{a b}{}^c=f_{[a b]}{}^c$ and $C_c{}^{a b}=C_c{}^{[a b]}$ where square brackets denote skewsymmetrisation.

{\bf Lie algebroids.} A Lie algebroid is a triplet $(E,\rho,\brdot_E)$ where:
\begin{itemize}\itemsep=0pt
\item $E\overset{\pi}{\to} \M$ is a vector bundle over the manifold $\M$,
\item $\rho\colon E\to T\M$ is a morphism of vector bundles called the \textit{anchor},\footnote{We will use the same symbol $\rho$ to denote the induced map of sections $\rho\colon \fE\to\Gamma(T\M)$.}
\item $\brdot_E \colon \Gamma(E)\otimes \Gamma(E)\to \Gamma(E)$ is a $\corps$-bilinear map called the \textit{bracket},
\end{itemize}
such that the following conditions are satisfied for all $f\in\fonc{\M}$ and $X,Y,Z\in\Gamma(E)$:
\begin{enumerate}\itemsep=0pt
\item \textit{skewsymmetry}: $\br{X}{Y}_E=-\br{Y}{X}_E$,
\item \textit{Leibniz rule}: \hspace{0.57cm}$\br{X}{f\cdot Y}_E=\rho_X[f]\cdot Y+f\cdot \br{X}{Y}_E$,
\item \textit{Jacobi identity}: \hspace{0.05cm}$\br{X}{\br{Y}{Z}_E}_E+\br{Y}{\br{Z}{X}_E}_E+\br{Z}{\br{X}{Y}_E}_E=0$.
\end{enumerate}
The previous conditions ensure that the map $\rho\colon \Gamma(E)\to\Gamma(T\M)$ defines a morphism of Lie algebras between the Lie algebra $(\Gamma(E),\brdot_E)$ and the Lie algebra of vector fields on~$\M$, i.e., $\rho_{\br{X}{Y}_E}=\br{\rho_X}{\rho_Y}$ for all $X,Y\in\Gamma(E)$.

\begin{Proposition}
Let $(E,\rho,\brdot_E)$ be a Lie algebroid.
The following statements hold:
\begin{enumerate}\itemsep=0pt
\item[$1.$] Let $\pset{x^\mu}|_{\mu\in\pset{1,\dots, \dim\M}}$ be a set of coordinates of~$\M$ and $\pset{e_a}|_{a\in\pset{1,\dots, \dim E}}$ be a basis of~$\Gamma(E)$. Setting $\rho_{e_a}[f]=\rho_a{}^\mu(x) \p_\mu f$ and $\br{e_a}{e_b}_E=f_{ab}{}^c(x) e_c$, the defining conditions of a Lie algebroid can be expressed in components as
\begin{displaymath}
f_{ab}{}^c=-f_{ba}{}^{c} ,\qquad 2 \rho_{[a}{}^\nu\p_\nu\rho_{b]}{}^\mu=\rho_c{}^\mu f_{ab}{}^c ,\qquad\rho_{[c}{}^\nu\p_\nu f_{ab]}{}^d=f_{e[c}{}^df_{ab]}{}^e.
\end{displaymath}
Acting on generic sections of $E$, the Lie algebroid bracket reads
\[
\br{X}{Y}_E=\big(\rho_X[Y^c]-\rho_Y[X^c]+f_{ab}{}^cX^aY^b\big)e_c.
\]
\item[$2.$] The exterior algebra $\Gamma(\w^\bullet E^*)$ is naturally endowed with a structure of dg commutative algebra with differential ${\rm d}_E\colon \Gamma(\w^\bullet E^*)\to\Gamma(\w^{\bullet+1} E^*)$ defined by
\begin{itemize}\itemsep=0pt
\item $({\rm d}_E f)(X)=\rho_X[f]$,
\item $({\rm d}_E \eta)(X,Y)=\rho_X[\eta(Y)]-\rho_Y[\eta(X)]-\eta(\br{X}{Y}_E)$,
\item ${\rm d}_E(\alpha\w\beta)=({\rm d}_E\alpha)\w\beta+(-1)^{|\alpha|} \alpha\w ({\rm d}_E\beta)$,
\end{itemize}
for all $X,Y\in\Gamma(E)$, $f\in\fonc{\M}$, $\eta\in\Gamma(E^*)$ and $\alpha,\beta\in\Gamma(\w^\bullet E^*)$.
\item[$3.$] The dual exterior algebra $\Gamma(\w^\bullet E)$ is endowed with a structure of Gerstenhaber algebra with graded bracket $\pbdot_E\colon \Gamma(\w^\bullet E)\otimes\Gamma(\w^\circ E)\to \Gamma(\w^{\bullet+\circ-1} E)$ defined as follows
\begin{gather*}
\pb{f}{g}_E=0,\qquad\pb{X}{f}_E=\rho_X[f],\qquad \pb{X}{Y}_E=\br{X}{Y}_E,\\
\pb{P}{Q}_E=-(-1)^{(|p|-1)(|q|-1)}\pb{Q}{P}_E,\\
\pb{P}{Q\w R}_E=\pb{P}{Q}_E\w R+(-1)^{|q|(|p|-1)}Q\w \pb{P}{R}_E,
\end{gather*}
for all $f,g\in\fonc{\M}$, $X,Y\in\Gamma(E)$ and $P,Q,R\in\Gamma(\w^\bullet E)$.

These conditions can be checked to ensure the graded Jacobi identity:
\begin{gather*}
\pb{\pb{P}{Q}_E}{R}_E+(-1)^{(|P|-1)(|Q|+|R|)}\pb{\pb{Q}{R}_E}{P}_E\\
\qquad{} +(-1)^{(|R|-1)(|P|+|Q|)}\pb{\pb{R}{P}_E}{Q}_E=0.
\end{gather*}
\end{enumerate}
\end{Proposition}

\begin{Example}\label{exLiebi}\quad
\begin{itemize}\itemsep=0pt
\item A Lie algebra is a Lie algebroid whose base manifold is a point.
\item Given a manifold $\M$, the \textit{standard Lie algebroid} is defined as the tangent bundle $T\M$ together with the identity map as anchor and the usual Lie bracket of vector fields as bracket. The corresponding differential on the space of differential forms $\Om^\bullet(\M)\simeq \Gamma(\w^\bullet T^*\M)$ coincides with the \textit{de Rham differential} while the induced Gerstenhaber bracket on $\Gamma(\w^\bullet T\M)$ identifies with the \textit{Schouten bracket} on polyvector fields.
\item Let $(\M,\pi)$ be a Poisson manifold. The dual tangent bundle $T^*\M$ is naturally endowed with a Lie algebroid structure with anchor $\pi^\sharp\colon \Gamma(T^*\M)\to\Gamma(T\M)\colon \alpha_\mu {\rm d}x^\mu\mapsto \pi^{\mu\nu}\alpha_\nu\p_\mu$ and bracket $\br{\alpha}{\beta}_\pi=\half\big(\Lag_{\pi^\sharp(\alpha)}\beta-\Lag_{\pi^\sharp(\beta)}\alpha+i_{\pi^\sharp(\alpha)}{\rm d}\beta -i_{\pi^\sharp(\beta)}{\rm d}\alpha\big)\text{ for all }\alpha,\beta\in\Gamma(T^*\M)$.
\end{itemize}
\end{Example}

{\bf Lie bialgebroids.} The concept of Lie bialgebroid was introduced by Mackenzie--Xu in \cite{Mackenzie1994} as the infinitesimal variant of a Poisson groupoid. We follow the modern definition of~\cite{Kosmann-Schwarzbach1995} and define a Lie bialgebroid as a vector bundle $E\overset{\pi}{\to} \M$ endowed with two dual Lie algebroid structures satisfying a natural compatibility condition. Denoting $(\rho,\brdot_E)$ the Lie algebroid structure on~$E$ and $(R,\brdot_{E^*})$ the one on $E^*$, the pair $(E,E^*)$ is a Lie bialgebroid if ${\rm d}_{E^*}$ is a~derivation of $\pbdot_E$, where we denoted $\pbdot_E$ the Gerstenhaber bracket on $\Gamma(\w^\bullet E)$ induced by $(\rho,\brdot_E)$ and ${\rm d}_{E^*}$ the differential on $\Gamma(\w^\bullet E)$ induced by $(R,\brdot_{E^*})$.\footnote{Note that the defining condition of a~Lie bialgebroid can equivalently be stated as the fact that the differential ${\rm d}_{E}$ on $\Gamma(\w^\bullet E^*)$ induced by $(\rho,\brdot_E)$ is a derivation of the Gerstenhaber bracket~$\pbdot_{E^*}$ induced by $(R,\brdot_{E^*})$, i.e., the notion of Lie bialgebroid is self-dual.}

\begin{Example}\label{example:Lie bialgebroids}\quad
\begin{itemize}\itemsep=0pt
\item A Lie bialgebra is a Lie bialgebroid whose base manifold is a point.
\item Letting $\M$ be a Poisson manifold, the two Lie algebroid structures on $T\M$ and $T^*\M$ as defined in Example~\ref{exLiebi} are compatible in the above sense and hence define a Lie bialgebroid structure on~$(T\M,T^*\M)$.
\end{itemize}
\end{Example}

Lie bialgebroids thus generalise both Lie bialgebras and Poisson manifolds. In fact, the base manifold $\M$ of any Lie bialgebroid is endowed with a canonical Poisson bracket\footnote{More generally, for a proto-Lie bialgebroid, the Jacobi identity for the bivector is deformed as \begin{displaymath}
\pi^{\rho[\lambda} \p_{\rho}\pi^{\mu \nu]} = \tfrac{1}{3}R^{a [\lambda} R^{b| \mu} R^{c| \nu]} \psi_{a b c}+\tfrac{1}{3}\rho_{a}{}^{[\lambda} \rho_{b}{}^{\mu} \rho_{c}{}^{\nu]} \varphi^{a b c}.
\end{displaymath}\vspace*{-5mm}}
defined as $\pb{f}{g}=\dl {\rm d}_{E^*}f, {\rm d}_{E}g\dr$ \big(or in components as $\pb{f}{g}=R^{a[\mu}\rho_a{}^{\nu]}\p_\mu f\p_\nu g$\big).

The exterior algebra $\big(\Gamma(\w^\bullet E),\w\big)$ of a Lie bialgebroid is endowed with both a structure of Gerstenhaber bracket $\pbdot_E$ and of a differential $d_{E^*}$ being a derivation for both the graded commutative product $\w$ and the Gerstenhaber bracket. Such a quadruplet $\big(\Gamma(\w^\bullet E),\w,{\rm d}_{E^*},\pbdot_E\big)$ is called a strong differential Gerstenhaber algebra. It was in fact shown in~\cite{Xu1997} (see also \cite{Kosmann-Schwarzbach1995}) that Lie bialgebroid structures on a vector bundle $E\overset{\pi}{\to} \M$ are in one-to-one correspondence with strong differential Gerstenhaber structures on $\big(\Gamma(\w^\bullet E),\w\big)$ \big(or equivalently with dg Poisson structures on $(\Gamma(\w^\bullet E)[1],\w)$ \big).

{\bf Coboundary Lie bialgebroids.}
Letting $E$ be a Lie algebroid, an \textit{$r$-matrix} is a section $\Lambda\in\Gamma\big({\w}^2E\big)$ satisfying $\pb{X}{\pb{\Lambda}{\Lambda}_E}_E=0$ for all $X\in\Gamma(E)$. An {$r$-matrix} endows~$E$ with a~structure of Lie bialgebroid by defining $d_{E^*}=\pb{\Lambda}{\cdot}_E$. The defining condition on $\Lambda$ is necessary and sufficient to ensure that the inner derivation $d_{E^*}$ squares to zero. A Lie bialgebroid defined in this way is called a coboundary Lie bialgebroid. Whenever the stronger condition $\pb{\Lambda}{\Lambda}_E=0$ holds, the induced Lie bialgebroid is said to be \textit{triangular}~\cite{Mackenzie1994} (cf.\ for example the Lie bialgebroid on Poisson manifolds defined in Example~\ref{exLiebi}). If furthermore, $\Lambda$ is of constant rank, the triangular Lie bialgebroid is said to be \textit{regular}.

{\bf Quasi-Lie, Lie-quasi and proto-Lie bialgebroids.} The above mentioned characterisation of Lie bialgebroids as dg Poisson structures on $\big(\Gamma(\w^\bullet E)[1],\w\big)$ calls for several natural generalisations, as summarised in the following table\footnote{As is transparent from the correspondence below, there is a series of inclusions of bialgebroids: Lie $\subset$ quasi-Lie $\subset$ proto-Lie and Lie $\subset$ Lie-quasi $\subset$ proto-Lie. Note furthermore that the notions of Lie bialgebroids and proto-Lie bialgebroids are self-dual \big(and thus can be defined on both $\Gamma(\w^\bullet E)[1]$ and $\Gamma(\w^\bullet E^*)[1]$\big) while the notions of Lie-quasi bialgebroids and quasi-Lie bialgebroids are dual to each other.} (see, e.g., \cite{Antunes2019,Fregier2015, Kosmann-Schwarzbach,Roytenberg2001}):
\begin{itemize}\itemsep=0pt
\item Lie bialgebroids on $(E,E^*)$ $\Leftrightarrow$ dg Poisson algebras on $\Gamma(\w^\bullet E)[1]$,
\item quasi-Lie bialgebroids on $(E,E^*)$ $\Leftrightarrow$ homotopy Poisson algebras\footnote{Recall that a (curved) homotopy Poisson structure on a graded commutative algebra $(\alg,\w)$ is a (curved) \Lieinf-structure on $\alg$ such that all brackets are multi-derivations with respect to $(\alg,\w)$, see, e.g., \cite{Lang2013, Mehta2007}. A dg Poisson algebra is thus a (flat) homotopy Poisson algebra for which the brackets of arity above 2 vanish.} on $\Gamma(\w^\bullet E)[1]$,
\item Lie-quasi bialgebroids on $(E,E^*)$ $\Leftrightarrow$ homotopy Poisson algebras on $\Gamma(\w^\bullet E^*)[1]$,
\item proto-Lie bialgebroids on $(E,E^*)$ $\Leftrightarrow$ curved homotopy Poisson algebras on $\Gamma(\w^\bullet E)[1]$.
\end{itemize}

We now focus on the most general case, namely proto-Lie bialgebroids. Apart from the usual data (two anchors $\rho$, $R$ and two brackets $\brdot_E$, $\brdot_{E^*})$, a proto-Lie bialgebroid contains two elements $\varphi\in\Gamma\big({\w}^3E\big)$ and $\psi\in\Gamma\big({\w}^3 E^*\big)$ which play the r\^ole of various obstructions to the usual Lie bialgebroid identities. Letting $E^*[1]$ be the shifted dual bundle coordinatised by $\{x^\mu,\zeta_a\}$ of respective degree $0$ and~$1$, the graded commutative algebra $\big(\fonc{E^*[1]},\cdot\big)$ is isomorphic to the exterior algebra of sections $\big(\Gamma(\w^\bullet E),\w\big)$. The most general curved homotopy Poisson structure~$l$ on $\Gamma(\w^\bullet E)[1]$ thus takes the form\footnote{\label{footnotedegree}The bracket $l_p$ being of degree $3-2p$ on $\Gamma(\w^\bullet E)\simeq\fonc{E^*[1]}$, the fact that each bracket is a multi-derivation for the underlying graded commutative algebra constrains all brackets of arity higher than 3 to vanish. Pulling back the brackets along the suspension map $s\colon \Gamma(\w^\bullet E)[1]\to\Gamma(\w^\bullet E)$ of degree 1 yields a series of brackets on $\Gamma(\w^\bullet E)[1]$ with the usual degree $2-p$.}
\begin{gather*}
\bullet \ l_0:=\frac16\varphi^{abc}\zeta_a\zeta_b\zeta_c,\\
\bullet \ l_1(f):={\rm d}_{E^*}f=R^{a|\mu}\zeta_a\frac{\p f}{\p x^\mu}-\frac{1}{2} C_c{}^{ab}\zeta_a\zeta_b\frac{\p f}{\p \zeta_c},\\
\bullet \ l_2(f,g):=\pb{f}{g}_E=-\rho_a{}^\mu\left( \frac{\p f}{\p \zeta_a}\frac{\p g}{\p x^\mu}+(-1)^{|f|} \frac{\p f}{\p x^\mu}\frac{\p g}{\p \zeta_a}\right)+(-1)^{|f|}f_{ab}{}^{c}\zeta_c\frac{\p f}{\p \zeta_a}\frac{\p g}{\p \zeta_b},\\
\bullet \ l_3(f,g,h):=(-1)^{|g|}\psi_{abc}\frac{\p f}{\p \zeta_a}\frac{\p g}{\p \zeta_b}\frac{\p h}{\p \zeta_c}.
\end{gather*}
Imposing the defining quadratic condition $\br{l}{l}_\mathsf{NR}=0$ of a (curved) \Lieinf-algebra yields a series of identities which precisely reproduce the components conditions \eqref{consHam1}--\eqref{consHam9} as
\begin{gather}
\bullet \ \br{l_0}{l_1}_\mathsf{NR}=0 \ \Leftrightarrow \ \mathcal C_8=0,\label{curvproto1}\\
\bullet \ \br{l_0}{l_2}_\mathsf{NR}+\tfrac{1}{2} \br{l_1}{l_1}_\mathsf{NR}=0 \ \Leftrightarrow \ \mathcal C_3=\mathcal C_4=0, \label{curvproto2}\\
\bullet \ \br{l_0}{l_3}_\mathsf{NR}+\br{l_1}{l_2}_\mathsf{NR}=0 \ \Leftrightarrow \ \mathcal C_5=\mathcal C_6=\mathcal C_7=0, \label{curvproto3}\\
\bullet \ \br{l_1}{l_3}_\mathsf{NR}+\tfrac{1}{2} \br{l_2}{l_2}_\mathsf{NR}=0 \ \Leftrightarrow \ \mathcal C_1=\mathcal C_2=0, \label{curvproto4}\\
\bullet \ \br{l_2}{l_3}_\mathsf{NR}=0 \ \Leftrightarrow \ \mathcal C_9=0. \label{curvproto5}
\end{gather}
Imposing $\psi\equiv0$ (resp. $\varphi\equiv0$) yields a Lie-quasi (resp.\ quasi-Lie) bialgebroid and Lie bialgebroids are recovered by setting $\psi\equiv0$, $\varphi\equiv0$. Assuming that the base manifold $\M$ is the one-point manifold and denoting the vector space $\Gamma(E)$ as $\alg$ allows to define the counterparts of these notions in the ``bialgebra'' realm, cf.\ Section~\ref{section:Lie bialgebras}.

{\bf Twisting.} The previous formulation allows us to introduce the notion of \textit{twist} of a proto-Lie bialgebroid in terms of twist of (curved) \Lieinf-algebras.\footnote{\label{footnoteGetzler}The twisting procedure has been introduced by Quillen~\cite{Quillen1969} for dg Lie algebras and later generalised by Getzler \cite{Getzler2004} to \Lieinf-algebras.
Letting $(\alg,l)$ be a nilpotent (curved) \Lieinf-algebra and $\mathfrak m\in\alg^1$ be an arbitrary element of degree~$1$, one defines the \textit{twisted brackets} $l^{\mathfrak m}_n\colon \alg^{\w n}\to\alg$ of degree $2-n$ for $n\geq0$ as
 \begin{displaymath}
l^{\mathfrak m}_n(a_1,\dots,a_n):=\sum_{k\geq0}\frac{1}{k!} l_{k+n}(\mathfrak m^{\w k},a_1,\dots,a_n),\qquad \text{where} \quad a_i\in\alg.
 \end{displaymath}
The twisted brackets can be checked to endow $\alg$ with a structure of \textit{curved} \Lieinf-algebra. Whenever the curvature $l^{\mathfrak m}_0\in\alg^2$ vanishes, the element $\mathfrak m\in\alg^1$ is called a \textit{Maurer--Cartan element} and $(\alg,l^{\mathfrak m})$ is then a \textit{flat} \Lieinf-algebra. }
Given a proto-Lie bialgebroid on $(E,E^*)$ defined by the curved homotopy Poisson algebra $\big(\Gamma(\w^\bullet E)[1],\w, l\big)$ where $ l:=\pset{ l_p}_{0\leq p \leq 3}$, one can define a new proto-Lie bialgebroid structure on $(E,E^*)$ via twisting the (curved) \Lieinf-algebra structure $ l$ by an arbitrary bivector $\Lambda\in\Gamma\big({\w}^2E\big)$ \big(being the most general element of degree $1$ in $\Gamma(\w^\bullet E)[1]$\big). In components, such a twist amounts to the following shift of the components of the proto-Lie bialgebroid:
\begin{gather}
\rho_a{}^\mu\overset{\Lambda}{\longrightarrow}\rho_a{}^\mu,\qquad
f_{a b}{}^c\overset{\Lambda}{\longrightarrow} f_{a b}{}^c+\psi_{abe}\Lambda^{ec},\qquad
\psi_{abc}\overset{\Lambda}{\longrightarrow}\psi_{abc},\nonumber\\
R^{a|\mu}\overset{\Lambda}{\longrightarrow} R^{a|\mu}-\rho_b{}^\mu\Lambda^{ba},\qquad
C_c{}^{a b}\overset{\Lambda}{\longrightarrow} C_c{}^{a b}+\rho_c{}^\mu\p_\mu\Lambda^{ab}+2 \Lambda^{e[a}f_{ec}{}^{b]}-\psi_{cef}\Lambda^{ea}\Lambda^{fb},
\nonumber\\
\varphi^{abc}\overset{\Lambda}{\longrightarrow}\varphi^{abc}-3 R^{[c|\mu}\p_\mu\Lambda^{ab]}+3 \Lambda^{e[c}C_e{}^{ab]}+3 \rho_d{}^\mu \Lambda^{d[a}\p_\mu\Lambda^{bc]}\nonumber\\
\hphantom{\varphi^{abc}\overset{\Lambda}{\longrightarrow}}{}
-3 f_{ef}{}^{[c}\Lambda^{e|a}\Lambda^{f|b]}-\psi_{def}\Lambda^{da}\Lambda^{eb}\Lambda^{fc}. \label{twistLambdaexp}
\end{gather}
The latter can be checked to form the components of a proto-Lie bialgebroid without extra assumption on the bivector $\Lambda$. Note that, while a twist by a generic $\Lambda$ does preserve the subspace of Lie-quasi bialgebroids (for which $\psi\equiv0$), only bivectors satisfying the Maurer--Cartan equation $ l_1(\Lambda)+\half l_2(\Lambda,\Lambda)+\frac{1}{6} l_3(\Lambda,\Lambda,\Lambda)=0$ in $\Gamma(\w^\bullet E)[1]$ preserve the subspace of quasi-Lie bialgebroids\footnote{Equivalently, the Maurer--Cartan condition ensures that the twisting preserves the flatness of the \Lieinf-algebra associated to a given quasi-Lie bialgebroid, see footnote \ref{footnoteGetzler}.} (for which $\varphi\equiv0$), see~\cite{Liu1995} for the Lie bialgebroid case and \cite{Roytenberg2001} for quasi-Lie bialgebroids.

\begin{Example}[twisted Poisson manifolds]\label{example:TWP}
Let $\M$ be a manifold and $H\in\Om^3(\M)$ be a closed 3-form.
There is a quasi-Lie bialgebroid structure
on $(T\M,T^*\M)$ defined in components as
\begin{gather*}
\rho_\nu{}^\mu=\delta_\nu{}^\mu,\qquad f_{\mu \nu}{}^\lambda=0 ,\qquad R^{\nu|\mu}=0,\qquad C_\lambda{}^{\mu \nu}=0 ,\qquad \psi_{\lambda \mu \nu}=H_{\lambda \mu \nu}.
\end{gather*}
\noindent Letting $\pi\in\Gamma\big({\w}^2T\M\big)$ be a bivector, the components of the twisted proto-Lie bialgebroid read
\begin{gather*}
\tilde\rho_\nu{}^\mu=\delta_\nu{}^\mu
,\qquad
\tilde f_{\mu \nu}{}^\lambda=H_{\mu \nu \gamma}\pi^{\gamma\lambda}
,\qquad
\tilde R^{\nu|\mu}=-\pi^{\mu \nu},\qquad
\tilde C_\lambda{}^{\mu \nu}=\p_\lambda\pi^{\mu\nu}-H_{\lambda\rho\sigma}\pi^{\rho\mu}\pi^{\sigma\nu},\\
 \tilde\varphi^{\lambda\mu \nu}=3\, \pi^{\rho[\lambda}\p_\rho\pi^{\mu\nu]}-H_{\rho\sigma\gamma}\pi^{\rho\lambda}\pi^{\sigma\mu}\pi^{\gamma\nu},\qquad
\tilde\psi_{\lambda \mu \nu}=H_{\lambda \mu \nu}.
\end{gather*}
The resulting proto-Lie bialgebroid is again a quasi-Lie bialgebroid if and only $\tilde\varphi\equiv0$, that is, if $\pi$ is a {twisted Poisson bivector}~\cite{Roytenberg2001}.\footnote{Recall that \textit{twisted Poisson bivectors} \cite{Klimcik2002,Severa2001} are bivectors satisfying the twisted Jacobi identity $\br{\pi}{\pi}_{{\mathsf S}}=\frac{1}{3}H(\pi,\pi,\pi)$.}
Whenever $H$ vanishes, we recover the Lie bialgebroid structure on the (co)-tangent bundle of a Poisson manifold, cf.\ Example~\ref{example:Lie bialgebroids}.
\end{Example}

Dually, one can consider twisting the curved homotopy Poisson algebra structure on
$\Gamma({\w}^\bullet\! E^*\!)[1]\!$ by a 2-form field $\om\in\Gamma\big({\w}^2 E^*\big)$ which amounts to the following shift of the components of the proto-Lie bialgebroid
\begin{gather}
\rho_a{}^\mu\overset{\om}{\longrightarrow}\rho_a{}^\mu-R^{b|\mu}\om_{ba},\qquad
f_{a b}{}^c\overset{\om}{\longrightarrow}f_{a b}{}^c+R^{c|\mu}\p_\mu\om_{ab}-2 \om_{d[a}C_{b]}{}^{cd}-\varphi^{cde}\om_{da}\om_{eb},\nonumber\\
R^{a|\mu}\overset{\om}{\longrightarrow} R^{a|\mu},\qquad C_c{}^{a b}\overset{\om}{\longrightarrow} C_c{}^{a b}+\varphi^{abd}\om_{dc},\qquad
\varphi^{abc}\overset{\om}{\longrightarrow}\varphi^{abc},\nonumber\\
\psi_{abc}\overset{\om}{\longrightarrow}\psi_{abc}-3 \rho_{[a}{}^{\mu}\p_\mu\om_{bc]}+3 \om_{e[c}f_{ab]}{}^e+3\, R^{d|\mu}\om_{d[a}\p_\mu\om_{bc]}\nonumber\\
\hphantom{\psi_{abc}\overset{\om}{\longrightarrow}}{}
-3 C_{[c}{}^{de}\om_{d|a}\om_{e|b}-\varphi^{def}\om_{da}\om_{eb}\om_{fc}.\label{twistomexp}
\end{gather}
Twisting by a 2-form field does preserve quasi-Lie bialgebroids but fails to preserve Lie-quasi bialgebroids unless $\om$ satisfies the associated Maurer--Cartan equation on $\Gamma(\w^\bullet E^*)[1]$.

{\bf Courant algebroids.} A Courant algebroid structure on a pseudo-Euclidean vector bundle\footnote{We remind the reader that a~\textit{pseudo-Euclidean vector bundle} is a vector bundle $\mE\to\M$ endowed with a~symmetric, non-degenerate and $\fonc{\M}$-bilinear form on the space of sections of $\mE$, denoted $\dldot_\mE\colon \fmE\vee\fmE\to\foncm$ and referred to as the \textit{fiber-wise metric} } $\big({\mathcal E},\dldot_{\mathcal E}\big)$ is a pair $\big(\rho_{\mathcal E},\brdot_{\mE}\big)$ where
$\rho_\mE\colon \fmE\to\vf$ is a $\foncm$-linear map called the \textit{anchor} while $\brdot_\mE\colon \fmE\otimes\fmE\to\fmE$ is a $\corps$-bilinear form on the fibers of $E$ referred to as the \textit{Dorfman bracket}.

The latter satisfy the following minimal set of axioms:
\begin{enumerate}\itemsep=0pt
\item The Dorfman bracket satisfies the Jacobi identity in its Leibniz form
\begin{displaymath}
\br{e_1}{\br{{e_2}}{{e_3}}_\mE}_\mE=\br{\br{e_1}{{e_2}}_\mE}{{e_3}}_\mE+\br{{e_2}}{\br{e_1}{{e_3}}}_\mE \qquad \text{for all }e_1,{e_2},{e_3}\in\Gamma(\mE).
\end{displaymath}
\item The symmetric part of the Dorfman bracket is controlled by the anchor
\begin{displaymath}
\dl\br{e_1}{e_1}_\mE,{e_2}\dr_\mE=\tfrac{1}{2}\rho_\mE|_{e_2}\big[\dl e_1,e_1\dr_\mE\big] \qquad \text{for all }e_1,{e_2}\in\fmE.
\end{displaymath}
\item The fiber-wise metric is compatible with the Courant algebroid structure
\begin{displaymath}
\rho_\mE|_{e_1}\big[\dl {e_2},{e_3}\dr_\mE\big]=\dl\br{e_1}{{e_2}}_\mE,{e_3}\dr_\mE+\dl {e_2},\br{e_1}{{e_3}}_\mE\dr_\mE\qquad \text{for all }e_1,{e_2},{e_3}\in\fmE.
\end{displaymath}
\end{enumerate}
Letting $\big({\mathcal E},\dldot_{\mathcal E},\rho_{\mathcal E},\brdot_{\mE}\big)$ and $\big({\mathcal E'},\dldot_{\mathcal E'},\rho_{\mathcal E'},\brdot_{\mE'}\big)$ be two Courant algebroids over the same manifold\footnote{More general notions of morphisms of Courant algebroids over different base manifolds can be defined, see~\cite{Vysoky2019} for a thorough treatment.} $\M$, a morphism of Courant algebroids is defined as a morphism of the underlying vector bundles $\mathcal F\in\Hom(\mathcal E,\mathcal E')$ preserving the additional structures,\footnote{A morphism of vector bundles $\mathcal F\in\Hom(\mathcal E,\mathcal E')$ satisfying $\dlr{e_1}{e_2}_{\mE}=\dlr{\mathcal F(e_1)}{\mathcal F(e_2)}_{\mE'}$ will be called a morphism of pseudo-vector bundles.} i.e.,
\begin{gather*}
\dlr{e_1}{e_2}_{\mE}=\dlr{\mathcal F(e_1)}{\mathcal F(e_2)}_{\mE'} ,\qquad\rho_{\mE}=\rho_{\mE'}\circ\mathcal F,\\
 \mathcal F(\br{e_1}{e_2}_{\mE})=\br{\mathcal F(e_1)}{\mathcal F(e_2)}_{\mE'}\qquad \text{for all }e_1,{e_2}\in\Gamma(\mE).
\end{gather*}
Whenever the morphism $\mathcal F$ is invertible, it is called an isomorphism of Courant algebroids.

Courant algebroids first appeared implicitly in the work of I.~Dorfman \cite{Dorfman1987} and T.~Cou\-rant~\cite{Courant1986+,Courant1986} on integrable Dirac structures before their precise geometric structure was abstracted away in~\cite{Liu1995} to account for the concept of double of Lie bialgebroids.
More generally, the double ${\mathcal E}:=E\oplus E^*$ of a proto-Lie bialgebroid $(E,E^*)$ carries a natural structure of pseudo-Euclidean vector bundle with fiber-wise metric $\dl e_1,e_2\dr_{\mathcal E}:=\alpha(Y)+\beta(X)$, where $e_1:=X\oplus\alpha$ and $e_2:=Y\oplus\beta$ for all $X,Y\in\Gamma(E)$ and $\alpha,\beta\in\Gamma(E^*)$. One can furthermore endow ${\mathcal E}$ with the anchor $\rho_{\mathcal E}(e_1):=\rho(X)\oplus R(\alpha)$ while the Dorfman bracket is defined through the following explicit expression
\begin{gather*}
\br{e_1}{e_2}_{\mathcal E}:=\big(\br{X}{Y}_E+\Lag^{E^*}_{\alpha}Y-\iota_\beta {\rm d}_{E^*}X-\varphi(\alpha,\beta,\cdot)\big)\\
\hphantom{\br{e_1}{e_2}_{\mathcal E}:=}{}
\oplus\big(\br{\alpha}{\beta}_{E^*}+\Lag^E_X\beta-\iota_Yd_E\alpha-\psi(X,Y,\cdot)\big),
\end{gather*}
where $\Lag^E_X$ stands for the unique derivative operator extending the action of $\br{X}{\cdot}_E$ to the tensor algebra of $E$ \big(and similarly for $\Lag^{E^*}_{\alpha}$\big). The axioms of a proto-Lie bialgebroid thus ensure that the pair $\big(\rho_{\mathcal E},\brdot_{\mE}\big)$ defines a Courant algebroid structure on $\big({\mathcal E},\dldot_{\mathcal E}\big)$, cf.\ \cite[Theorem~3.8.2]{Roytenberg2007}.

\begin{Example}[exact Courant algebroids]\label{example:Exact Courant algebroids}
A Courant algebroid such that the following sequence\footnote{The vector bundle morphism $\rho_{\mE}^*\colon T^*\M\to \mE$ is defined via $\dlr{\rho^*_\mE|_\alpha}{e}_\mE=\dlr{\alpha}{\rho_\mE|_e}$ for all $\alpha\in\ff$ and $e\in\Gamma(\mE)$, and satisfies $\rho_{\mE}\circ\rho_{\mE}^*=0$.}
\begin{displaymath}
0\longrightarrow T^*\M \overset{\rho_{\mE}^*}{\longrightarrow} \mE \overset{\rho_{\mE}}{\longrightarrow} T\M \longrightarrow0
\end{displaymath}
is exact is called an exact Courant algebroid. It is an important result due to P.~{{\v{S}}evera} \cite{Severa2017} that exact Courant algebroids over a fixed base manifold $\M$ are classified\footnote{More precisely, any exact Courant algebroid is isomorphic to $\mE:=T\M\oplus T^*\M$ with canonical fiber-wise metric, projection $T\M\oplus T^*\M\to T\M$ as anchor and Dorfman bracket $\br{e_1}{e_2}_{\mathcal E}:=\br{X}{Y}\oplus\big(\Lag_X\beta-\iota_Yd\alpha-H(X,Y,\cdot)\big)$, where $H\in\Om^3(\M)$ is a closed 3-form. The latter Courant algebroid structure coincides with the one induced by the quasi-bialgebroid structure of Example~\ref{example:TWP}. Letting $\om\in\Gamma\big({\w}^2E^*\big)$, the isomorphism $X\oplus \alpha\overset{\om}{\longrightarrow} X\oplus\big(\alpha+\half \iota_X\om\big)$ corresponding to the twist \eqref{canonicaltwistom} amounts to shift the closed 3-form $H$ by an exact 3-form as $H\overset{\om}{\longrightarrow}H-d\om$. Hence isomorphism classes of exact Courant algebroids over $\M$ are in bijective correspondence with elements of~$H^3_\textsf{dR}(\M)$.} by the third de Rham cohomology $H^3_\textsf{dR}(\M)$ of $\M$.
\end{Example}

{\bf Dirac structures.}
A subbundle $L\subset \mathcal E$ of a Courant algebroid $\mathcal E$ is called a {Dirac structure}~if:
\begin{enumerate}\itemsep=0pt
\item $L$ is maximally isotropic with respect to the fiber-wise metric $\dldot_{\mathcal E}$, i.e., $\dlr{\Gamma(L)}{\Gamma(L)}_\mE=0$, $\dim L=\half \dim \mathcal E$.
\item The space of sections $\Gamma(L)$ is closed under the Courant bracket $\brdot_{\mathcal E}$, i.e., $\br{\Gamma(L)}{\Gamma(L)}_{\mathcal E}\subseteq \Gamma(L)$.
\end{enumerate}
These conditions ensure in particular that the restrictions of the Courant anchor and bracket to~$L$ endow the vector bundle~$L$ with a structure of Lie algebroid.
\begin{Example}[Dirac structure]
Let $\M$ be a manifold and $[H]\in H^3_\textsf{dR}(\M)$. Given a representative $H\in\Om^3(\M)$ and a twisted Poisson structure~$\pi$ with respect to~$H$, the graph of the map $\pi^\sharp\colon T^*\M\to T\M$ can be checked to be a Dirac structure for the exact Courant algebroid associated to $\big(\M,[H]\big)$ \cite{Severa2017,Severa2001}.
\end{Example}

A pair $(\mathcal E,L)$ where $L$ is a Dirac structure for the Courant algebroid $\mathcal E$ is referred to as a~\textit{Manin pair}. A triplet $(\mathcal E,L,M)$ where $L$, $M$ are two Dirac structures of $\mathcal E$ such that $\mathcal E=L\oplus M$ is referred to as a~\textit{Manin triple}. Both the notions of Manin pairs and triples recover their usual extension in the context of Lie algebras when the base manifold $\M$ is a point. Letting $\big({\mathcal E},L\big)$ and $\big({\mathcal E'},L'\big)$ be two Manin pairs over the same manifold $\M$, a morphism of Manin pairs\footnote{Similarly, letting $\big({\mathcal E},L,M\big)$ and $\big({\mathcal E'},L' M'\big)$ be two Manin triples over the same manifold $\M$, a morphism of Manin triples is a~morphism of Courant algebroids $\mathcal F\in\Hom(\mathcal E,\mathcal E')$ such that $\mathcal F(L)\subseteq L'$ and $\mathcal F(M)\subseteq M'$.} is a~morphism of Courant algebroids $\mathcal F\in\Hom(\mathcal E,\mathcal E')$ such that $\mathcal F(L)\subseteq L'$.

 Letting $\mathcal E=E\oplus E^*$ be the pseudo-Euclidean vector bundle associated to the vector bundle $E$ over $\M$ and endowed with the canonical fiber-wise metric $\dl \cdot,\cdot\dr_{\mathcal E}$, we have the following hierarchy of identifications:
\begin{itemize}\itemsep=0pt
\item proto-Lie bialgebroids on $(E,E^*)$ $\Leftrightarrow$ Courant algebroids structures on $\big(\mathcal E,\dl \cdot,\cdot\dr_{\mathcal E}\big)$,
\item Lie-quasi bialgebroids on $(E,E^*)$ $\Leftrightarrow$ Manin pairs $\big((\mathcal E,\dl \cdot,\cdot\dr_{\mathcal E}), E\big)$,
\item quasi-Lie bialgebroids on $(E,E^*)$ $\Leftrightarrow$ Manin pairs $\big((\mathcal E,\dl \cdot,\cdot\dr_{\mathcal E}), E^*\big)$,
\item Lie bialgebroids on $(E,E^*)$ $\Leftrightarrow$ Manin triples $\big((\mathcal E,\dl \cdot,\cdot\dr_{\mathcal E}), E,E^*\big)$.
\end{itemize}
Such identification allows to define morphisms of proto/Lie-quasi/quasi-Lie/Lie bialgebroids as morphisms of Courant algebroids preserving additional Dirac structures.
Explicitly, morphisms of pseudo-vector bundles $E\oplus E^*\overset{\mathcal F}{\longrightarrow} E'\oplus E'{}^*$ with canonical fiber-wise metrics generically take the form
\begin{gather*}
\mathcal F(e)^{a'}=\mathcal F^{a'}{}_b X^b+\mathcal F^{a'b}\alpha_b ,\qquad \mathcal F(e)_{a'}=\mathcal F_{a'b}X^b+\mathcal F_{a'}{}^b\alpha_b
\end{gather*}
satisfying
\begin{gather*}
\mathcal F_{c'(a}\mathcal F^{c'}{}_{b)}=0 ,\qquad \mathcal F^{c'(a}\mathcal F_{c'}{}^{b)}=0 ,\qquad \mathcal F_{c'}{}^a \mathcal F^{c'}{}_b+\mathcal F^{c'a}\mathcal F_{c'b}=\delta^a{}_b.
\end{gather*}
Such morphisms map the Dirac structure $E$ to $E'$ if and only if $\mathcal F_{a'b}=0$ and the Dirac structure~$E^*$ to~$E'{}^*$ if and only if $\mathcal F^{a'b}=0$.
Infinitesimal endomorphisms of the pseudo-Euclidean vector bundle $(E\oplus E^*,\dl \cdot,\cdot\dr_{\mathcal E})$ read\footnote{Restricting to Lie-quasi (resp. quasi-Lie) bialgebroids yields $\om\equiv0$ (resp.~$\Lambda\equiv0$), so that the subalgebra $\End\big({\fE}\big)$ acts on the abelian ideal $\Gamma\big({\w}^2E\big)$ \big(resp.~$\Gamma\big({\w}^2E^*\big)$\big) by rotations.}
\begin{gather*}
\delta\mathcal F(e)=\big(\lambda(X)+\tfrac{1}{2}\iota_\alpha\Lambda\big)\oplus\big({-}\lambda^T(\alpha)+\tfrac{1}{2}\iota_X\om\big),
\end{gather*}
where
\begin{itemize}\itemsep=0pt
\item $\lambda\in\End({\fE})$ generates infinitesimal rotations of the fibers of $(E,E^*)$,
\item $\Lambda\in\Gamma\big({\w}^2E\big)$ generates the infinitesimal version of the twist of Lie-quasi bialgebroids~\eqref{twistLambdaexp},
\item $\om\in\Gamma\big({\w}^2E^*\big)$ generates the infinitesimal version of the twist of quasi-Lie bialgebroids~\eqref{twistomexp}.
\end{itemize}

\section[Incarnation of the Theta-graph in d=3]{Incarnation of the $\boldsymbol{\Theta}$-graph in $\boldsymbol{d=3}$}\label{section:Incarnation}

\begin{figure}[ht]\centering\vspace*{-5mm}
 $\raisebox{-4ex}{\hbox{\begin{tikzpicture}[scale=0.5, >=stealth']
\tikzstyle{w}=[circle, draw, minimum size=4, inner sep=1]
\tikzstyle{b}=[circle, draw, fill, minimum size=2, inner sep=0.02]
\node [b] (b2) at (-30:3) {2};
\node [b] (b1) at (-90:1.5) {1};
\node [b] (b3) at (-90:3) {3};
\node [b] (b4) at (0,0) {4};
\node [b] (b6) at (-150:3) {6};
\draw[black,>=latex] (b3) to (b6);
\draw[black,>=latex] (b4) to (b2);
\draw[black,>=latex] (b4) to (b6);
\draw[black,>=latex] (b4) to (b1);
\draw[black,>=latex] (b2) to (b1);
\draw[black,>=latex] (b3) to (b2);
\end{tikzpicture}}}_A$
\hspace{2cm}
 $\raisebox{-4ex}{\hbox{\begin{tikzpicture}[scale=0.5, >=stealth']
\tikzstyle{w}=[circle, draw, minimum size=4, inner sep=1]
\tikzstyle{b}=[circle, draw, fill, minimum size=2, inner sep=0.02]
\node [b] (b2) at (-30:3) {2};
\node [b] (b1) at (-90:1.5) {1};
\node [b] (b3) at (-90:3) {3};
\node [b] (b4) at (0,0) {4};
\node [b] (b6) at (-150:3) {6};
\draw[black,>=latex] (b3) to (b6);
\draw[black,>=latex] (b4) to (b2);
\draw[black,>=latex] (b4) to (b6);
\draw[black,>=latex] (b4) to (b1);
\draw[black,>=latex] (b3) to (b1);
\draw[black,>=latex] (b3) to (b2);
\end{tikzpicture}}}_B$
\hspace{2cm}
 $\raisebox{-5ex}{\hbox{\begin{tikzpicture}[rotate=90,scale=0.5, >=stealth']
\tikzstyle{w}=[circle, draw, minimum size=4, inner sep=1]
\tikzstyle{b}=[circle, draw, fill, minimum size=2, inner sep=0.02]
\node [b] (b2) at (30:2) {2};
\node [b] (b1) at (150:2) {1};
\node [b] (b3) at (-90:2) {3};
\node [b] (b4) at (90:2) {4};
\node [b] (b5) at (-30:2) {5};
\node [b] (b6) at (-150:2) {6};
\draw[black,>=latex] (b4) to (b1);
\draw[black,>=latex] (b4) to (b2);
\draw[black,>=latex] (b3) to (b4);
\draw[black,>=latex] (b3) to (b6);
\draw[black,>=latex] (b6) to (b1);
\draw[black,>=latex] (b5) to (b2);
\draw[black,>=latex] (b3) to (b5);
\end{tikzpicture}}}_C$

 \caption{Shape of graphs involved in $\vartheta_3$ ($A$ and $B$) and $\Theta_3$ ($C$).} \label{figshape}
\end{figure}
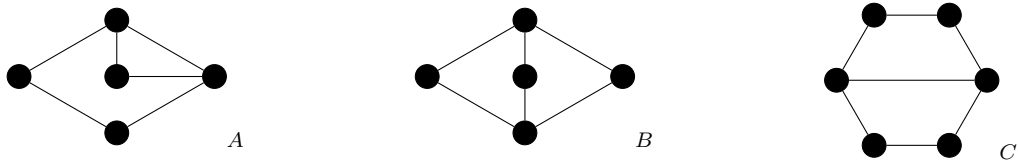

The present appendix is devoted to collect some additional results regarding the exotic \Lieinf-structure $\theta$ of Theorem~\ref{TheoremMain} generated by the cocycle class $[\Theta_3]\in H^1({\odfGC{3}{0}{2}})$. For concreteness, we fix a representative of the class $[\Theta_3]$ as follows:

\begin{Proposition}\label{propTHETA3}
There is a unique pair of combination of graphs $\Theta_3\in\odfGC{3}{0}{2}$ and $\vartheta_3\in\odfGCblack{3}{1}{1}$ such that:
\begin{enumerate}\itemsep=0pt
\item[$1.$] $\Theta_3=-\delta\vartheta_3$, i.e., $\Theta_3$ is exact in $\odfGCblack{3}{1}{1}$.
\item[$2.$] $\Theta_3$ contains only graphs of the shape $C$, cf.\ Figure~{\rm \ref{figshape}}.
\item[$3.$] Each graph of~$\vartheta_3$ contains at least one red cycle, i.e., $\vartheta_3\notin\odfGC{3}{0}{2}$.
\end{enumerate}
\end{Proposition}

\begin{Remark}\label{remApp}\quad
\begin{itemize}\itemsep=0pt
\item The combination of graphs $\vartheta_3$ contains 68 black-oriented graphs (48 graphs of shape $A$ and 20 graphs of shape $B$) while the combination $\Theta_3$ contains 288 bi-oriented graphs of shape $C$.
\item Although $\Theta_3$ is exact in $\odfGCblack{3}{1}{1}$, it is crucial to note that $\Theta_3$ is {\it not} exact in $\odfGC{3}{0}{2}$, i.e., there is no combination of graphs $\kappa_3\in\odfGC{3}{0}{2}$ such that $\Theta_3=-\delta\kappa_3$. Hence $\Theta_3$ is a non-trivial cocycle in~$\odfGC{3}{0}{2}$.
\end{itemize}
\end{Remark}

Let $E\overset{\pi}{\to} \M$ be a vector bundle. To each Lie bialgebroid structure on $(E,E^*)$ (represented by the Hamiltonian function $\cH\in \A^E_{\rm Lie}$), we will associate the functions
\begin{itemize}\itemsep=0pt
\item $\vartheta_3(\cH):=\Rep^E_5(\vartheta_3)(\cH^{\otimes 5})\in \A^E_\text{Lie-quasi}$,
\item $\Theta_3(\cH):=\Rep^E_6(\Theta_3)(\cH^{\otimes 6})\in \A^E_{\rm Lie}$.
\end{itemize}

Note that the condition $\Theta_3=-\delta\vartheta_3$ ensures that $\Theta_3(\cH)\sim\Q[\vartheta_3(\cH)]$~-- where the differential~$\Q$ is defined as $\Q:=\pb{\cH}{\cdot}_\symp^E$~-- so that $\Theta_3(\cH)$ is a coboundary in the complex $(\A^E_\text{Lie-quasi}|\Q)$. However, $\Theta_3(\cH)$ is generically not exact in $\A^E_{\rm Lie}$ and the obstruction for $\Theta_3(\cH)$ to be a coboundary in $\A^E_{\rm Lie}$ is precisely given by the component of $\vartheta_3(\cH)$ proportional to $\zeta^3$. We will denote this obstruction as $\Ob(\cH)^{abc}$, so that 
\[ 
\Ob(\cH)^{abc}=0\ \Rightarrow \ \vartheta_3(\cH)\in\A^E_{\rm Lie}
\] and $\Theta_3(\cH)$ is a trivial cocycle in $(\A^E_\text{Lie-quasi}|\Q)$.
A straightforward computation gives
\begin{gather*}
\Ob(\cH)^{abc}= R^{d| \mu} R^{e| \nu} \big(2 \rho_{f}{}^{\lambda} \p_{\mu}{C_{e}{}^{[a| f}} \p_{\lambda \nu}{C_{d}{}^{b c]}}-\rho_{d}{}^{\beta} \p_{\beta \nu}{R^{[a| \lambda}} \p_{\mu \lambda}{C_{e}{}^{b c]}}-\rho_{d}{}^{\lambda} \p_{\mu}{C_{e}{}^{[a| f}} \p_{\lambda \nu}{C_{f}{}^{b c]}}\\
\hphantom{\Ob(\cH)^{abc}=}{}
-\p_{\mu}{R^{f| \lambda}} f_{e f}{}^{[a} \p_{\nu \lambda}{C_{d}{}^{b c]}}
-2 f_{d f}{}^{[a} \p_{\mu}{C_{e}{}^{b |g}} \p_{\nu}{C_{g}{}^{c] f}}+2 f_{f g}{}^{[a} \p_{\mu}{C_{e}{}^{b |f}} \p_{\nu}{C_{d}{}^{c] g}}\big)\\
\hphantom{\Ob(\cH)^{abc}=}{}
+R^{d| \mu} \p_{\mu}{R^{e| \nu}} \big(f_{d f}{}^{[a} C_{e}{}^{f g} \p_{\nu}{C_{g}{}^{b c]}}-2 f_{e f}{}^{[a} C_{d}{}^{f g} \p_{\nu}{C_{g}{}^{b c]}}
+ \rho_{d}{}^{\lambda} \p_{\lambda}{C_{e}{}^{[a| f}} \p_{\nu}{C_{f}{}^{b c]}}
\\
\hphantom{\Ob(\cH)^{abc}=}{}
-2 \rho_{e}{}^{\lambda} \p_{\lambda}{C_{d}{}^{[a| f}} \p_{\nu}{C_{f}{}^{b c]}}-\rho_{e}{}^{\lambda} \p_{\lambda}{C_{f}{}^{[a b}} \p_{\nu}{C_{d}{}^{c] f}}\big)\\
\hphantom{\Ob(\cH)^{abc}=}{}
+\rho_{d}{}^{\beta} \p_{\beta}{R^{e |\nu}} R^{d| \mu} \p_{\mu \nu}{R^{[a| \lambda}} \p_{\lambda}{C_{e}{}^{b c]}}
+2 \rho_{e}{}^{\lambda} R^{d| \mu} \p_{\lambda \mu}{R^{[a| \nu}} C_{d}{}^{e f} \p_{\nu}{C_{f}{}^{b c]}}\\
\hphantom{\Ob(\cH)^{abc}=}{}
+\rho_{g}{}^{\nu} R^{d| \mu} C_{d}{}^{e f} \p_{\nu}{C_{e}{}^{[a b}} \p_{\mu}{C_{f}{}^{c] g}}.
\end{gather*}

The latter encodes the first obstruction for $\cH$ to define a ``quantizable Lie bialgebroid''. Although the obstruction does not vanish for a generic Lie bialgebroid, the following proposition displays two important examples:
\begin{Proposition}\label{propob}
The obstruction vanishes for
\begin{itemize}\itemsep=0pt
\item Lie bialgebras,
\item coboundary Lie bialgebroids.
\end{itemize}
\end{Proposition}
\begin{proof}
Setting $R^{a| \mu}\equiv0$ yields $\Ob(\cH)^{abc}=0$ hence the obstruction vanishes for Lie bialgebras. More generally, it can be checked that each graph appearing in the combinations $\vartheta_3$ and $\Theta_3$ contains at least one arrow of the type ${\raisebox{-0.5ex}{\hbox{\begin{tikzpicture}[scale=0.5, >=stealth']
\tikzstyle{w}=[circle, draw, minimum size=4, inner sep=1]
\tikzstyle{b}=[circle, draw, fill, minimum size=2, inner sep=0.02]
\draw (0,0) node [ext] (b1) {\rm \tiny{$i$}};
\draw (1.3,0) node [exttiny] (b2) {\rm \tiny{$j$}};
\draw(0.6,0.3) node {};
\draw[black,->,>=latex,line width=0.2mm][decoration={markings,mark=at position .55 with {\arrow[scale=1,linkcolor]{latex}}},
 postaction={decorate},
 shorten >=0.4pt] (b1) to (b2);
\end{tikzpicture}}}}$ so that both $\vartheta_3(\cH)$ and $\Theta_3(\cH)$ vanish identically on the graded Poisson subalgebra\footnote{See footnote~\ref{footpoint}.} $\A^\alg_{\rm Lie}\subset\A^E_{\rm Lie}$.

For coboundary Lie bialgebroids, we perform the replacement $R^{a|\mu}=\rho_b{}^\mu\Lambda^{b a}$ and $C_c{}^{ab}=-\rho_c{}^\mu\p_\mu\Lambda^{ab}-2\Lambda^{d[a}f_{dc}{}^{b]}$, with $\Lambda^{ab}=\Lambda^{[ab]}$ a bivector, see Appendix \ref{section:Appendix}. Under this replacement, it can be checked that $\Ob(\cH)^{abc}$ identically vanishes modulo the defining conditions $\mathcal C_1\equiv0$, $\mathcal C_2\equiv0$ ensuring that the maps $(\rho,f)$ define a Lie algebroid.
\end{proof}

\subsection*{Acknowledgements}

The author is grateful to T.~Basile, D.~Lejay, H.~Y.~Liao and P.~Xu for discussions as well as to Y.~Kosmann--Schwarzbach, S.~Merkulov and T.~Willwacher for correspondence. The author would also like to thank J.~H.~Park for support. Finally, the author is grateful to anonymous referees whose suggestions greatly helped to improve the quality of the present paper. This work was supported by Brain Pool Program through the National Research Foundation of Korea (NRF) funded by the Ministry of Science and ICT (2018H1D3A1A01030137) and by Basic Science Research Program through the National Research Foundation of Korea (NRF) funded by the Ministry of Education (NRF-2020R1A6A1A03047877).

\pdfbookmark[1]{References}{ref}
\LastPageEnding

\end{document}